\documentclass[journal, final]{IEEEtran}
\usepackage{setspace}

\usepackage[table]{xcolor}
\usepackage{graphicx}
\usepackage{grffile}
\usepackage{tikz}
\usepackage[noadjust]{cite}

\usetikzlibrary{positioning}
\definecolor{highlight}{HTML}{cdeeff} 
\definecolor{Gray}{gray}{0.9}
\usepackage{amsmath, amsthm, amssymb,bm} 

\usepackage{color}
\usepackage{framed}
\usepackage[short]{optidef} 
\usepackage{hyperref,hypcap}\hypersetup{colorlinks=true,allcolors=blue}
\usepackage[capitalise]{cleveref}

\usepackage{algorithm}\usepackage[noend]{algpseudocode} 

\usepackage[font=footnotesize]{caption} 
\usepackage[font=footnotesize]{subcaption}

\graphicspath{{./figures/}}

\usepackage{color}

 \newcommand{\removed}[1]{}

\newcommand{\R}{\mathbb{R}} 

\newcommand{\cV}{{\mathcal{V}}} 

\newcommand{\ds}{{\displaystyle}}

\renewcommand{\theequation}{\thesection.\arabic{equation}}

\newtheorem{pro}{Proposition}[section]

\newcommand{\mv}[1]{\boldsymbol{\mathbf{#1}}} 

\DeclareMathOperator*{\argmin}{arg\,min}

\begin{document}

\title{Simultaneous Drift Calibration and Reconstruction for Scanning-Probe X-Ray Tomography}


\author{Xiang~Huang,~Parth~Brahmbhatt,~Chang~Meng,~and~Zichao~Wendy~Di%
\thanks{X.~Huang and Z.~W.~Di are with the Mathematics and Computer Science Division, Argonne National Laboratory, Lemont, IL 60439 USA.}%
\thanks{Z.~W.~Di is also with the Advanced Photon Source, Argonne National Laboratory, Lemont, IL 60439 USA.}%
\thanks{P.~Brahmbhatt is with the Department of Chemical and Biological Engineering, University of Wisconsin--Madison, Madison, WI 53706 USA.}%
\thanks{C.~Meng is with the Department of Mathematics, Emory University, Atlanta, GA 30322 USA.}}

\maketitle

\begin{abstract}
  Scanning-probe x-ray tomography is useful for imaging nanoscale structures of a sample. The image resolution,  which can reach down to 10 nm by the increased brightness and coherence of the x-ray optics, however, is highly susceptible to experimental error. Failure to
  address these errors can lead to a smeared image and, in the worst case, to misinterpretation of the imaged object's structure. In this work, we
  present a novel optimization-based approach to calibrate a common yet challenging source of experimental error, the drifts of the scanning positions,
  while simultaneously reconstructing the object. This approach uses the coupled and complementary information from different
  measurements to enforce consistency between the measurements and the reconstruction. We illustrate the proposed approach on
  both synthetic and real tomography images and show its superior performance compared with reconstruction without explicit error
  calibration.
\end{abstract}


\section{Introduction}
\label{sec:intro}

X-ray tomography is a powerful technique for revealing the internal three-dimensional structure of materials and biological specimens, and it plays a critical role in fields ranging from materials science to environmental and life sciences~\cite{baruchel2000x,paulus2000high}. By collecting a series of x-ray projections at multiple viewing angles, tomographic reconstruction enables quantitative visualization of the internal morphology of a sample. To achieve nanoscale resolution, x-ray tomography is often implemented in a \emph{scanning-probe} geometry, in which a tightly focused x-ray beam is raster-scanned across the specimen while the object or the beam is rotated~\cite{de2010hard,dierolf2010ptychographic}. Recent advances in x-ray optics and synchrotron instrumentation have pushed the spatial resolution of scanning-probe x-ray tomography to 10~nm and below~\cite{yan2018multimodal}. At this scale, however, even small mechanical or thermal instabilities can lead to significant geometric inconsistencies that degrade reconstruction quality.

\begin{figure}[htb]
  \begin{minipage}{0.15\textwidth}
    \begin{tikzpicture}
        \node (img)  {\includegraphics[width=\textwidth]{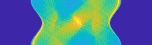}};
        \node[below=of img, node distance=0cm, yshift=1.4cm,font=\color{black}] {$\tau$};
        \node[left=of img, node distance=0cm, rotate=90, anchor=center,yshift=-0.9cm,font=\color{black}] {$\theta$};
    \end{tikzpicture}
  \end{minipage}%
  \quad
  \begin{minipage}{0.15\textwidth}
    \begin{tikzpicture}
        \node (img)  {\includegraphics[width=\textwidth]{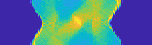}};
    \end{tikzpicture}
  \end{minipage}%
  \begin{minipage}{0.15\textwidth}
    \begin{tikzpicture}
        \node (img)  {\includegraphics[width=\textwidth]{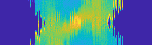}};
    \end{tikzpicture}
  \end{minipage}
  \caption{Comparison of sinograms showing the effects of scanning-position drift.
  Left: drift-free sinogram; middle: sinogram with a maximum drift of one expected spatial resolution; right: sinogram with a maximum drift of five times the expected spatial resolution.}
  \label{fig:sinograms}
\end{figure}

In practice, a variety of experimental factors---including thermal expansion, sample drift, vibration, and rotation-stage misalignment---introduce errors in the assumed scanning positions. For instance, H{\"u}e~\emph{et~al.}~\cite{hue2011extended,hue2010wave} observed position errors of the same order of magnitude as the target spatial resolution over the duration of a single acquisition. As illustrated in Fig.~\ref{fig:sinograms}, even modest position drifts can visibly distort the sinogram and, consequently, blur or smear the reconstructed image. Because such instabilities are intrinsic to long-duration, high-resolution experiments, numerical correction of these geometric errors is essential.

When the drift magnitude is small relative to the beam width, it can sometimes be mitigated by incorporating prior knowledge through Bayesian or variational regularization frameworks~\cite{van2014super,charbonnier1997deterministic,sauer1994bayesian}. However, as imaging resolution improves, these generic treatments become insufficient. In full-field tomography, research has focused primarily on correcting the dominant geometric error---the drift of the center of rotation---through feature-tracking or optimization-based approaches~\cite{sisc2018,cao2011automatic,gursoy2017rapid,ali2022accelerating}. In scanning-probe modalities, a number of techniques have been developed to explicitly address scanning-position errors, including simulated annealing, model-based correction, and cross-correlation of overlapping or shadow images~\cite{beckers2013drift,tripathi2014ptychographic,zhang2013translation,hurst2010probe,sang2014revolving}. While effective in specialized contexts, these methods are often model-specific, computationally intensive, and difficult to generalize for automated tomographic reconstruction.

In this work, we build upon our earlier optimization-based framework~\cite{huang2019calibrating,di2019optimization} and present an extended, more general approach that simultaneously calibrates scanning-position drift and reconstructs the object. The key idea is to exploit the coupling and complementary information among multiple projection measurements to enforce a coherent relationship between the reconstruction and the observed data. We explicitly model the scanning drift within the tomographic forward operator and solve a joint inverse problem over both the object and the drift parameters. This formulation preserves the flexibility to incorporate prior knowledge, yet enables greater automation and robustness than previous approaches. Moreover, it can handle both systematic (Type~I) and spatially varying (Type~II) drift patterns, providing reliable performance even under limited or noisy data.

The remainder of this paper is organized as follows. In Sec.~\ref{sec:model}, we introduce the mathematical formulation of the drifted tomographic model. Section~\ref{sec:algorithm} describes the alternating optimization algorithm for simultaneous reconstruction and drift calibration. Numerical and experimental validations are presented in Sec.~\ref{sec:simulation}, followed by concluding remarks and future directions in Sec.~\ref{sec:conclusion}. For brevity, we focus on two-dimensional reconstruction, while extension to the three-dimensional parallel-beam case is straightforward.

\section{Problem Formulation}\label{sec:model}

Given a 2D object, we discretize it as $N\times N$ pixels whose indices are represented in column-major-order as an array $\cV$, and denote by $\mv{w}=\{\mv{w}_v:v\in\cV\}\in \R^{N^2 \times 1}$, where
$\mv{w}_v$ denotes the pixel intensity at $v\in\cV$,  that is, the value of the object property we intend to recover such as attenuation coefficient for x-ray
transmission and refractive index in ultrasound tomography \cite{devaney1982filtered}.
In the ideal case without drifts of the scanning positions, given a fixed step size $\Delta$, the object is rotated with $N_\theta,\theta \in \{1,\dots,N_\theta\}$ angles, where at each angle, it is raster scanned with $N_\tau,\tau \in \{1,\dots,N_\tau\}$ steps. Let
$L_v^{\theta,\tau}$ be the intersection length of the beam $(\theta,\tau)$ with the pixel $v\in\cV$.
Then the forward model of tomography is governed by the Radon transform \cite{kak_principles_1988} which is formulated as
\[R(\mv{w})=\mv{L}\mv{w}.\]
Let $\mv{S} \in \R^{N_\theta \times N_\tau}$ denote the 2D measurement data (i.e., the sinogram) and $\mv{s}\in \R^{N_\theta N_\tau \times 1}$ denote the vector formed by stacking the columns of $\mv{S}$. We obtain the traditional tomographic reconstruction of $\mv{w}$ by solving the following
optimization problem:
\begin{equation}
  \label{eq:P0}
  \ds
  \min_{\mv{w} \geq \mv{0}}\frac{1}{2}\|\mv{L}\mv{w}-\mv{s}\|^2+ \lambda \Phi (\mv{w}),
\end{equation}
where $\Phi (\mv{w})$ is a regularizer (e.g., total variation (TV) \cite{engl1996regularization}) to incorporate prior knowledge, and
$\lambda\geq 0$ is a scalar balancing the data misfit (e.g., noise, model error) and 
regularizer. Notice we also have the nonnegative constraint on $\mv{w}$ based on the fact that these physical properties can not be negative. 

\begin{figure}[h]
  \centering
  \includegraphics[width=0.35\textwidth]{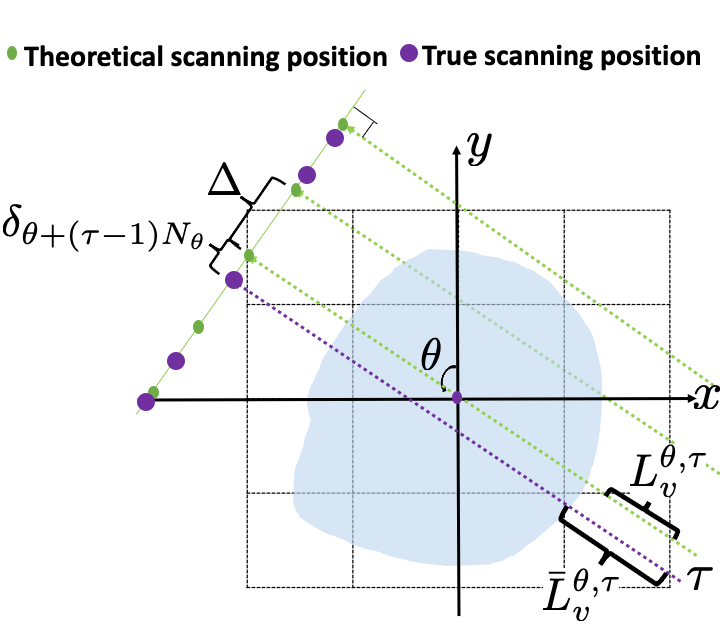}
  \caption{Geometric sketch of the drifts of the scanning positions in tomography.}
  \label{fig:geometry}
\end{figure}

However, in the presence of the drifts of the scanning positions, the true forward model $\mv{L}$ is no longer static, instead, we can write it as a function of the drifts $\mv{L}(\mv{\delta})$, where $\mv{\delta}=[\delta_{(\theta+(\tau-1)N_\theta}] \in \R^{N_\theta N_\tau \times 1}$ denotes the collection of drifts for each measurement. Equivalently, we extend the Radon transform as a function of both $\mv{W}$ and $\mv{\delta}$: \[R(\mv{w},\mv{\delta})=\mv{L}(\mv{\delta})\mv{w}.\]
Fig.~\ref{fig:geometry} illustrates the geometric sketch of the experiment with the drifts of the scanning positions. Beam $\tau$ without drift
(which is referred as the theoretical scanning position) is shown as the green line, given the constant scanning step $\Delta$. Correspondingly, its drifted
position with drift $\delta_{\theta+(\tau-1)N_\theta}$ is shown as the purple line. In practice, drift is a consequence of both systematic error and random error. One can expect a manual calibration before each rotation during the course
of the experiment, so that the systematic drifts of the scanning positions are the same after any rotation. However, random drifts are expected to be different for each individual scanning. 

\section{Optimization-Based Reconstruction Algorithm}\label{sec:algorithm}
 Our goal is to reconstruct object $\mv{w}$ from the measurement $\mv{s}$, while simultaneously calibrating unknown drift $\mv{\delta}$. It is well known that traditional tomographic reconstruction is ill-posed due to the limited number of measurements, the added degree of freedom in the current setting with respect to drift only makes the reconstruction even more ill-posed, therefore, we propose to solve the following regularized optimization problem:
\begin{equation}
  \label{eq:optwd}
\min_{\mv{w}, \mv{\delta}}\frac{1}{2}\|R(\mv{w}, \mv{\delta}) - \mv{s} \|^2+ \lambda \Phi_w (\mv{w}) + \mu \Phi_d (\mv{\delta}),
\end{equation}
where the regularization functions $\Phi_w (\mv{w})$, $\Phi_d
(\mv{\delta})$ allow flexible incorporation of prior knowledge for both the object and the drifts, respectively, and $\lambda, \mu \geq 0$ are the corresponding weights. 

We propose to optimize \cref{eq:optwd},  by alternating between $\mv{\delta}$ and $\mv{w}$. In practice, we often have a rough idea on how big the drift can be, i.e., have an estimate of maximum drift $\Delta_{\max}$. Therefore it is reasonable to restrict the drifts to fall in $[-\Delta_{\max}, \Delta_{\max}]$.
\cref{alg:driftedTomography} details the proposed alternating minimization. In step 4, we propose an adaptive way to update the regularizer weight $\lambda$ to better compensate for the model/data discrepancy (details in \cref{sec:autoLambda}) and accelerate the convergence: generally it uses a large $\lambda$ in the beginning and gradually reduce it in later iterations. In step 5, \cref{eq:optw} is a classic regularized linear inverse problem, and we explore iterative shrinkage-thresholding algorithms such as
TwIST\cite{bioucas2007new} to efficiently solve it.
At each iteration of TwIST,  we enforce the nonnegative constraint by projecting the negative components of $\mv{w}$ to 0, which yields faster convergence and better reconstruction. To solve \cref{eq:optd} in step 6, we consider two different scenarios, as detailed in  \cref{subsec:SpecialDrift} and \cref{subsec:GeneralDrift} respectively. 

\begin{algorithm}[H]
\caption{Tomographic Reconstruction with Drift Calibration} \label{alg:driftedTomography}
\begin{algorithmic}[1]
\State Input: $\mv{s}, \mu,\mv{\delta}^0, \mv{w}^0 , \lambda_0$
\State $k=1$
\While{not converged}
\State  Update $\lambda_k$ using \cref{eq:lambdaupdate}
\State  Use TwIST to solve
\begin{equation}
  \label{eq:optw}
\ds \mv{w}^k=\argmin_{\mv{w}\geq 0}\frac{1}{2}\|\mv{L}(\mv{\delta}^{k-1})\mv{w}-\mv{s} \|^2+ \lambda_k \Phi_w (\mv{w}),
\end{equation}
\State Solve
\begin{equation}
  \label{eq:optd}
  \mv{\delta}^k=\argmin_{-\Delta_{\max}\leq \mv{\delta}\leq \Delta_{\max}}\frac{1}{2}\|R(\mv{\delta}; \mv{w}^k)-\mv{s} \|^2+ \mu \Phi_d (\mv{\delta}),
\end{equation}
\State $k \gets k+1$
\EndWhile

\State Output: $\mv{w}^k, \mv{\delta}^k $
\end{algorithmic}
\end{algorithm}

The computational cost of \cref{alg:driftedTomography} is dominated by the repeated TwIST solves in the image and drift updates. If $K$ outer iterations are performed and $T_k^w$ and $T_k^\delta$ denote the corresponding inner TwIST iteration counts at outer iteration $k$, then the total cost scales as $\sum_{k=1}^{K}(T_k^w C_w + T_k^\delta C_\delta)$, where $C_w$ and $C_\delta$ denote the cost of one forward/adjoint evaluation in the image and drift subproblems, respectively. In the numerical experiments, we use at most 10 outer iterations.

The iterative reconstruction process is terminated when the absolute change in the data fidelity term between consecutive iterations drops below a tolerance of $10^{-3}$. This stopping criterion is defined as $|\mathcal{J}_k - \mathcal{J}_{k-1}| \le 10^{-3}$, where $\mathcal{J}_k = \frac{1}{2} \|\mv{L}(\mv{\delta}^{k-1})\mv{w}^k - \mv{s} \|^2$ represents the squared error between the measured sinogram $\mv{s}$ and the forward projection of the reconstructed image $\mv{w}^k$ calibrated with the drift estimate $\mv{\delta}^{k-1}$.

\subsection{Type I: Systematic Drifts}
\label{subsec:SpecialDrift}

In this section, we consider a simplified case where the scanning-position drifts are systematic across all projection angles, that is,
\[
\delta_{(\theta-1)N_\tau + \tau} = \delta_{(\theta'-1)N_\tau + \tau}, \quad \theta' \neq \theta.
\]
In other words, the drift depends only on the scanning index $\tau$, and we denote it simply by $\delta_\tau$. This scenario corresponds to the common experimental condition where the primary drift source is a systematic misalignment that can be reset after each rotation.

We begin with the simplest case where each drift $\delta_\tau$ is an integer multiple of the theoretical scanning step size $\Delta$. In this case, the relationship between the drifted forward model $\bar{\mv{L}}$ and the nominal, drift-free model $\mv{L}$ can be formally described as follows.

\begin{pro}[Integer drifts]
  \label{prop:integer}
  If $\forall \tau \in \{1,\dots,N_\tau\}$, $\delta_\tau = k_\tau \Delta$ for some integer $k_\tau$, then $\bar{\mv{L}}\mv{w} \subseteq \mv{L}\mv{w}$.\footnote{Here $\mv{a}\subseteq \mv{b}$ for two vectors indicates that the components of $\mv{a}$ are a subset of those of $\mv{b}$.} Moreover, there exists a binary row-selection matrix $\mv{P}\in \mathbb{R}^{N_{\theta}N_{\tau}\times N_{\theta}N_{\tau}}$ such that $\bar{\mv{L}} = \mv{P}\mv{L}$, where $\forall (\theta,\tau)$,
  \[
    \sum_{j=1}^{N_\theta N_\tau} \mv{P}_{(\tau-1)N_\theta+\theta,j} \le 1,
  \]
  and
  \begin{equation*}
  \label{eqn:p_int}
  \mv{P}_{(\tau-1)N_\theta+\theta,(\tau'-1)N_\theta+\theta}=
  \begin{cases}
    1, & \text{if beam } (\theta,\tau) \text{ is drifted to } (\theta,\tau'),\\
    0, & \text{otherwise.}
  \end{cases}
  \end{equation*}
\end{pro}

Next, we turn to the more general and practical case where $\delta_\tau$ is not an integer multiple of $\Delta$. In this regime, $\bar{\mv{L}}\mv{w}$ is not necessarily a subset of $\mv{L}\mv{w}$, and hence $\bar{\mv{L}}$ cannot be obtained purely by row selection. Instead, assuming that the one-dimensional projection (i.e., the sinogram along the $\tau$ direction) is smooth,\footnote{This is a weaker condition than assuming that the underlying 2D object itself is smooth, since the Radon integration reduces high-frequency variations.} we approximate
\[
\bar{\mv{L}} \approx \mv{P}\mv{L},
\]
where $\mv{P}$ interpolates between adjacent rows of $\mv{L}$. For example, with linear interpolation, each row of $\mv{P}$ contains only two adjacent nonzero entries:
\begin{equation}
  \label{eqn:p_frac}
  \mv{P}_{(\tau-1)N_\theta+\theta,(\tau'-1)N_\theta+\theta}=
  \begin{cases}
    \alpha_\tau, & \text{if } \tau'+1 = \lfloor \delta_\tau/\Delta \rfloor, \\[4pt]
    1-\alpha_\tau, & \text{if } \tau' = \lfloor \delta_\tau/\Delta \rfloor, \\[4pt]
    0, & \text{otherwise,}
  \end{cases}
\end{equation}
where
\[
\alpha_\tau = \frac{\delta_\tau}{\Delta} - \left\lfloor \frac{\delta_\tau}{\Delta} \right\rfloor \in [0,1)
\]
is the fractional part of the drift, and $\pi_\tau = \lfloor \delta_\tau / \Delta \rfloor$ denotes its integer component. Quadratic or higher-order interpolation can be similarly incorporated by allowing each row of $\mv{P}$ to have three or more adjacent nonzero entries. A quantitative comparison of interpolation orders is provided in \cref{sec:interp}.

To compute the current drift operator $\mv{P}^k$ in step~5 of \cref{alg:driftedTomography}, we reshape $\mv{L}\mv{w}^k$ and $\mv{P}\mv{L}\mv{w}^k$ into 2D matrices ${\mv{S}}^\star$ and $\bar{\mv{S}}$ of size $\mathbb{R}^{N_\theta\times N_\tau}$, respectively. The corresponding data-fitting term $\|\mv{P}\mv{L}\mv{w}^k-\mv{s}\|^2$ can then be expressed as
\[
\|\bar{\mv{S}} - \mv{S}\|_F^2,
\]
where $\|\cdot\|_F$ denotes the Frobenius norm. The regularization term $\lambda_k\Phi(\mv{w}^k)$ is omitted here since it does not depend on $\mv{P}$. For the case of linear interpolation,
\[
  \bar{\mv{S}}_{\cdot\tau} = (1-\alpha_\tau)\,{\mv{S}}^\star_{\cdot,\tau+\pi_\tau}
  + \alpha_\tau\,{\mv{S}}^\star_{\cdot,\tau+\pi_\tau+1},
\]
where $\bar{\mv{S}}_{\cdot\tau}$ denotes the $\tau$th column of $\bar{\mv{S}}$, and $\pi_\tau \in [-\pi_{\max},\pi_{\max}]$ with $\pi_{\max} = \lceil \Delta_{\max}/\Delta \rceil$.

\Cref{alg:computeDrift} summarizes the computation of the $N_\tau$ unknown drift values by matching the measured data $\mv{s}$ to the model prediction $\mv{L}\mv{w}^k$ at iteration $k$. Once the recovered drifts $\delta_\tau$ stabilize, the final reconstruction can be performed using the explicitly constructed $\bar{\mv{L}}$, rather than the interpolated approximation.

The computational complexity of \cref{alg:computeDrift} is dominated by Step~3, where each $\delta_\tau$ is estimated independently by sweeping $\pi_\tau$ over $[-\pi_{\max},\pi_{\max}]$. The corresponding fractional component $\alpha_\tau$ can be computed analytically by minimizing~\eqref{eqn:submin} with respect to $\alpha_\tau$. Hence, the total cost of solving~\eqref{eqn:submin} for $N_\theta$ projection angles and $N_\tau$ scanning positions is $O(N_\theta N_\tau \pi_{\max})$.

\begin{algorithm}[H]
\caption{Calibration of Drift Error}
\label{alg:computeDrift}
\begin{algorithmic}[1]
\State \textbf{Input:} $\mv{s}, \mv{L}, \mv{w}^k$
\State Construct 2D matrices $\mv{S}$ and $\mv{S}^\star$ from $\mv{s}$ and $\mv{L}\mv{w}^k$
\For{$\tau = 1,\ldots,N_\tau$}
  \State Solve
  \begin{equation}
  \label{eqn:submin}
  \min_{\substack{\pi_\tau\in [-\pi_{\max},\pi_{\max}]\\ \alpha_\tau \in [0,1]}}
  \sum_{\tau=1}^{N_\tau}
  \|{\mv{S}}_{\cdot\tau} - (1-\alpha_\tau){\mv{S}}^\star_{\cdot,\tau+\pi_\tau}
   - \alpha_\tau{\mv{S}}^\star_{\cdot,\tau+\pi_\tau+1}\|^2
  \end{equation}
  \State Compute $\delta_\tau = \pi_\tau + \alpha_\tau$ and construct $\mv{P}$ using \cref{eqn:p_frac}
\EndFor
\State \textbf{Output:} $\{\delta_\tau\}_{\tau=1}^{N_\tau}, \mv{P}$
\end{algorithmic}
\end{algorithm}

The Type~I formulation described above focuses on systematic drifts that are identical across all projection angles. This simplification captures many practical scenarios in which experimentalists re-align the scanning stage between rotations, resulting in consistent lateral shifts along each projection. However, in more challenging settings—such as long-duration or thermally unstable acquisitions—the drift can vary for each individual measurement, producing a spatially and temporally nonuniform distortion. To address these cases, we next extend our framework to a fully general drift model, referred to as \emph{Type~II}, in which each projection angle and scan position $(\theta,\tau)$ has its own independent drift variable $\delta_{\theta,\tau}$. The resulting formulation preserves the same optimization structure as in the Type~I case but introduces additional flexibility to handle complex, nonuniform drift patterns.

\subsection{Type II: General Drifts}  \label{subsec:GeneralDrift}

In this section, we consider a general type of drift where each individual measurement $(\theta, \tau)$ has its own drift $\delta_{\theta,\tau}$. Therefore, the true scanning step is $\Delta+\delta_{\theta,\tau}$. In this case, the total number of unknown drifts are $N_\theta \times N_\tau$ which is the same as the size of sinogram $\mv{S}$. Similar to representing $\mv{S}$ as a vector $\mv{s}$ in  column-major-order, we represent the  drifts as a vector $\mv{\delta} = [\delta_1; \delta_2; \cdots \delta_i; \cdots \delta_{N_\theta N_\tau}]$, where $i=\theta+(\tau-1)N_\theta$ is the one dimensional linear index for the measurement at $(\theta, \tau)$.

In this case of general drift, reconstruction of $\mv{w}, \mv{\delta}$ is harder than \cref{subsec:SpecialDrift}, because we have $N_\theta \times N_\tau$ unknown drifts $\mv{\delta}$ rather than only $N_\tau$ unknown drifts $\mv{\delta}$ in \cref{subsec:SpecialDrift}. Problem (\ref{eq:optd}) is further complicated because the Radon transform $R(\mv{w},\mv{\delta})=\mv{L}(\mv{\delta})\mv{w}$ becomes a nonlinear multivariate vector function of $\mv{\delta}$. However, we found that $R(\mv{w},\mv{\delta})$ is seperable on $\mv{\delta}$, because the $i$th measurement only depends on the $i$th but not other drifts. Hence we can simplify the multivariate vector function $R(\mv{w},\mv{\delta})$ as a vector composition of univariate scalar functions:
$R(\mv{\delta}; \mv{w}) = [R_1(\delta_1; \mv{w}), R_2(\delta_2; \mv{w}), \cdots, R_{N_\theta N_\tau}(\delta_{N_\theta N_\tau}; \mv{w})]$, where $R_i(\delta_i; \mv{w})$ is a scalar function depends only on variable $\delta_i$ given the parameter $\mv{w}$. 
Next, $R_i(\delta_i; \mv{w})$ can be approximated to a linear function interpolated from two endpoints $R_i(-\Delta_{\max}; \mv{w})$ and
$R_i(\Delta_{\max}; \mv{w})$. Combining this linear approximation with separable property on $\mv{\delta}$, we transform the originally nonlinear $R(\mv{\delta}; \mv{w})$ as a linear function of $\mv{\delta}$ where the linear transform is
a diagonal matrix. Problem (\ref{eq:optd}) then can be transformed to the following regularized linear least square problem with constraint:
\begin{equation}
  \label{eq:optd_constrained}
  \ds
  \min_{\delta_i\in [-\Delta_{\max}, \Delta_{\max}]}\frac{1}{2} \sum_{i=1}^{N_\theta N_\tau}\left(a_i \delta_i+b_i\right)^2+ \mu \Phi_d (\mv{\delta}),
\end{equation}
where
\begin{eqnarray*}
  a_i &=& \frac{R_i(\Delta_{\max}; \mv{w}) - R_i(-\Delta_{\max}; \mv{w})}{2 \Delta_{\max}},  \\
  b_i &=& \frac{R_i(\Delta_{\max}; \mv{w}) + R_i(-\Delta_{\max}; \mv{w})}{2 } -s_i.
\end{eqnarray*}
We solve~\eqref{eq:optd_constrained} using a modified TwIST, by adding the projection: $\delta_i = -\Delta_{\max}$ when $\delta_i < -\Delta_{\max}$ and $\delta_i = \Delta_{\max}$ when $\delta_i > \Delta_{\max}$, at every iteration.

\section{Choosing the Regularizer Parameter}\label{sec:autoLambda}

One critical aspect of successfully reconstructing $\mv{w}$ and $\mv{\delta}$ is the choice of the reconstruction regularization parameter $\lambda$. In this section, we focus on $\lambda$ in the image-reconstruction subproblem; the drift regularization parameter $\mu$ is held fixed in the Type~II numerical experiments. Assuming the effective measurement discrepancy is $\mv{\epsilon}$, our forward model becomes
\begin{equation*}\label{eq:forwardModelWithNoise}
\mv{s} = \mv{L}\mv{w} +\mv{\epsilon}.
\end{equation*}
The fact that the noise is of the order of $\mv{\epsilon}$ motivates a widely used posteriori-stopping rule, the \textit{discrepancy principle}, i.e.,  $\|\mv{L}\mv{w}-\mv{s}\|\leq \eta \|\mv{\epsilon}\|$, where $\eta$ is a safety factor slightly larger than 1 \cite{morozov2012methods}.
Notice that the major source of discrepancy $\mv{\epsilon}$ includes both the measurement noise (e.g., the quantum noise) and the experimental error (e.g., the electronic noise and detector system drifts). In the context of this work where our goal is to calibrate the drift error, we expect that $\mv{\epsilon}$ will reduce along the outer iteration $k$ following Alg.~\ref{alg:driftedTomography}. Therefore, instead of following the standard trial and error to pick a sub-optimal parameter, we explore the adaptive regularizer parameter approach developed in \cite{gazzola2015,gazzola2020krylov} and propose a robust way for automating the choices of regularization parameter $\lambda$ in \cref{alg:driftedTomography} to account for the model/data discrepancy.

In practice, $\|\mv{\epsilon}\|$ is never exactly known. In this work, we assume $\mv{\epsilon}$ follows a Gaussian distribution with 0 mean and standard deviation $\sigma$.

Let
\begin{equation*} \label{eq:wstep2}
    \mv{w}^{k,\lambda}=\argmin_{\mv{w} \geq 0}\frac{1}{2}\|R(\mv{w},\mv{\delta}^k)- \mv{s}\|^2 + \lambda \Phi(\mv{w}),
\end{equation*}

and
\begin{equation*}
\phi_k(\lambda) = \|R(\mv{w}^{k,\lambda},\mv{\delta}^k) - \mv{s}\|.
\end{equation*}

Ideally when the forward model is error-free and $\mv{w}$ is solved accurately, $\phi_k(\lambda)$ should be equal to the noise in $\mv{s}$, that is, $\mv{\epsilon}$. Hence in each iteration $k$, we solve for $\lambda$ such that
\begin{equation*}
    \phi_k(\lambda) \approx s_{\max} \sqrt{N_\theta N_\tau} \sigma,
\end{equation*}
where $s_{\max}$ is the largest element in $\mv{s}$. This estimate is conservative for the experiments considered here, so we do not introduce an additional safety factor.

We further assume that $\phi_k(\lambda)$ can be well approximated by the linear function
\begin{equation*}
    \phi_k(\lambda) \approx \alpha_k + \beta_k\lambda,
\end{equation*}
and Taylor series suggest that
\begin{eqnarray*}
    \alpha_k = \phi_k(0) = \|R(\mv{w}^{k,0},\mv{\delta}^k) - \mv{s}\|.
\end{eqnarray*}

To compute $\beta_k$, consider the approximation $\phi_k(\lambda_{k-1}) \approx \alpha_k + \beta_k\lambda_{k-1}$, hence
\begin{equation*}
    \beta_k \approx \frac{\phi_k(\lambda_{k-1}) - \alpha_k}{\lambda_{k-1}}.
\end{equation*}
Therefore when we solve for $\phi_k(\lambda_k) = \alpha_k + \beta_k\lambda_k =  \|\mv{\epsilon}\|$, we obtain
\begin{equation*}
    \lambda_k = \frac{\|\mv{\epsilon}\| - \alpha_k}{\beta_k} = \frac{\|\mv{\epsilon}\| - \phi_k(0)}{\phi_k(\lambda_{k-1}) - \phi_k(0)}\lambda_{k-1}.
\end{equation*}

Note that in practice, one can use the following to preserve the nonnegative of $\lambda$,
\begin{equation}\label{eq:lambdaupdate}
    \lambda_k = \left\lvert \frac{\|\mv{\epsilon}\| - \phi_k(0)}{\phi_k(\lambda_{k-1}) - \phi_k(0)}\right\rvert\lambda_{k-1}.
\end{equation}

In order to have the initial $\lambda_0$ reflect the level of discrepancy from detection noise and experimental drifts, we propose the following formula:
\begin{equation*}\label{eq:setlam0}
        \lambda_0=0.001\left(1 + \sigma \right)\left(1+\Delta_{\max}\right).
\end{equation*}
This heuristic increases the starting regularization weight when either the sinogram noise level or the maximum admissible drift increases, while keeping $\lambda_0$ in the empirically useful range of approximately $10^{-3}$ to $10^{-2}$ for the experiments considered here. The subsequent adaptive update in~\eqref{eq:lambdaupdate} then decreases or increases $\lambda_k$ according to the observed discrepancy at each outer iteration.

\section{Interpolation Error Analysis}
\label{sec:interp}

In this section, we analyze the numerical effects of different interpolation schemes used to approximate the drifted forward operator $\bar{\mv{L}} \approx \mv{P}\mv{L}$. In particular, we compare the performance of linear and quadratic interpolation in both the forward simulation and the backward recovery processes. These tests isolate the interpolation error by evaluating the forward sinogram approximation and the backward drift-recovery step with the relevant ground-truth quantities fixed.

In the forward analysis, we evaluate how accurately the interpolated operator $\mv{P}\mv{L}(0)$ can reproduce the sinogram obtained from the true drifted model $\mv{P}^*\mv{L}(0)\mv{w}^*$. The corresponding error is measured by
\[
\|\mv{P}^*\mv{L}(0)\mv{w}^* - \mv{s}\|,
\]
where $\mv{P}^*$ encodes the true drift pattern $\mv{\delta}^*$, $\mv{w}^*$ is the ground-truth object, and $\mv{s}$ denotes the measured sinogram. This comparison quantifies how accurately the interpolation scheme approximates the true physical forward model.

In the backward analysis, we assess the accuracy of the recovered drifts $\mv{\delta}_k$ obtained by solving the optimization problem~\eqref{eq:optd_constrained} using the ground-truth object $\mv{w}^*$. The reconstruction error is defined as
\[
\|\mv{\delta}^* - \mv{\delta}_k\|,
\]
which measures the discrepancy between the true and estimated drift fields.

Table~\ref{tab:error_analysis} summarizes the quantitative results, reported as peak signal-to-noise ratio (PSNR, in dB), for two representative test objects (Brain and Phantom) at three spatial resolutions. The forward error reflects the fidelity of the interpolated sinograms, while the backward error indicates the accuracy of the recovered drifts. The results show that linear and quadratic interpolation yield very similar performance in both forward and backward analyses. Although quadratic interpolation provides slightly higher accuracy at higher resolutions, the improvement is marginal compared with the substantially increased computational cost. Therefore, we adopt linear interpolation for all subsequent reconstructions owing to its simplicity and efficiency.

\setlength{\tabcolsep}{.4em}
\begin{table*}[!t]
\caption{Quantitative comparison between linear and quadratic interpolation at three spatial resolutions for two test samples. 
Results are reported as PSNR (dB). The Brain dataset, being smoother, yields smaller forward errors but larger drift-recovery errors than the artificial Phantom. Given the negligible accuracy difference between interpolation orders, we employ linear interpolation throughout this work.}
\label{tab:error_analysis}
\begin{center}
\begin{tabular}{|c|ccc|ccc|ccc|ccc|}
\hline
Sample Error & \multicolumn{3}{c|}{Forward: Brain} & \multicolumn{3}{c|}{Forward: Phantom} & \multicolumn{3}{c|}{Backward: Brain} & \multicolumn{3}{c|}{Backward: Phantom} \\
\hline
Image Resolution & 50 & 100 & 200 & 50 & 100 & 200 & 50 & 100 & 200 & 50 & 100 & 200 \\
\hline \hline
\rowcolor{highlight} Linear & 40.18 & 45.00 & 51.09 & 29.15 & 34.56 & 40.77 & 20.95 & 28.66 & 25.91 & 24.56 & 22.10 & 23.69 \\
\rowcolor{Gray} Quadratic & 40.60 & 45.68 & 51.81 & 29.18 & 34.83 & 40.97 & 22.04 & 28.79 & 26.07 & 24.09 & 20.09 & 20.09 \\
\hline
\end{tabular}
\end{center}
\end{table*}

\section{Experimental Results} \label{sec:simulation}
\begin{figure}[htb]
  \centering
  \begin{subfigure}{0.32\linewidth}
    \includegraphics[width=\linewidth]{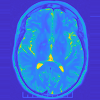}
    \caption{Brain}
  \end{subfigure}
 \begin{subfigure}{0.32\linewidth}
    \includegraphics[width=\linewidth]{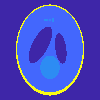}
    \caption{Phantom}
  \end{subfigure}
 \begin{subfigure}{0.055\linewidth}
    \includegraphics[width=\linewidth]{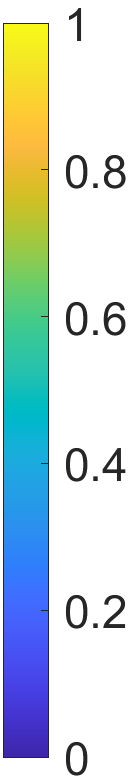}
    \caption*{}
  \end{subfigure}
\caption{Test images and colormap.}
\label{fig:testObj}
\end{figure}

To validate the efficacy of our proposed method, we compare its performance with the regularized baseline reconstruction that directly solves
problem~\eqref{eq:P0}, as well as the one proposed in \cite{huang2019calibrating} where the regularization parameter is chosen empirically. We confine ourselves to the TV regularization as one of the most popular choices in the community of tomographic
reconstruction \cite{li2011fast,mahmood2018adaptive}. We choose two images as the testing objects: the standard synthetic Shepp-Logan phantom
and the real MRI brain image, which are shown in Fig.~\ref{fig:testObj} as the ground truth. We map the original grayscale test images with the
colormap in Fig.~\ref{fig:testObj} to reveal more details.

We mainly compare the reconstruction performance between the baseline method (i.e., no drift calibration is performed) and our proposed reconstruction with drift calibration. Furthermore, notice that during the calibration process, the reconstruction regularizer parameter $\lambda$ plays a vital role in terms of the final solution quality since it gradually adjusts the discrepancy between model and data as the forward model is calibrated. Therefore, we also compare the performance between two approaches to sequentially update $\lambda_k$. The first approach is adopted from our previous work \cite{huang2019calibrating} where the $\lambda_k$ is linearly decreased as the drift is being calibrated, which we refer to as ``diminishing $\lambda_k$''. The second approach is based on the formula~\eqref{eq:lambdaupdate} we propose in this work, which we refer to as ``adaptive $\lambda_k$''.

For the experimental configuration, we fix the object discretization as $N=256$, and we choose the scan range $N_\tau\Delta$ slightly larger than $\sqrt{2}N \Delta + 10\Delta$, to
guarantee the full coverage of the object range $\sqrt{2}N \Delta$ with at least 5 scans outside of the object of interest on both sides.

In Alg.~\ref{alg:driftedTomography}, we set the stopping criteria as
\begin{equation}
    \label{eqn:stop}
    \frac{1}{2}\|R(\mv{w}^k,\mv{\delta}^k)-\mv{s}\|^2 - \frac{1}{2}\|R(\mv{w}^{k-1},\mv{\delta}^{k-1})-\mv{s}\|^2 \leq 10^{-3}.
\end{equation} Notice that as a comparison, in~\cite{huang2019calibrating}, the stopping criteria is set to reach the maximum number of iterations allowed which also determines the updating of $\lambda$, therefore, in the numerical experiments, in order to have a fair comparison, we will follow the stopping criteria set in~\cite{huang2019calibrating} but indicate where the algorithm would terminate if using the new stopping criteria.

\subsection{Type I Drift Result}

\setlength{\tabcolsep}{.4em}
\begin{table*}[!t]
\caption{Type I drift performance comparison.
Blue row: PSNR (dB); gray row: SSIM.  }
\label{tab:PSNR_SSIM}
  \begin{center}
    \begin{tabular}{|c|ccc|ccc|ccc|ccc|}
\hline
Max Drift &   & 1 &   &  & 2 &  &  & 3 &  &  & 5 &  \\
\hline
Noise & 0 & 0.01 & 0.02 & 0 & 0.01 & 0.02 & 0 & 0.01 & 0.02 & 0 & 0.01 & 0.02 \\
\hline \hline
\rowcolor{highlight}Brain (adaptive $\lambda_k$) & 29.88 & 27.19 & 25.17 & 29.12 & 26.83 & 25.19 & 27.14 & 25.74 & 24.49 & 25.73 & 25.43 & 23.98\\
\rowcolor{Gray} & 0.8778 & 0.8133 & 0.7474 & 0.8722 & 0.7989 & 0.7402 & 0.8241 & 0.7645 & 0.7177 & 0.7742 & 0.7725 & 0.7011 \\
\hline
\rowcolor{highlight} Brain (diminishing $\lambda_k$) & 27.65 & 23.95 & 25.79 & 27.65 & 26.86 & 25.44 & 23.75 & 23.59 & 23.41 & 23.50 & 23.37 & 23.25  \\
\rowcolor{Gray} & 0.8417 & 0.7143 & 0.6916 & 0.8399 & 0.7916 & 0.6642 & 0.7068 & 0.6991 & 0.6875 & 0.6968 & 0.6888 & 0.6808 \\
\hline
\rowcolor{white} Brain (baseline) & 26.24 & 23.23 & 19.15 & 21.15 & 20.21 & 17.65 & 18.95 & 18.07 & 16.08 & 15.51 & 15.06 & 13.75 \\
\rowcolor{Gray} & 0.6893 & 0.4548 & 0.2677 & 0.4229 & 0.3111 & 0.1979 & 0.2616 & 0.2031 & 0.1320 & 0.1167 & 0.0977 & 0.0672 \\
\hline \hline
\rowcolor{highlight} Phantom (adaptive $\lambda_k$) & 26.77 & 26.47 & 24.71 & 25.12 & 24.65 & 23.26 & 23.12 & 22.89 & 22.61 & 21.54 & 21.49 & 21.27\\
\rowcolor{Gray} & 0.9595 & 0.8136 & 0.8695 & 0.9557 & 0.8501 & 0.8514 & 0.9189 & 0.7043 & 0.8116 & 0.8725 & 0.7265 & 0.7241  \\
\hline
\rowcolor{highlight} Phantom (diminishing $\lambda_k$) & 25.95 & 25.99 & 24.73 & 25.00 & 24.36 & 23.44 & 22.78 & 22.37 & 19.72 & 19.62 & 19.6574 & 20.92 \\
\rowcolor{Gray} & 0.9519 & 0.9200 & 0.8476 & 0.9509 & 0.9094 & 0.8420 & 0.9118 & 0.8685 & 0.8146 & 0.7962 & 0.7672 & 0.7407 \\
\hline
\rowcolor{white} Phantom (baseline) & 21.34 & 20.87 & 19.71 & 16.98 & 16.76 & 16.21 & 14.58 & 14.47 & 14.08 & 11.81 & 11.75 & 11.41 \\
\rowcolor{Gray} & 0.7245 & 0.6143 & 0.4944 & 0.5487 & 0.4948 & 0.4235 & 0.4509 & 0.4160 & 0.3701 & 0.3639 & 0.3383 & 0.3071 \\
\hline
\end{tabular}
\end{center}
\end{table*}

We first demonstrate the performance on Type I drift. We perform tests for different levels of drifts separately; namely, $\max\limits_{\tau}\Delta_\tau \in \{\Delta, 2\Delta, 3\Delta,
5\Delta\}$.  For each case, to further test the robustness of our proposed method, we add two levels of Gaussian noise to the sinogram
(simulated with the assigned drifts) with standard deviation $\sigma$ (relative to the maximum intensity of its corresponding noise-free sinogram) as
$0.01$ and $0.02$, respectively.
We adopt two error metrics, the peak signal-to-noise ratio (PSNR) and the
structural similarity index (SSIM) \cite{hore2010image}, to quantify the reconstruction quality. All experiments use
$\lambda=1.3\times10^{-5}$, which is the weight that yields the best PSNR for the baseline algorithm for the brain sample with maximum
drift $\Delta$.

Table~\ref{tab:PSNR_SSIM} shows the full comparison results for Type I drift. Both drift-calibrated approaches substantially improve over the uncalibrated baseline in most noise and drift settings. The adaptive $\lambda_k$ update gives the highest PSNR or SSIM in many cases, while the diminishing $\lambda_k$ update remains competitive for several high-noise and high-drift settings.

\begin{figure}[!tb]
    \centering
    \setlength{\tabcolsep}{2pt}
    \begin{tabular}{ccc}
    \includegraphics[width=0.31\linewidth]{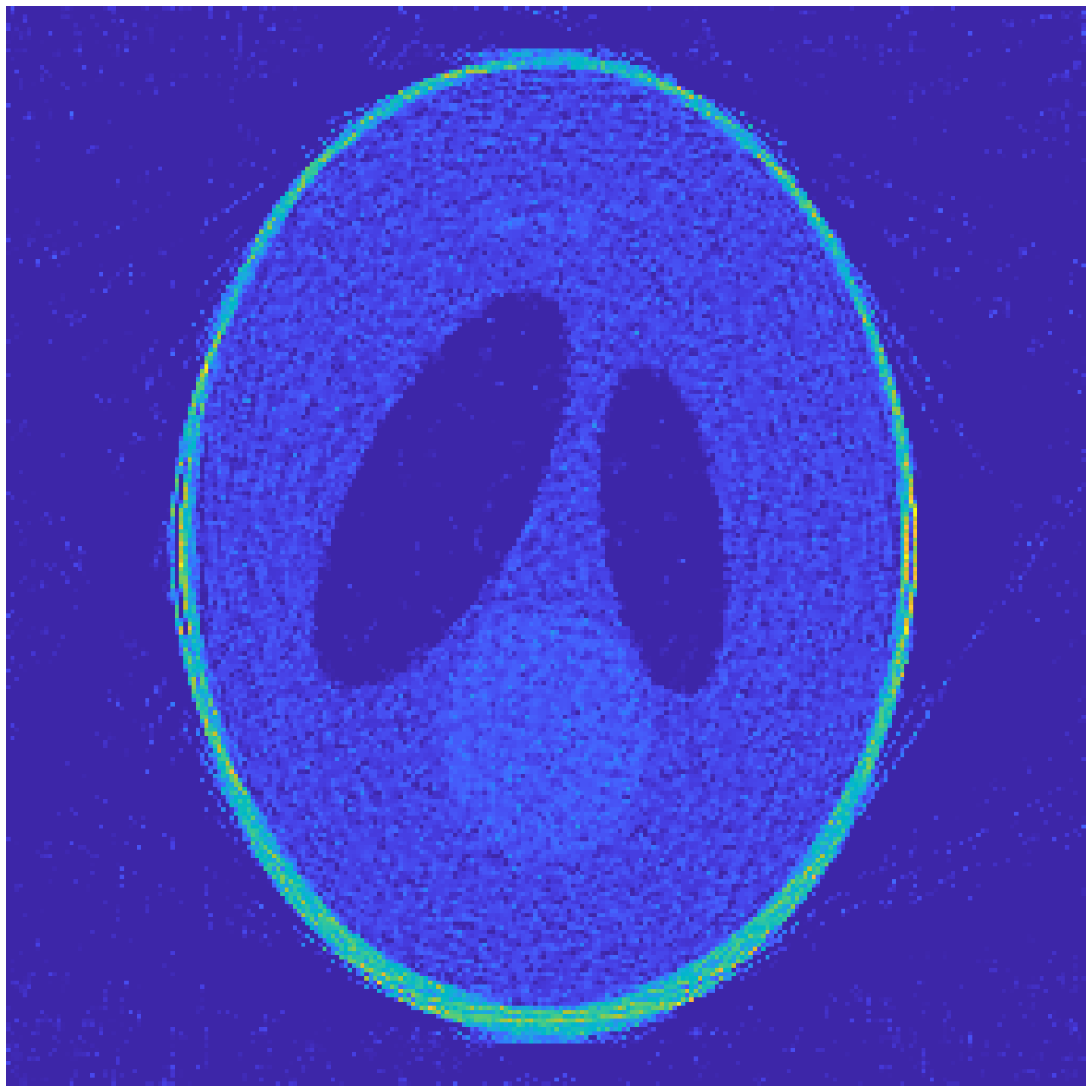} &
    \includegraphics[width=0.31\linewidth]{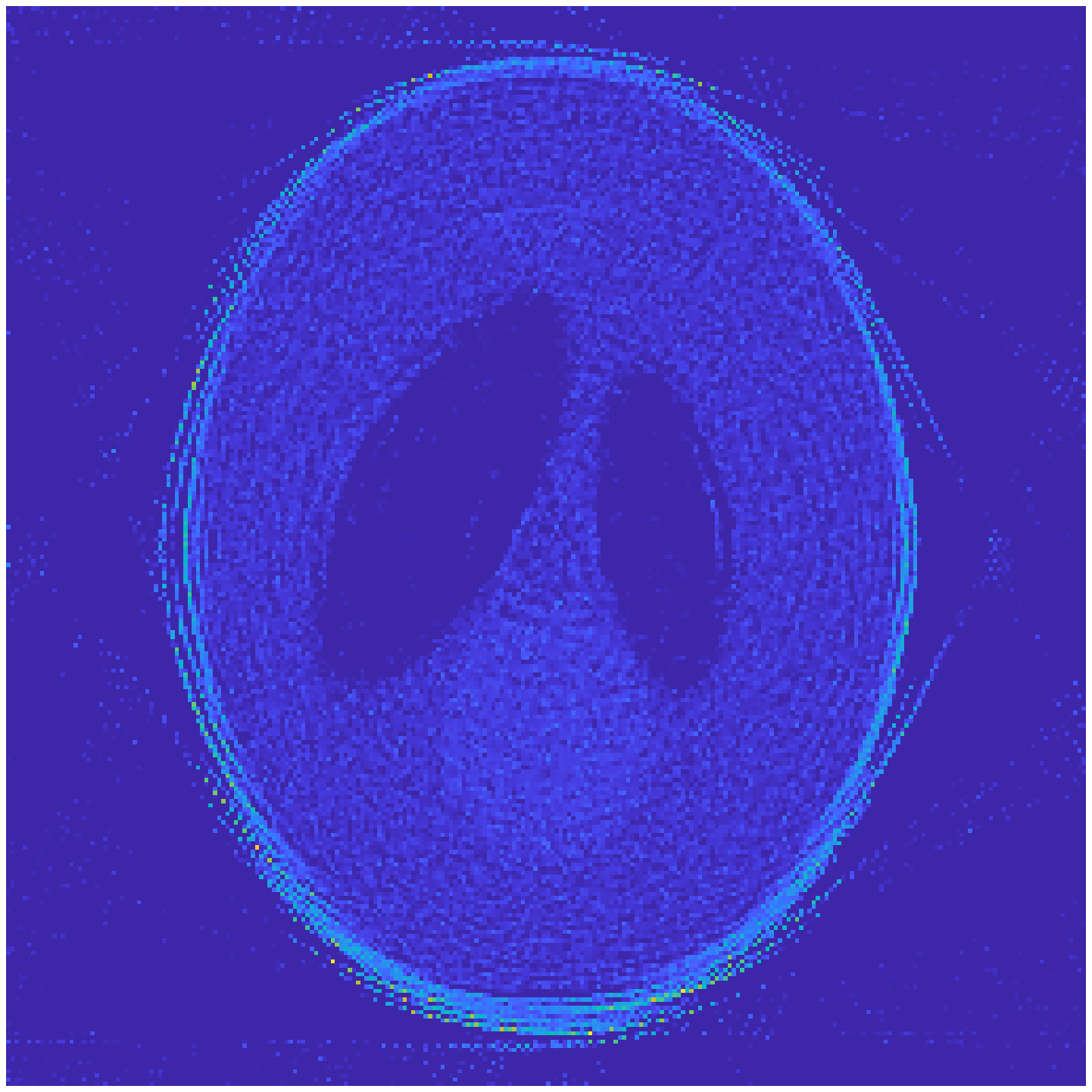} &
    \includegraphics[width=0.31\linewidth]{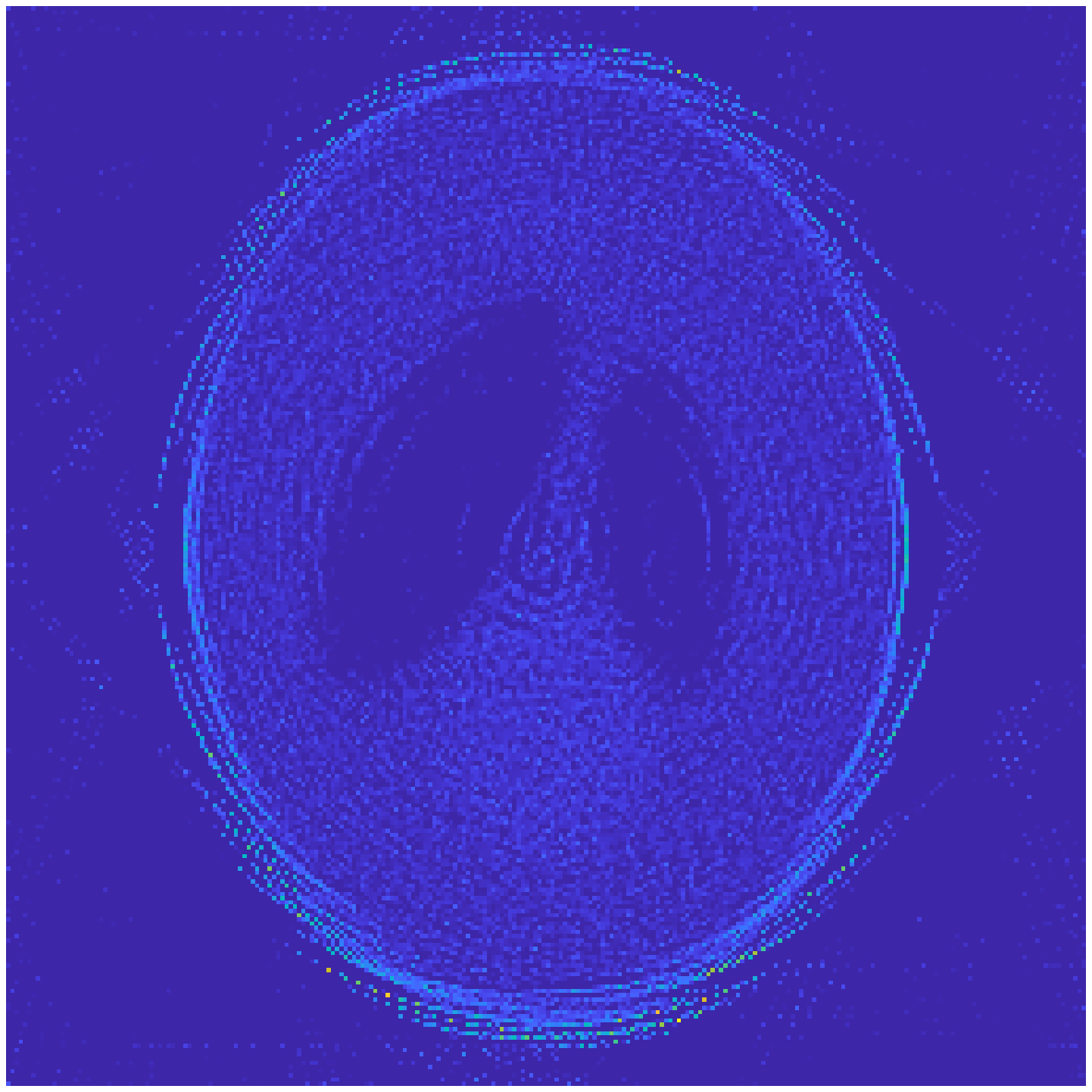} \\

    \includegraphics[width=0.31\linewidth]{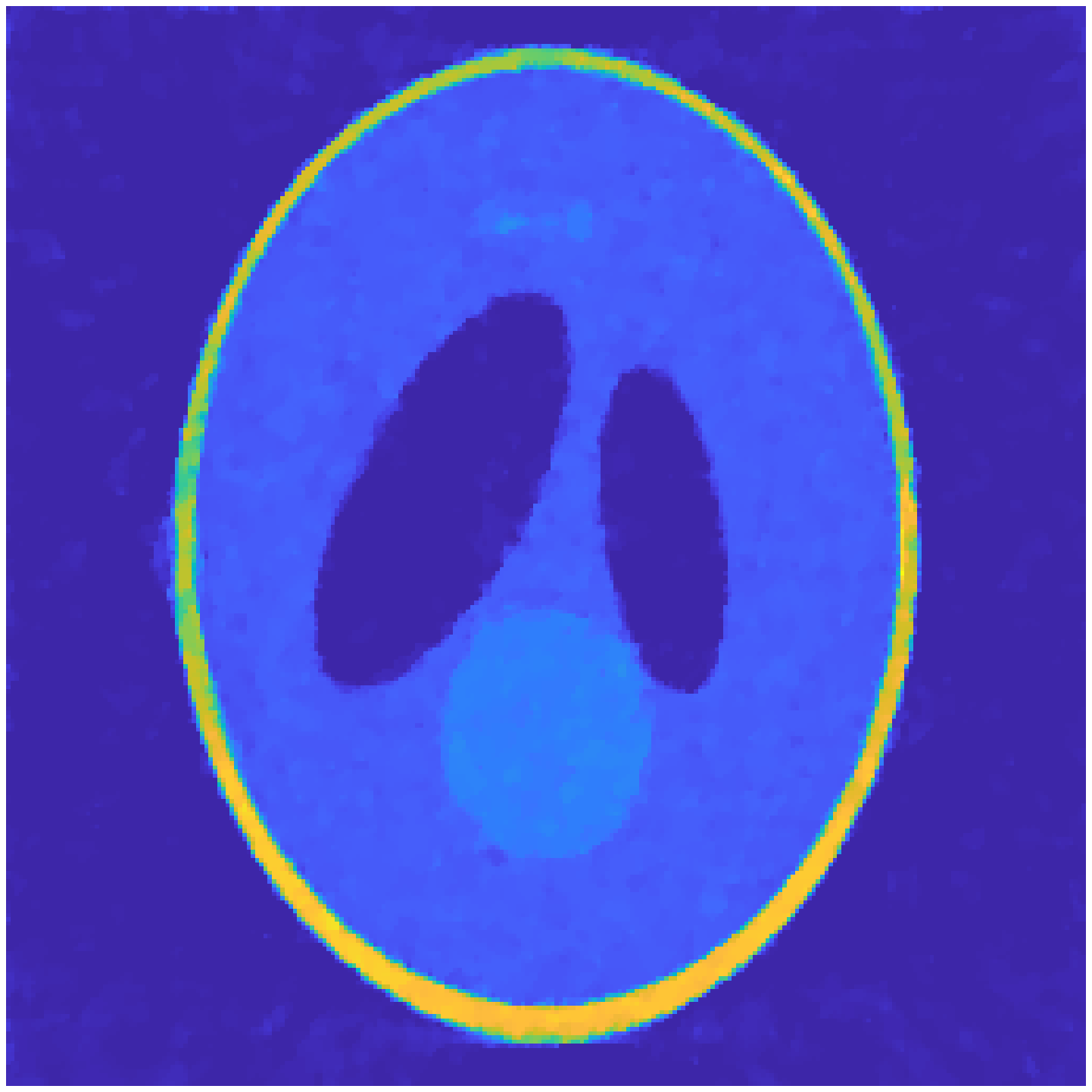} &
    \includegraphics[width=0.31\linewidth]{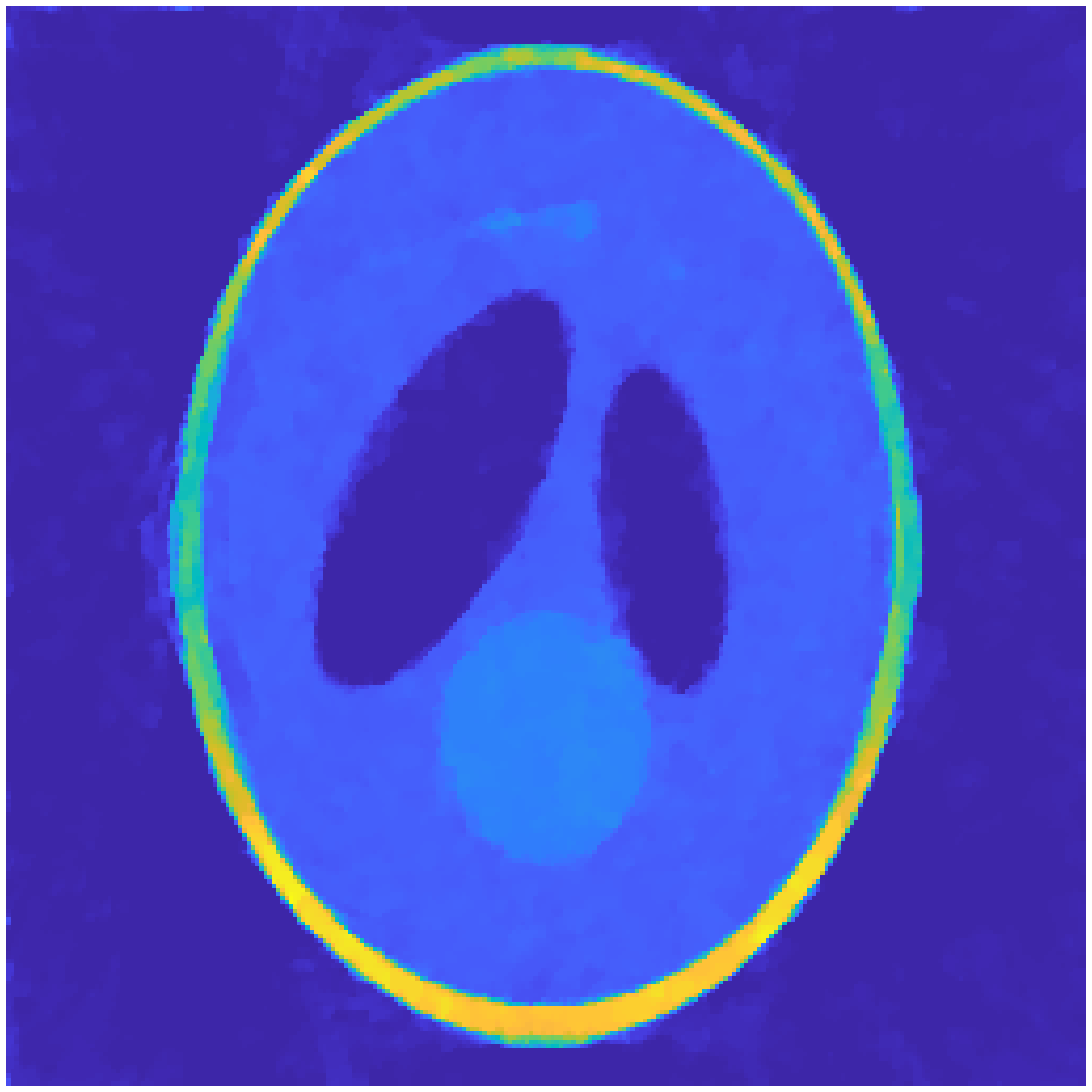} &
    \includegraphics[width=0.31\linewidth]{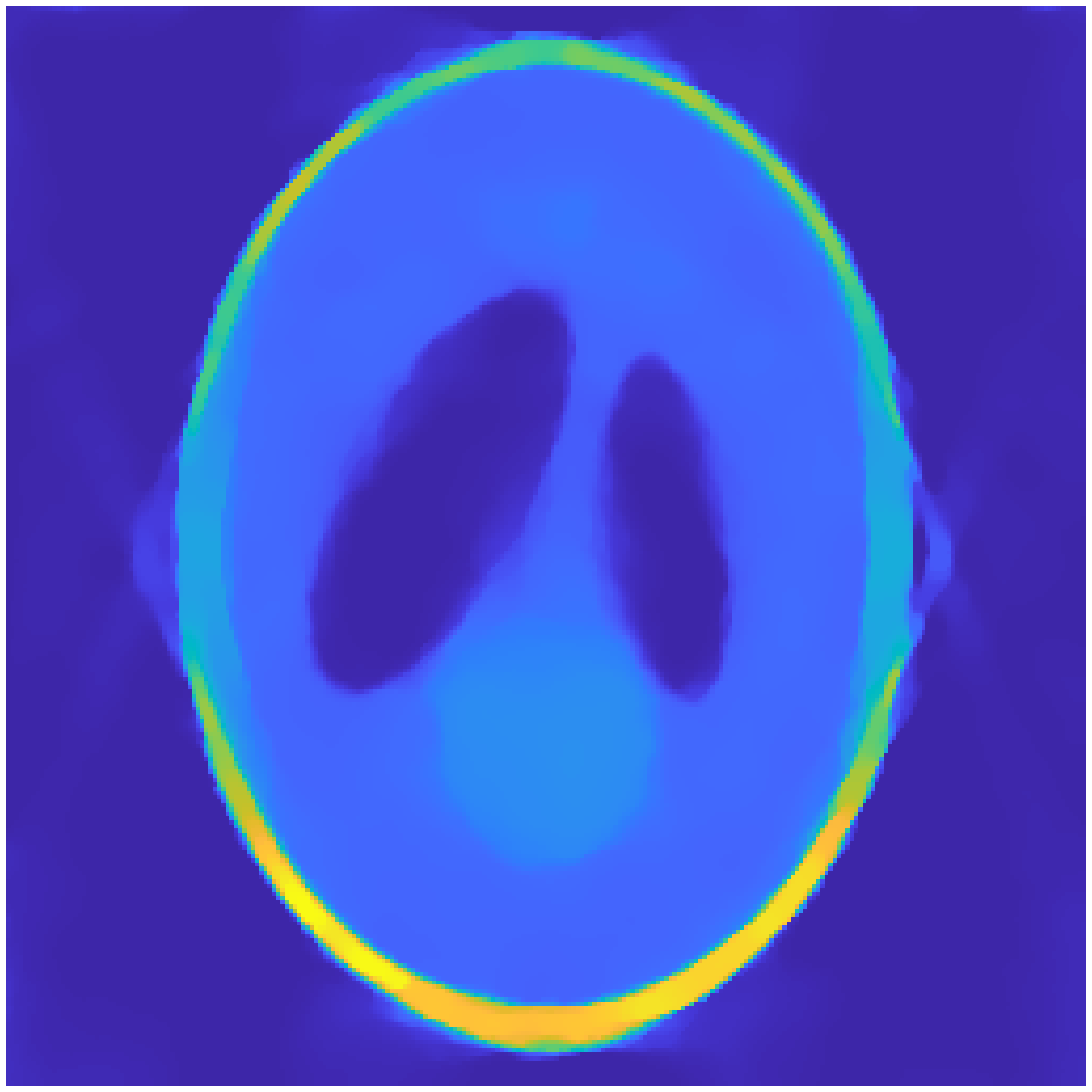} \\

    \includegraphics[width=0.31\linewidth]{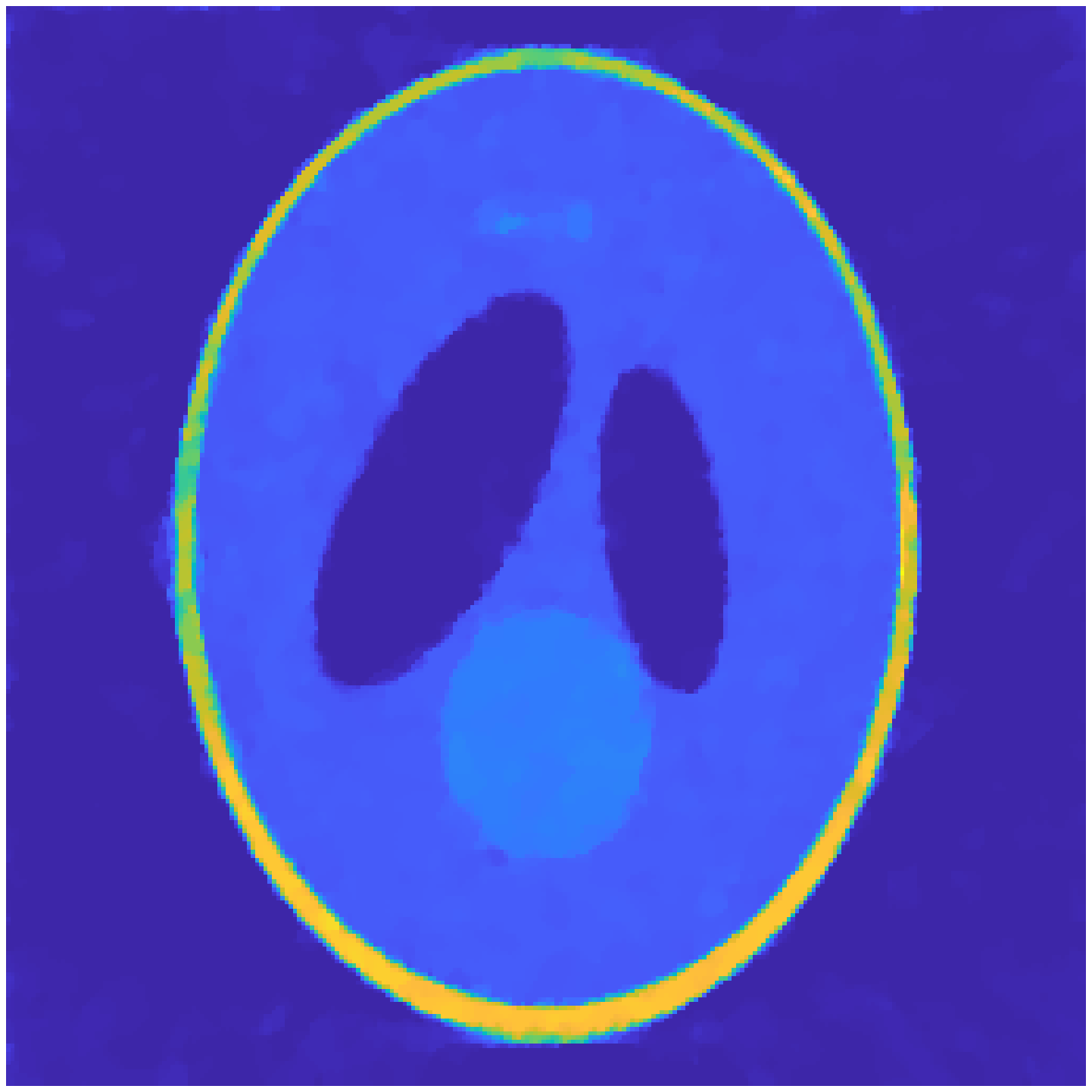} &
    \includegraphics[width=0.31\linewidth]{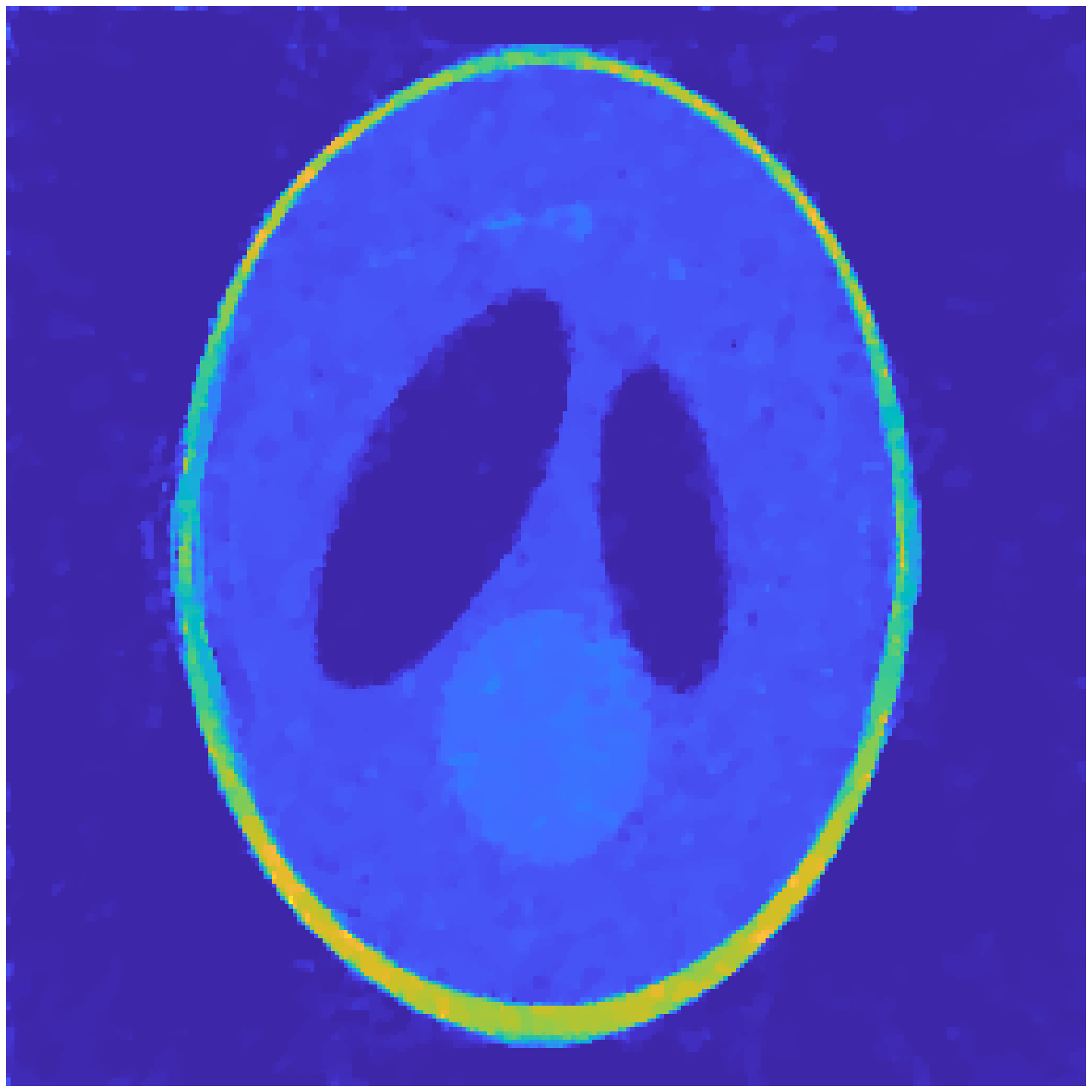} &
    \includegraphics[width=0.31\linewidth]{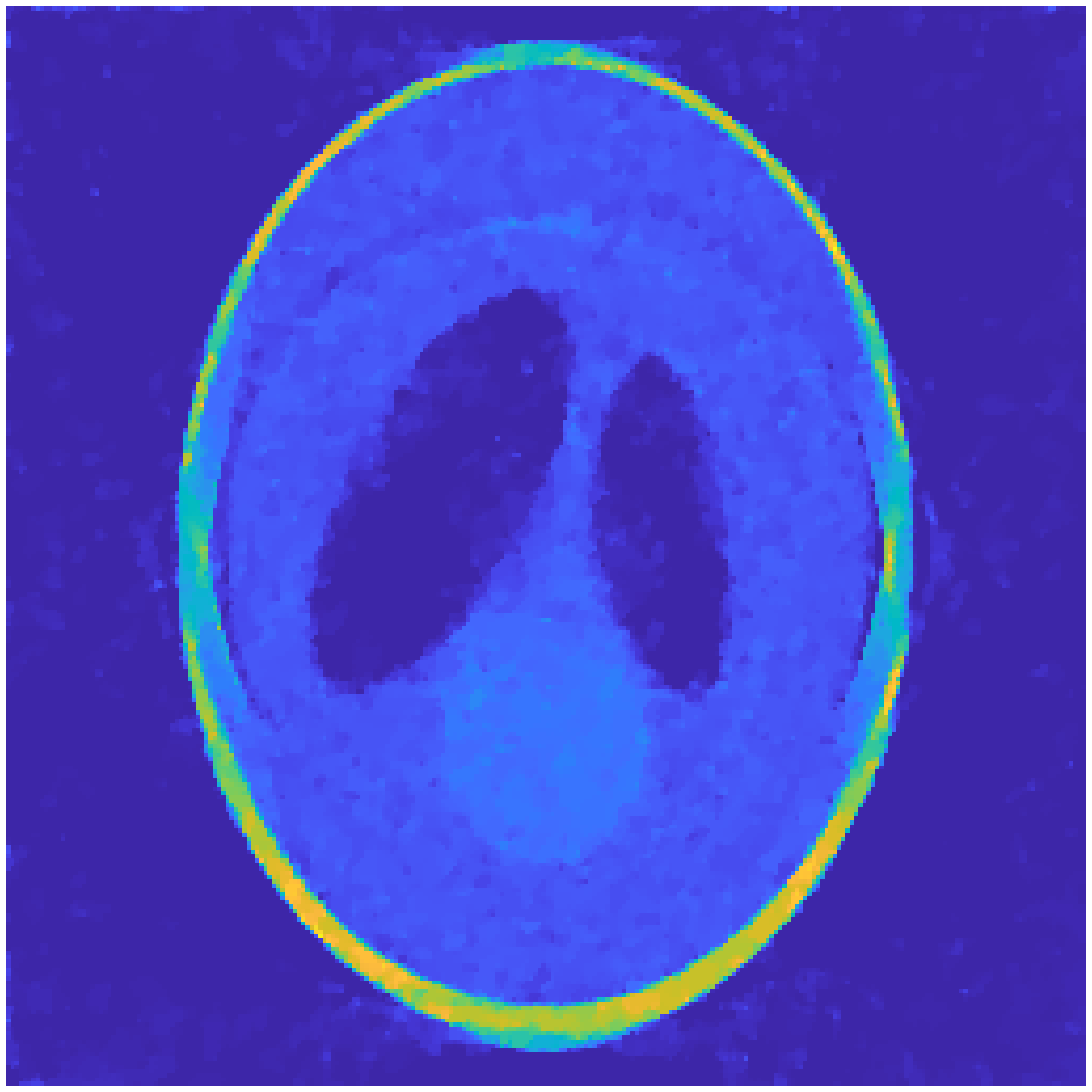} \\

    {\small $\Delta_{\max}=1,\, \sigma=0.02$} &
    {\small $\Delta_{\max}=3,\, \sigma=0.02$} &
    {\small $\Delta_{\max}=5,\, \sigma=0.02$} \\
    \end{tabular}

    \caption{256 by 256 phantom. Top row: baseline; middle row: diminishing $\lambda_k$ approach; bottom row: new approach}
    \label{fig:phantype1recs1-002}
\end{figure}

\begin{figure}[!tb]
    \centering
    \setlength{\tabcolsep}{2pt}
    \begin{tabular}{ccc}
    \includegraphics[width=0.31\linewidth]{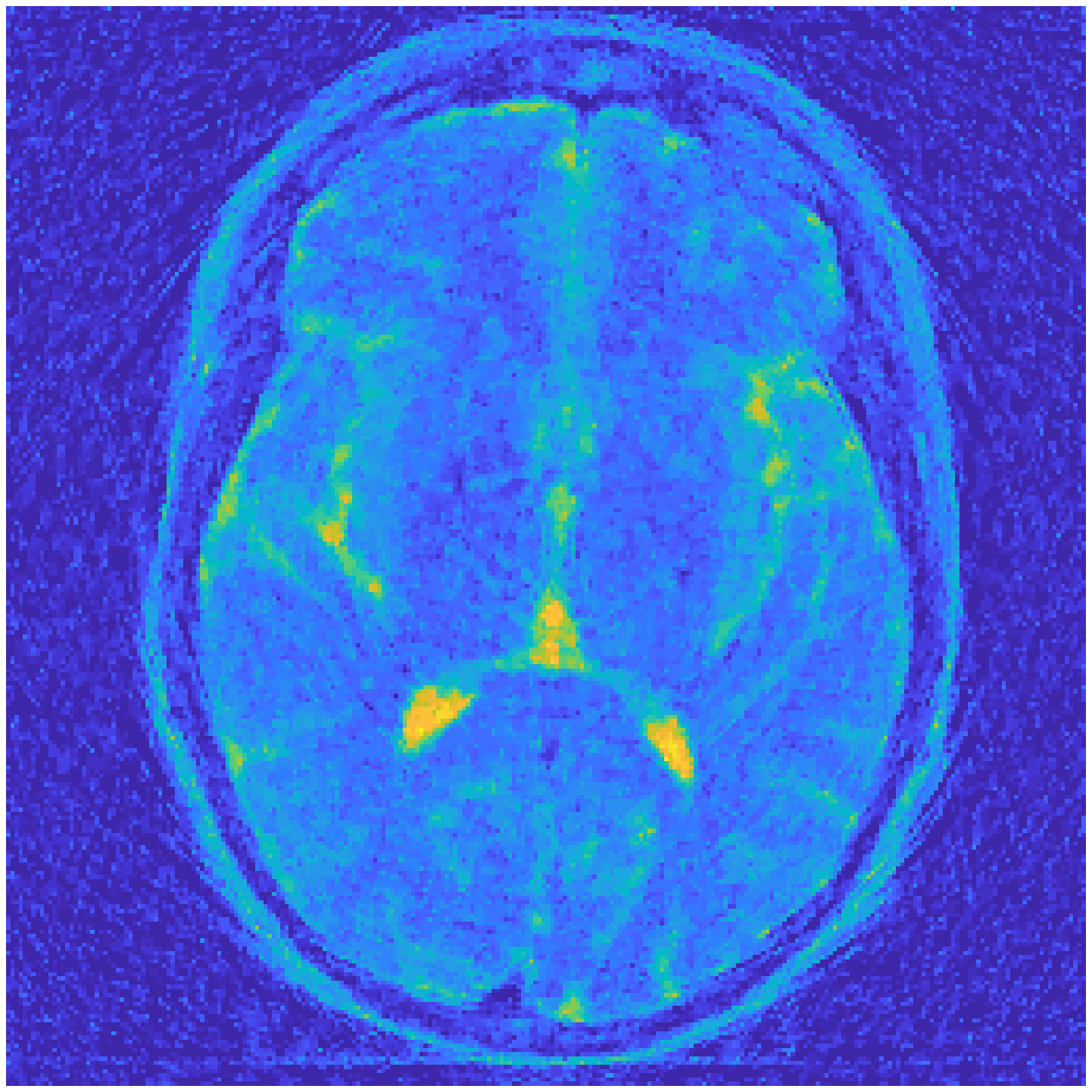} &
    \includegraphics[width=0.31\linewidth]{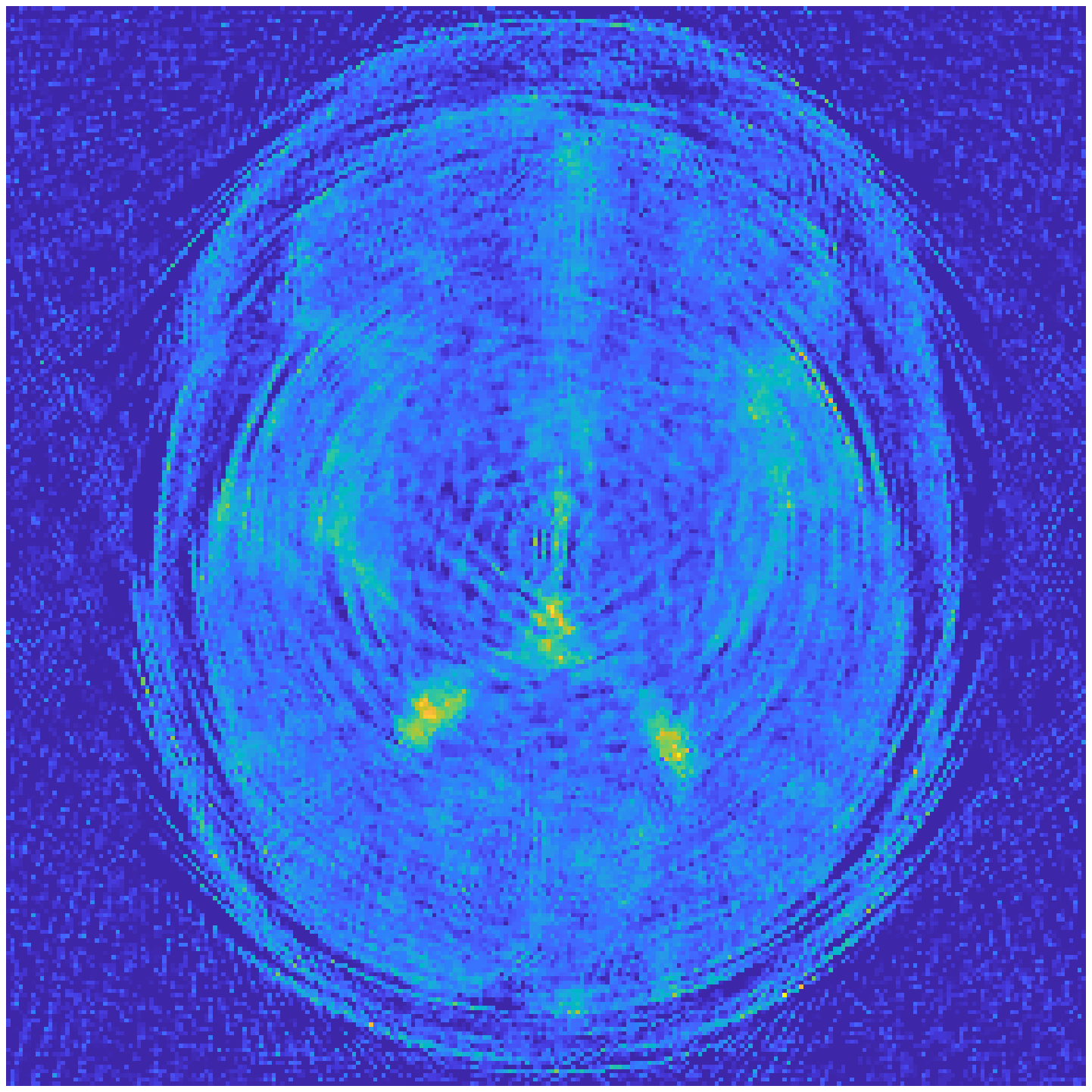} &
    \includegraphics[width=0.31\linewidth]{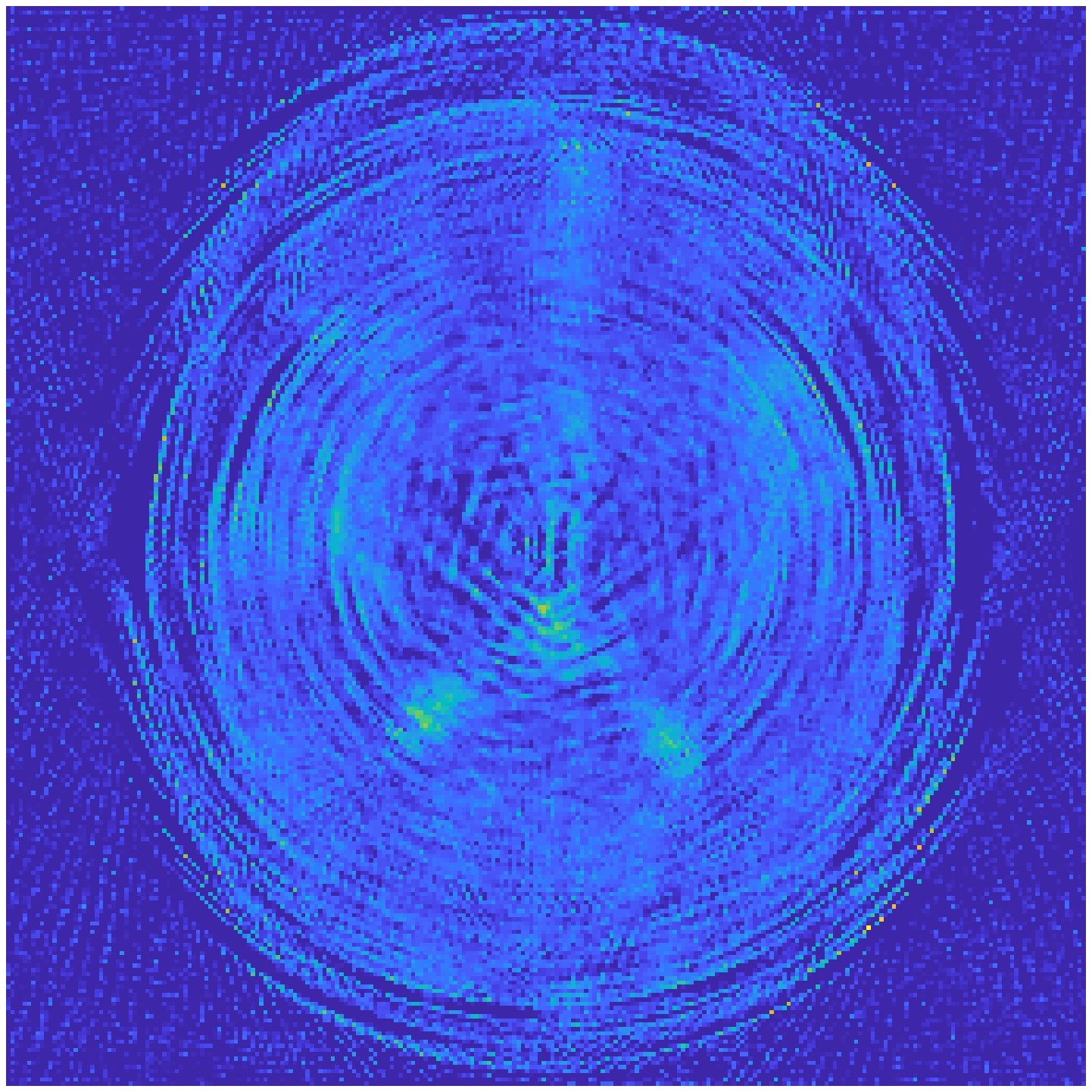} \\

    \includegraphics[width=0.31\linewidth]{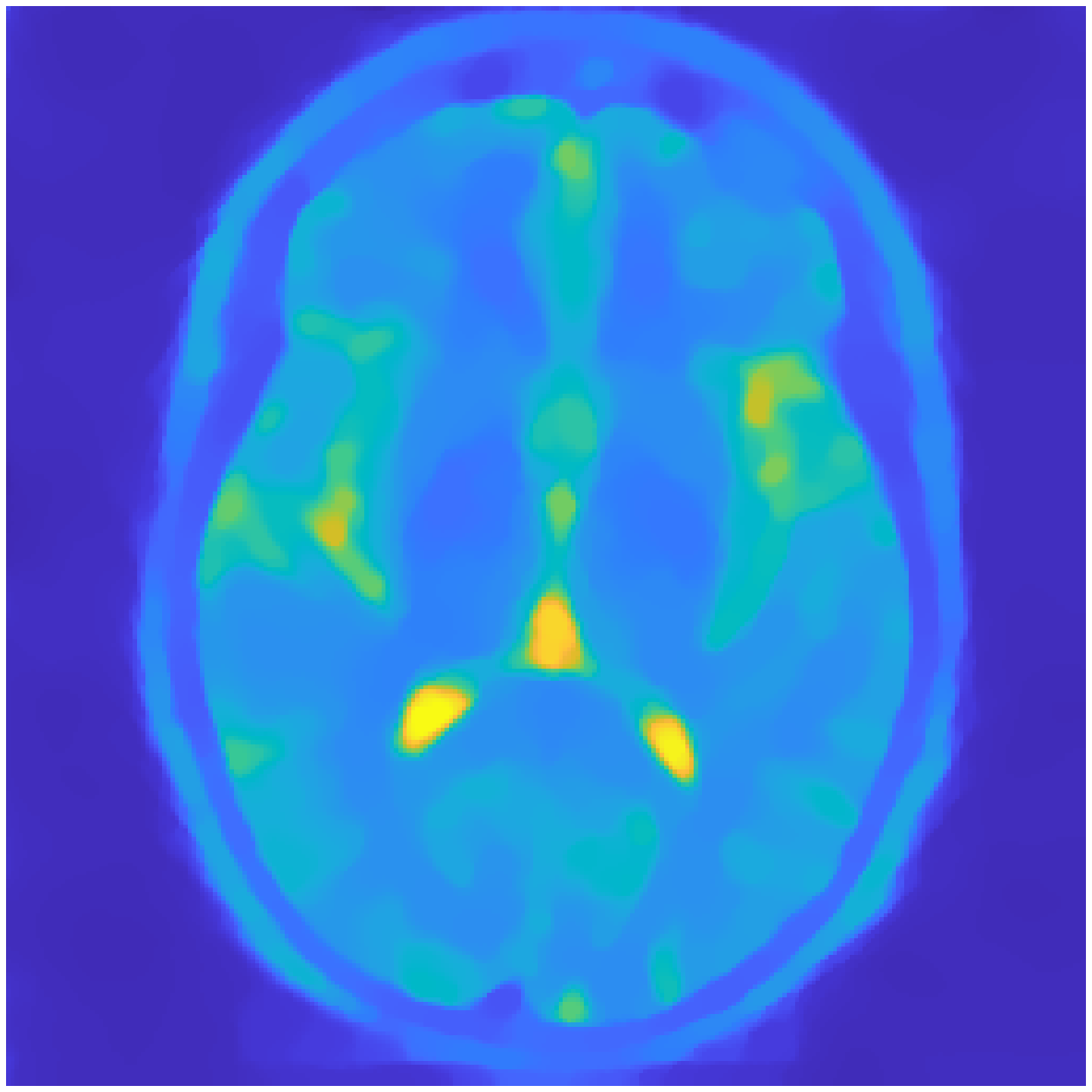} &
    \includegraphics[width=0.31\linewidth]{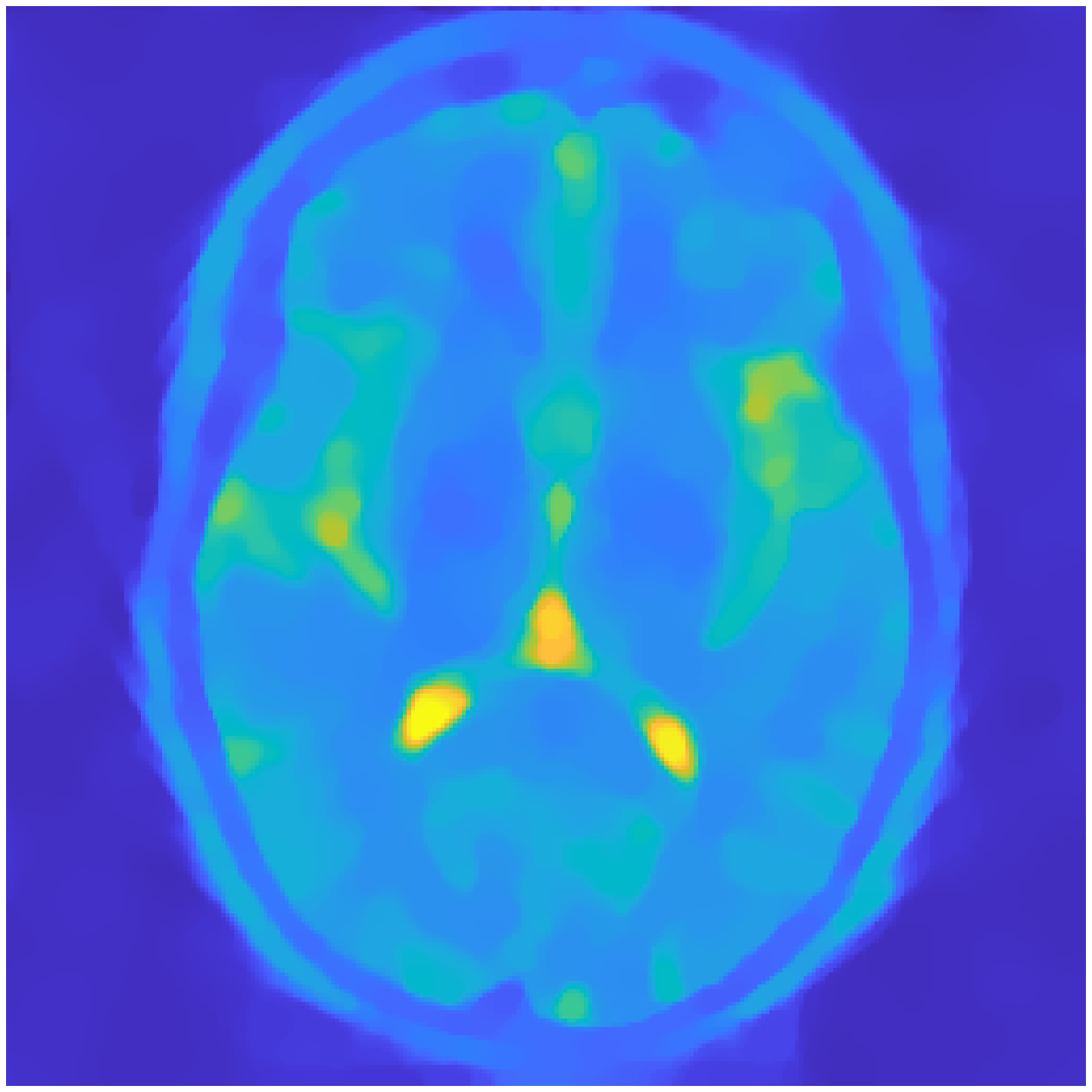} &
    \includegraphics[width=0.31\linewidth]{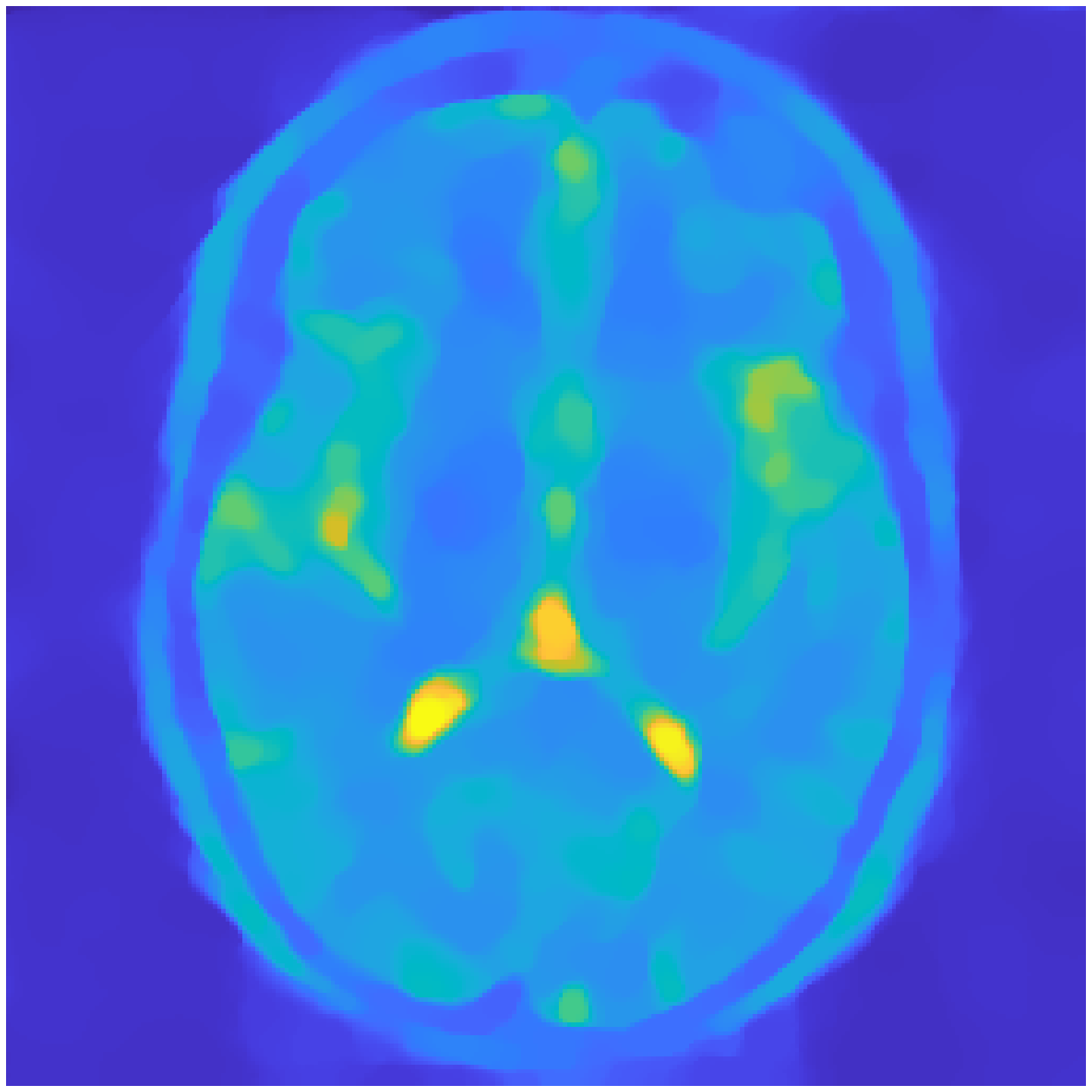} \\

    \includegraphics[width=0.31\linewidth]{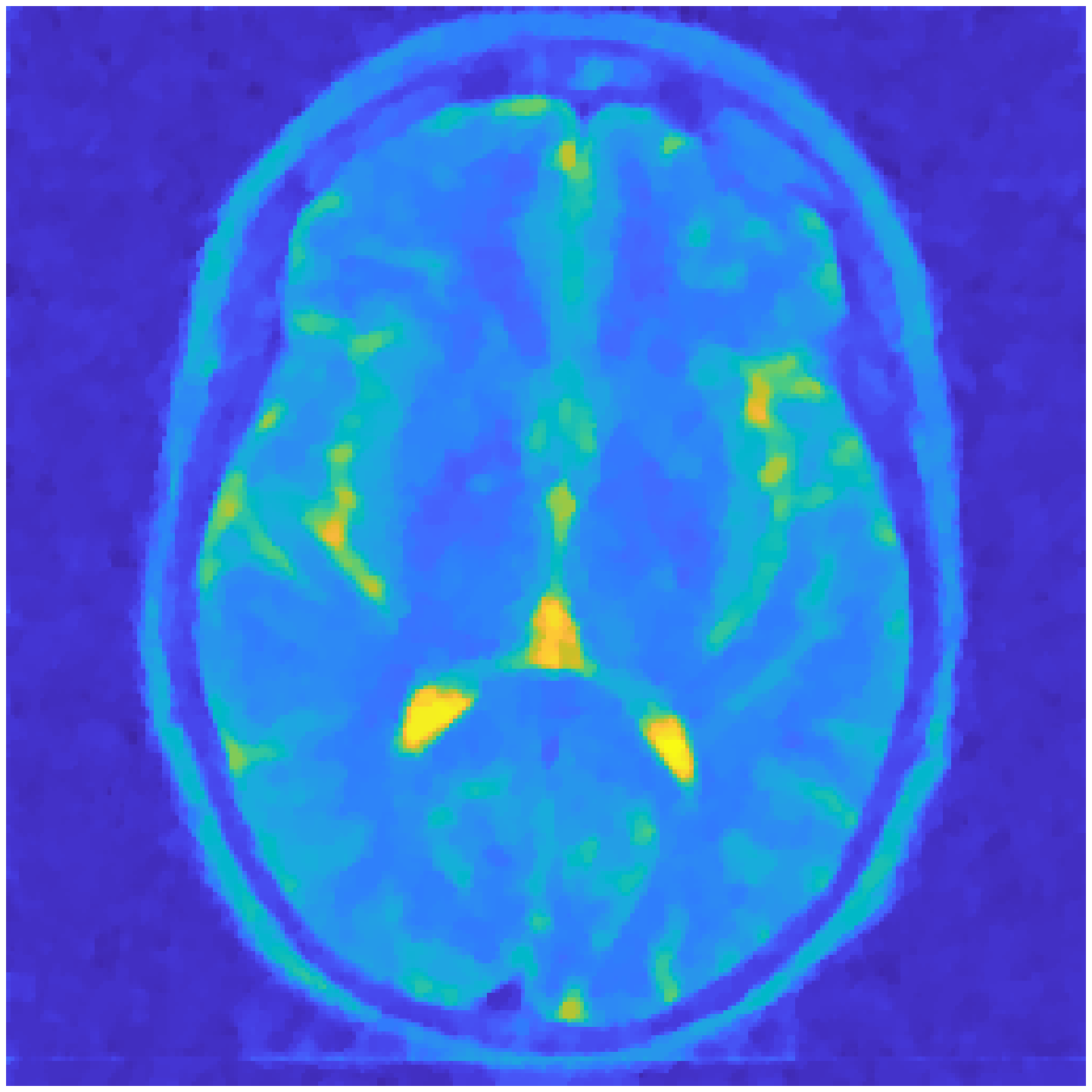} &
    \includegraphics[width=0.31\linewidth]{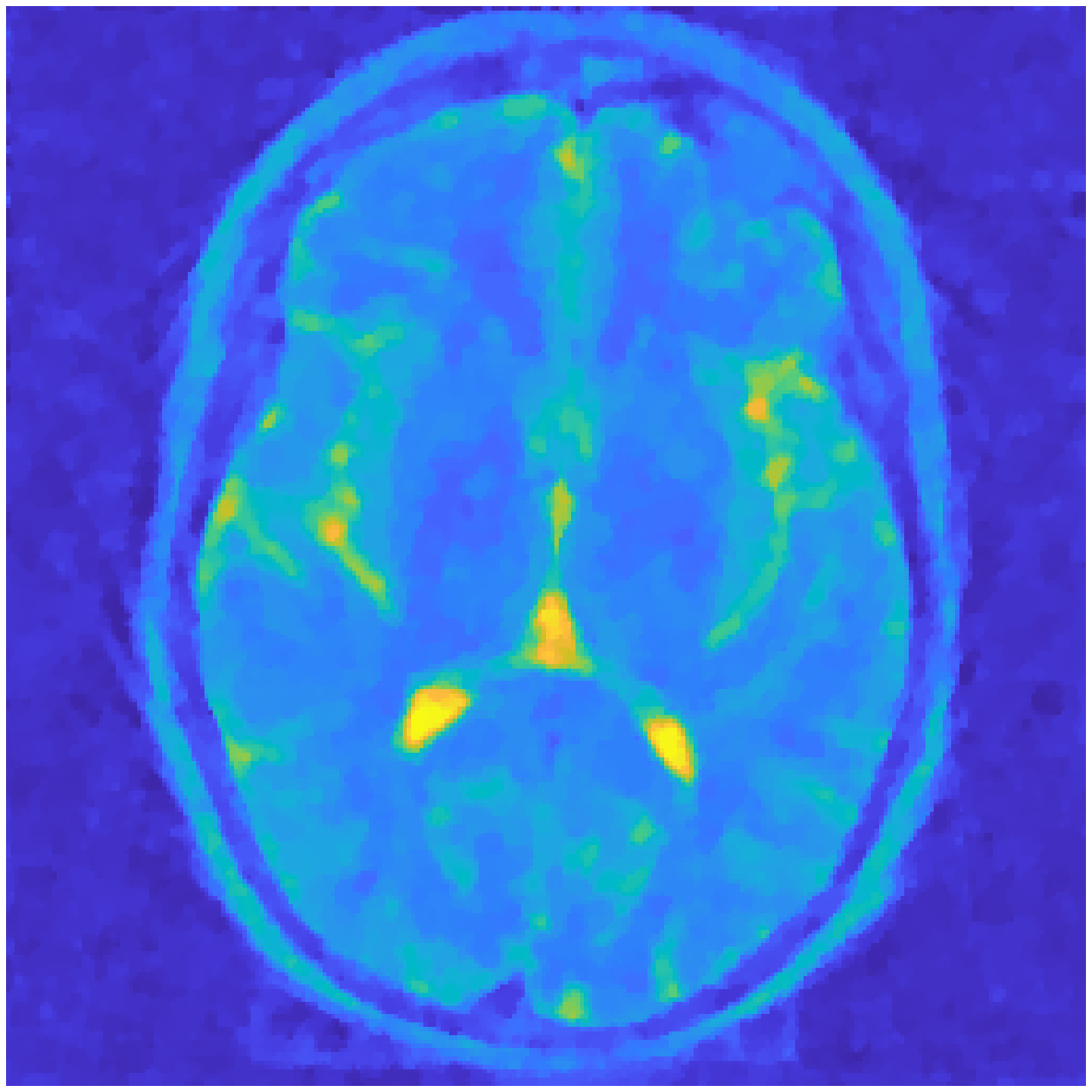} &
    \includegraphics[width=0.31\linewidth]{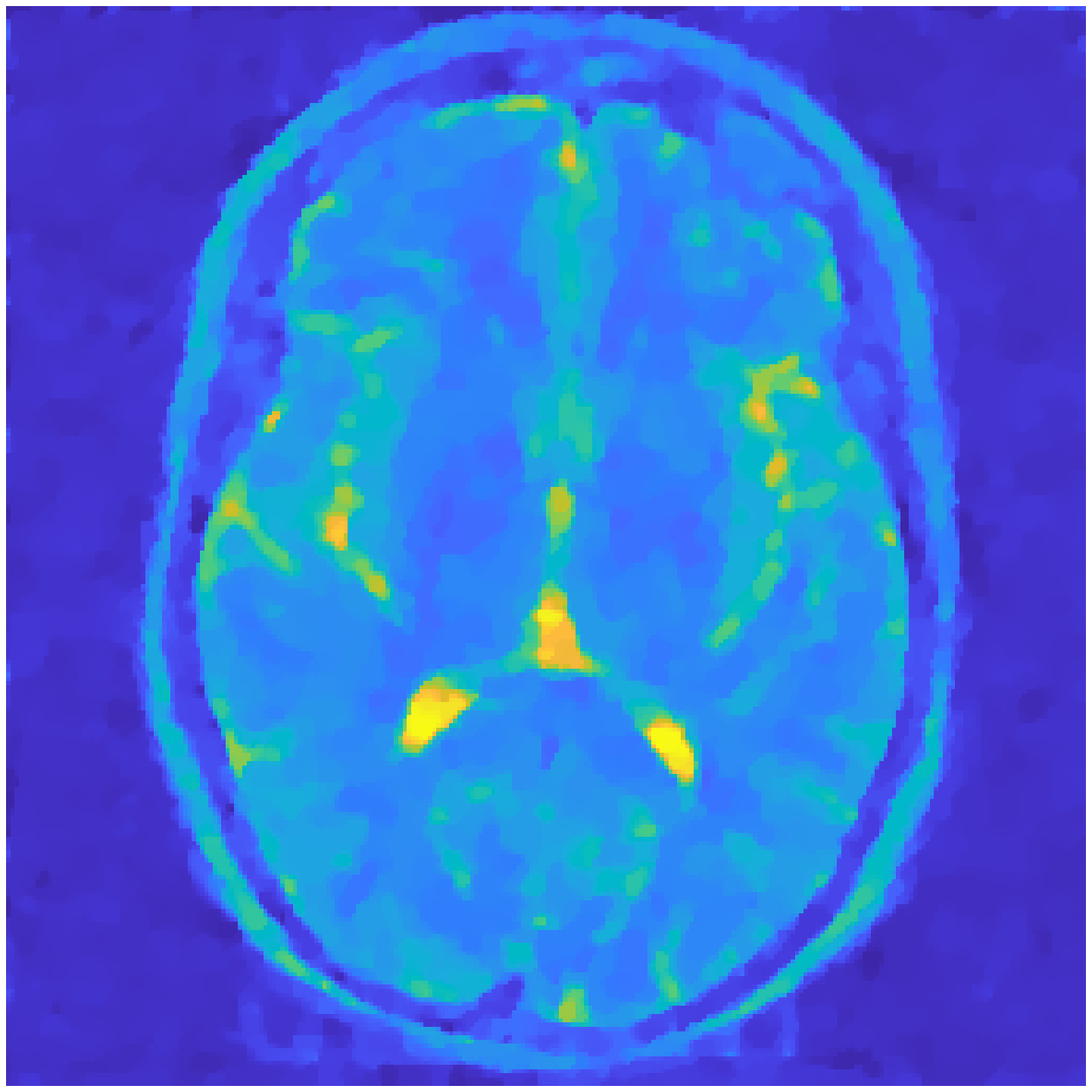} \\

    {\small $\Delta_{\max}=1,\, \sigma=0.01$} &
    {\small $\Delta_{\max}=3,\, \sigma=0.01$} &
    {\small $\Delta_{\max}=5,\, \sigma=0.01$} \\
    \end{tabular}

    \caption{256 by 256 brain. Top row: baseline; middle row: diminishing $\lambda_k$ approach; bottom row: new approach}
    \label{fig:braintype1recs1-001}
\end{figure}

\begin{figure}[!tb]
    \centering
    \setlength{\tabcolsep}{2pt}
    \begin{tabular}{cc}
         \includegraphics[width=0.48\linewidth]{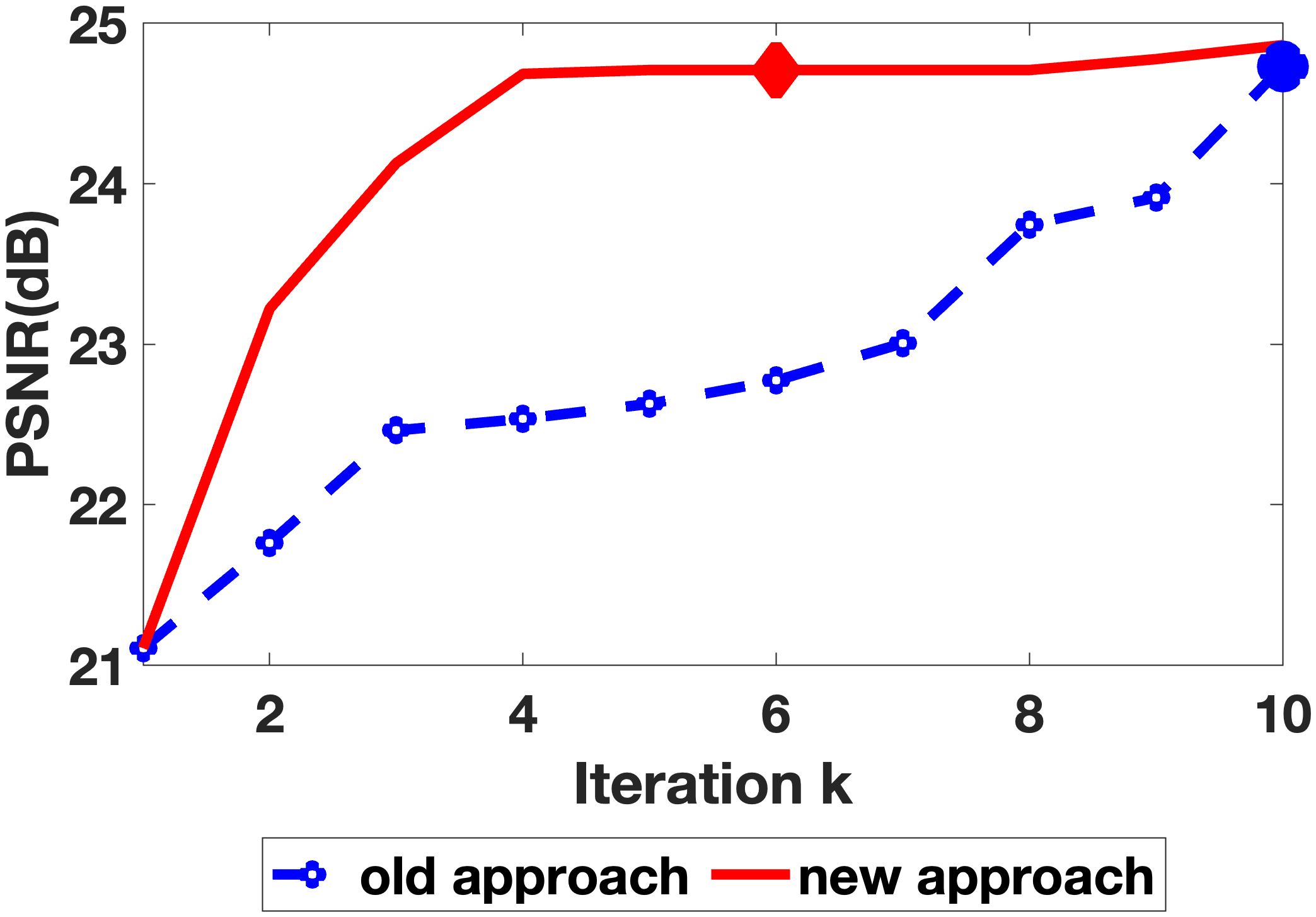} &
         \includegraphics[width=0.48\linewidth]{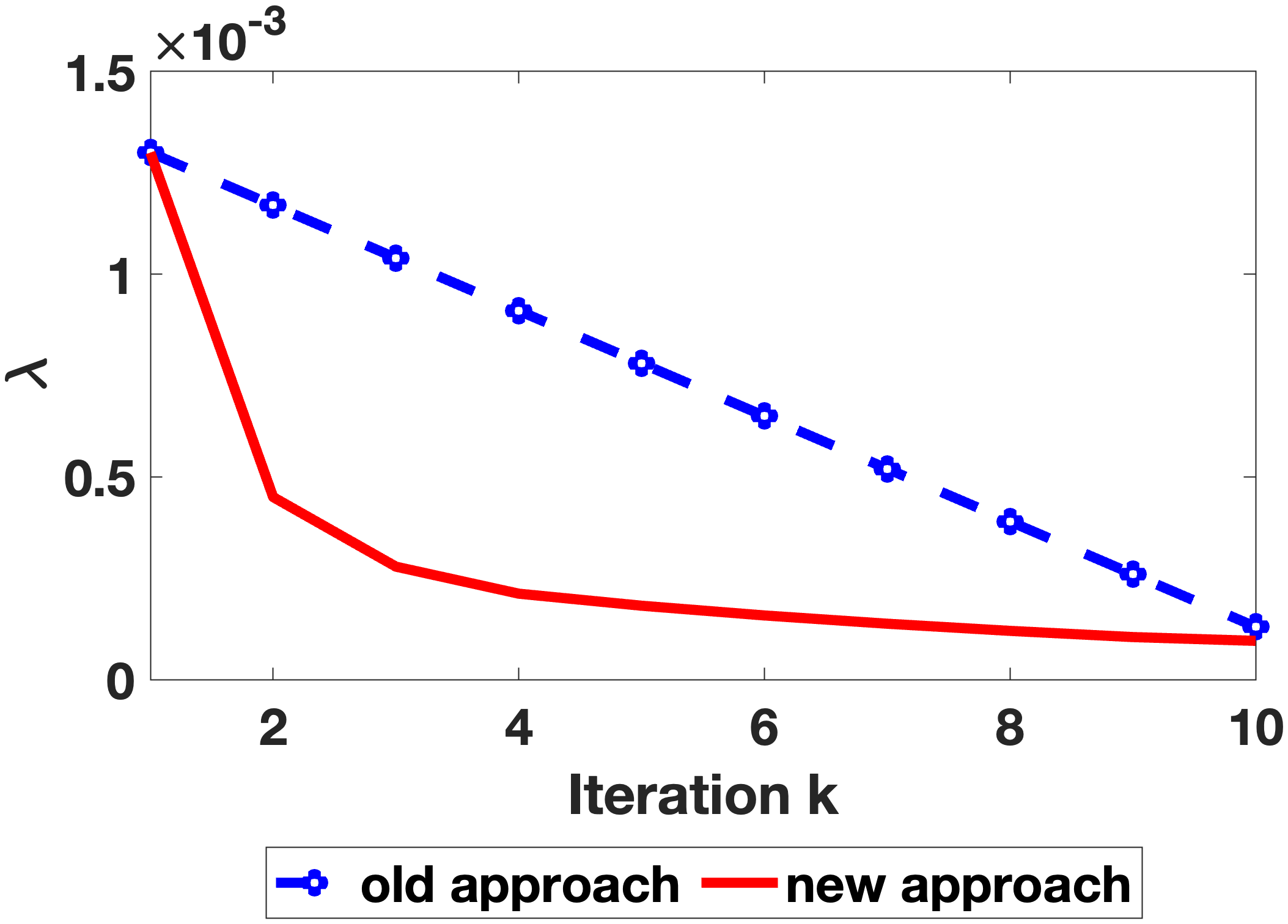} \\

        \includegraphics[width=0.48\linewidth]{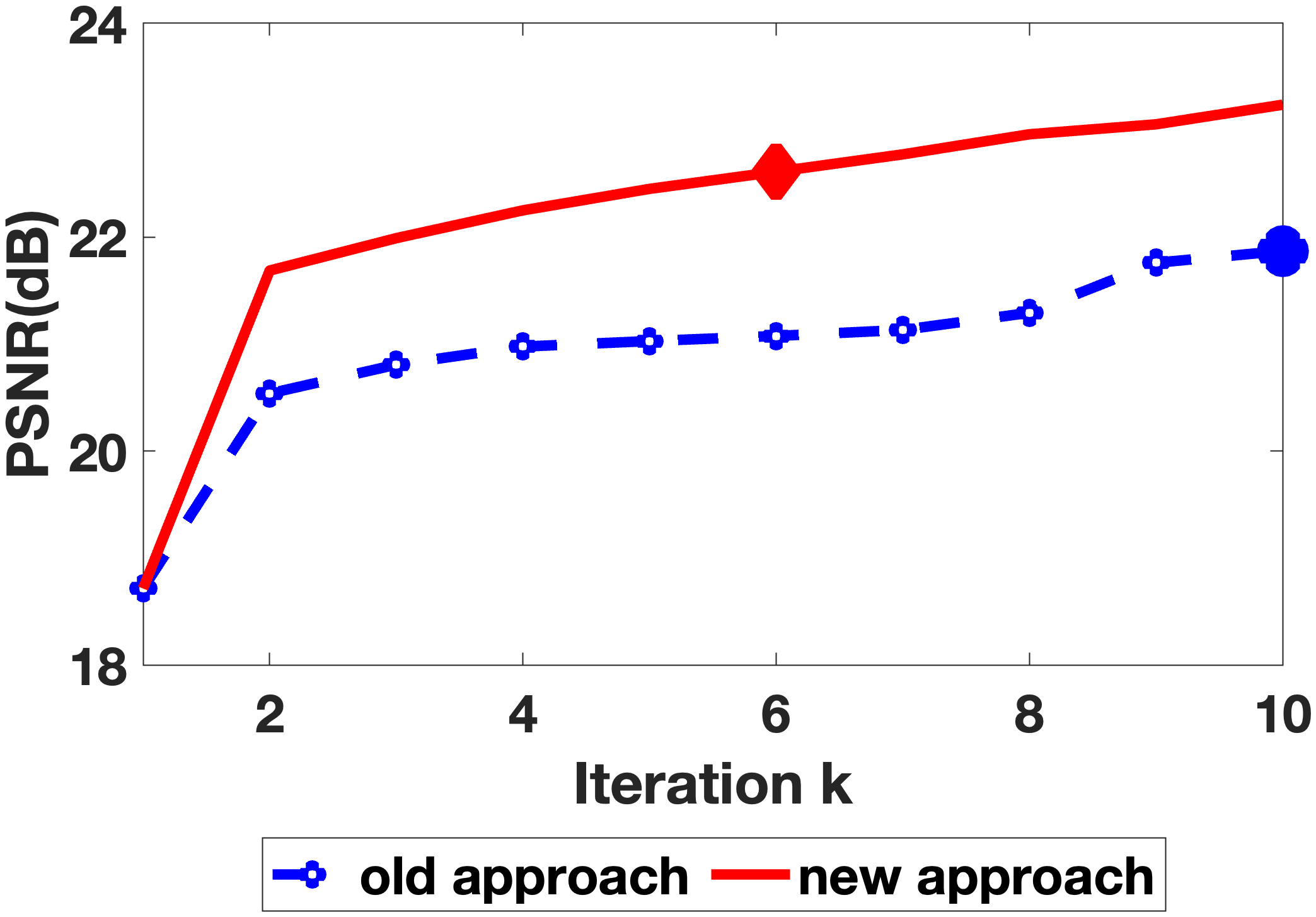} &
        \includegraphics[width=0.48\linewidth]{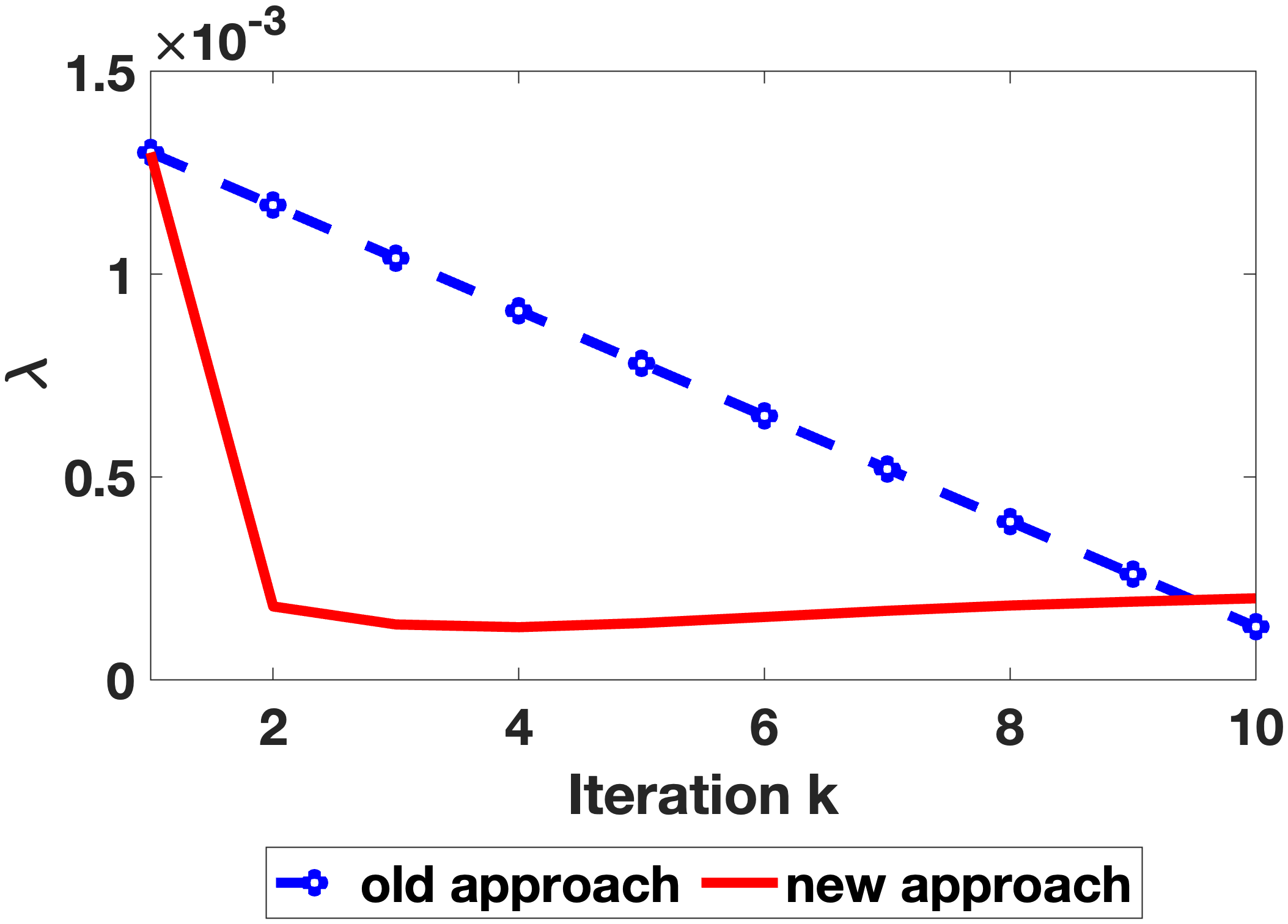} \\

        \includegraphics[width=0.48\linewidth]{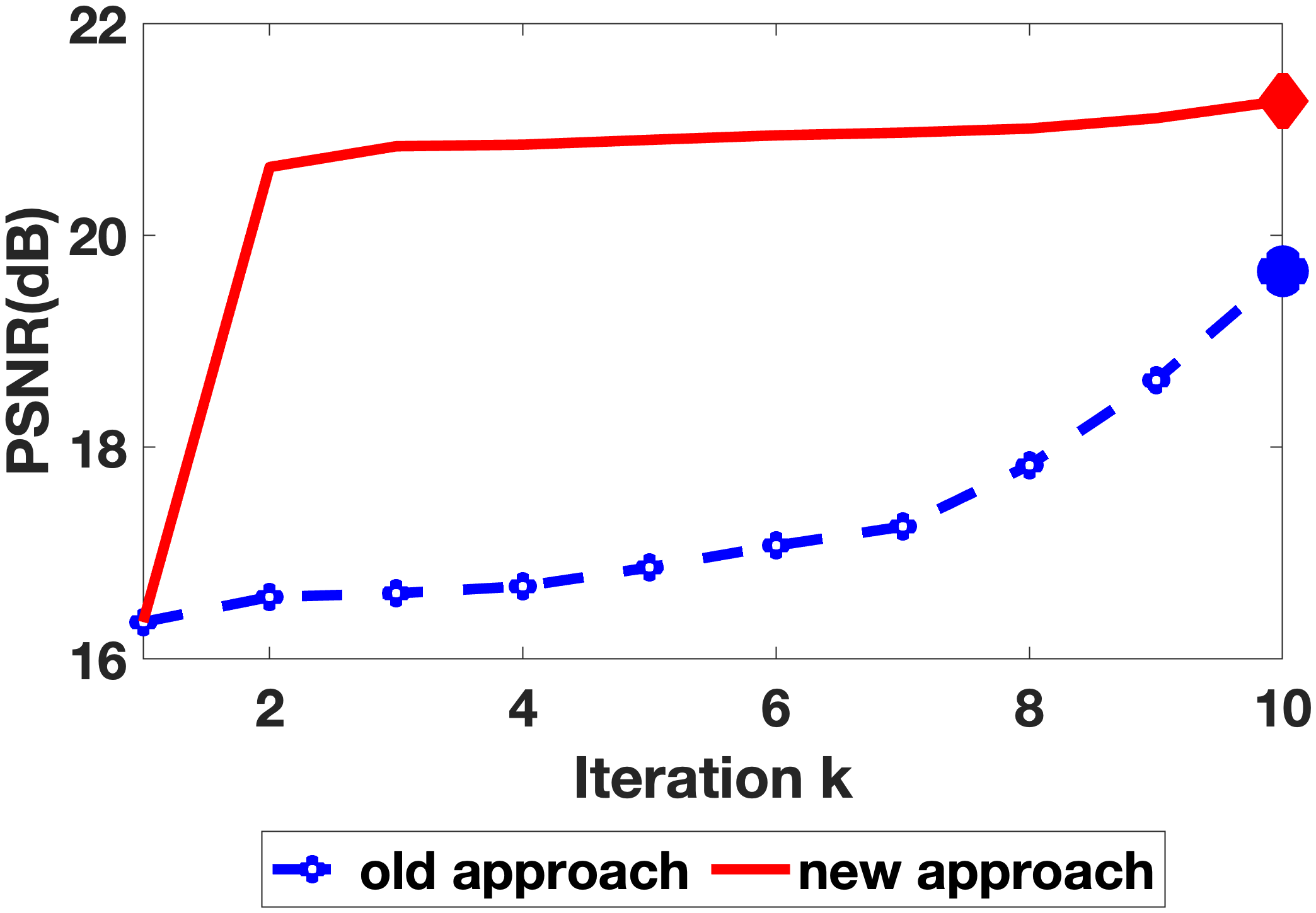} &
        \includegraphics[width=0.48\linewidth]{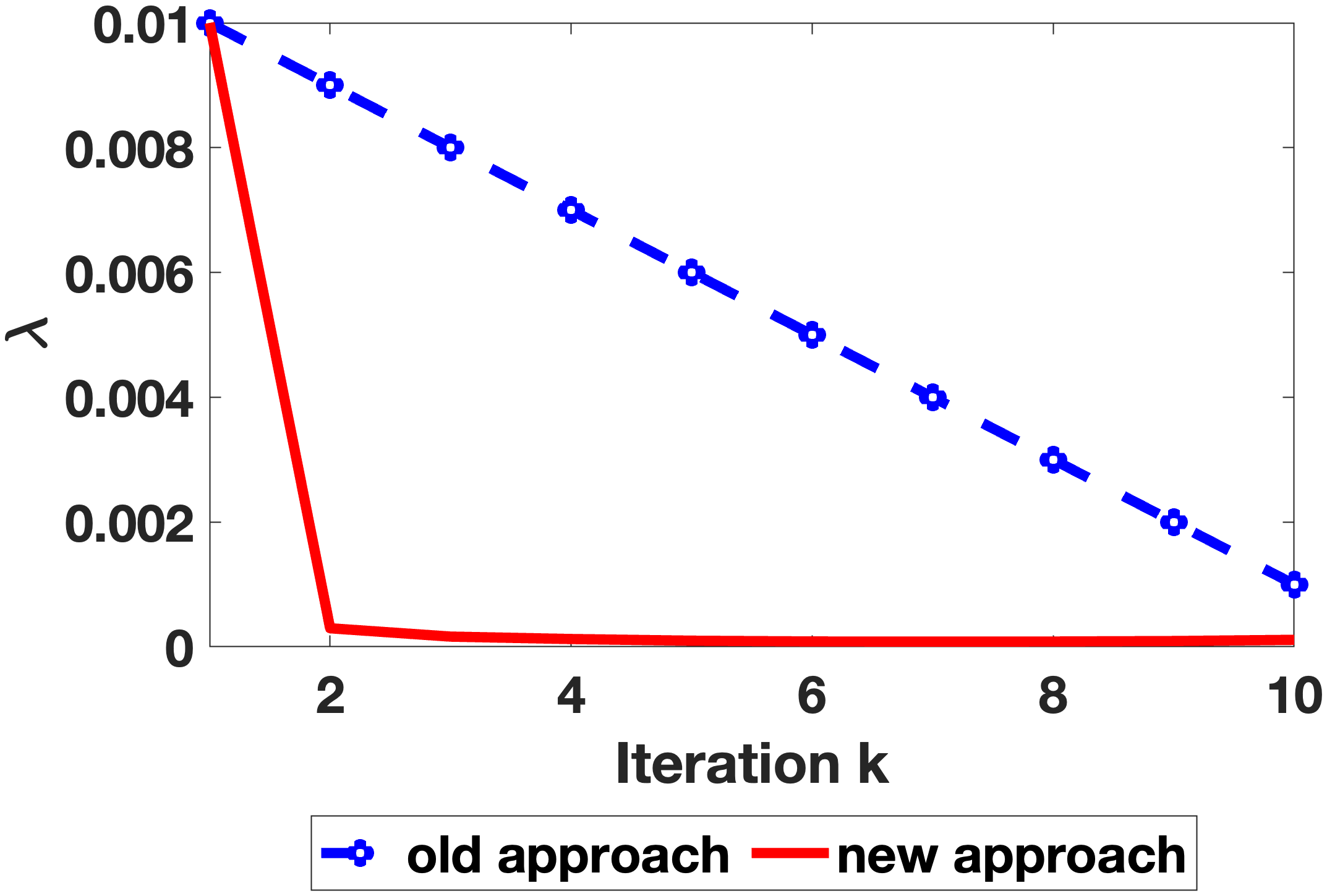} \\
    \end{tabular}
    \caption{256 by 256 phantom. Top row: maxDrift = 1, noise = 0.02; middle row: maxDrift = 3, noise = 0.02; bottom row: maxDrift = 5, noise = 0.02. Left column: PSNR; right column: $\lambda_k$'s}
    \label{fig:phantype1conv1-002}
\end{figure}

\begin{figure}[!tb]
    \centering
    \setlength{\tabcolsep}{2pt}
    \begin{tabular}{cc}
        \includegraphics[width=0.48\linewidth]{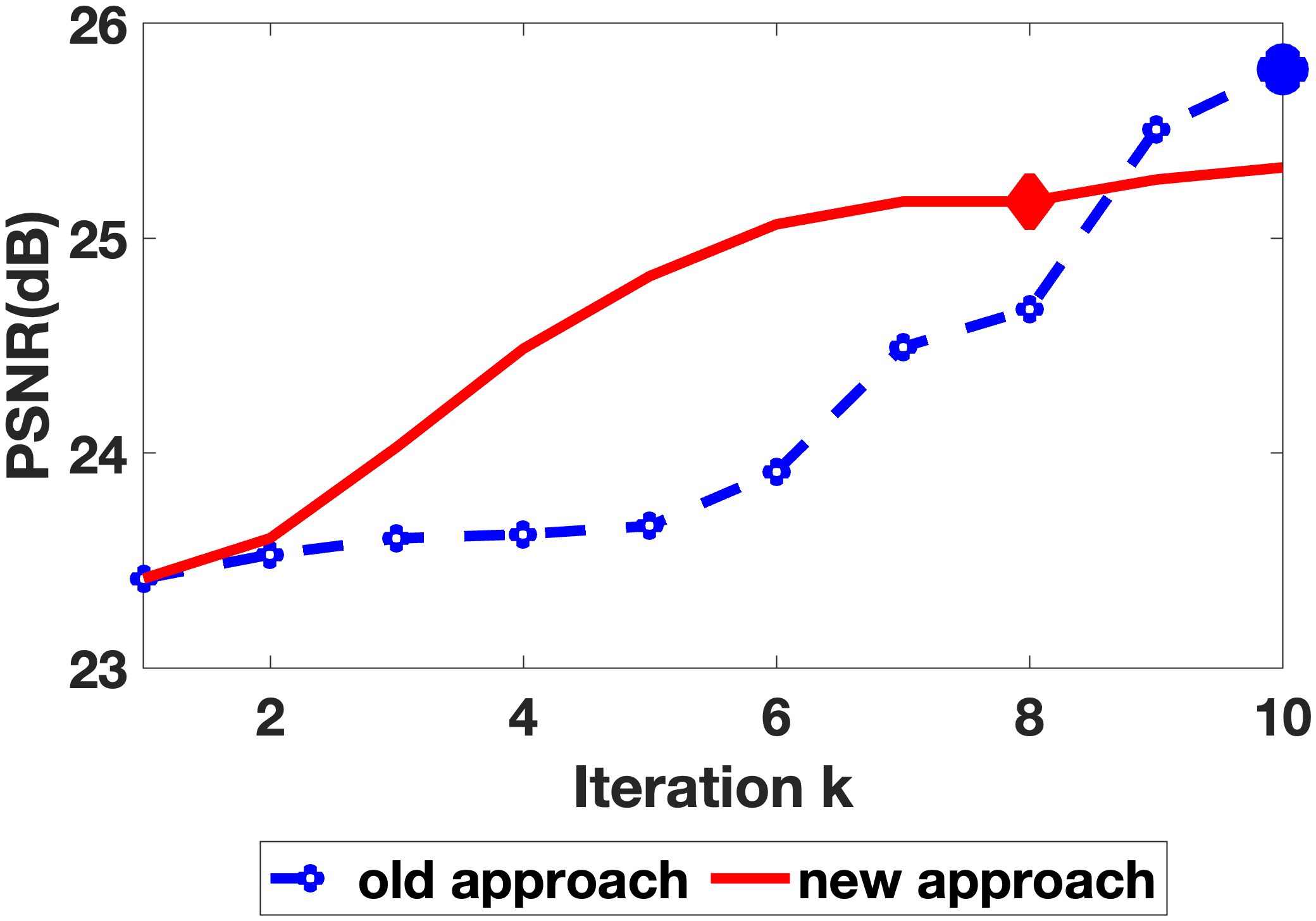} &
         \includegraphics[width=0.48\linewidth]{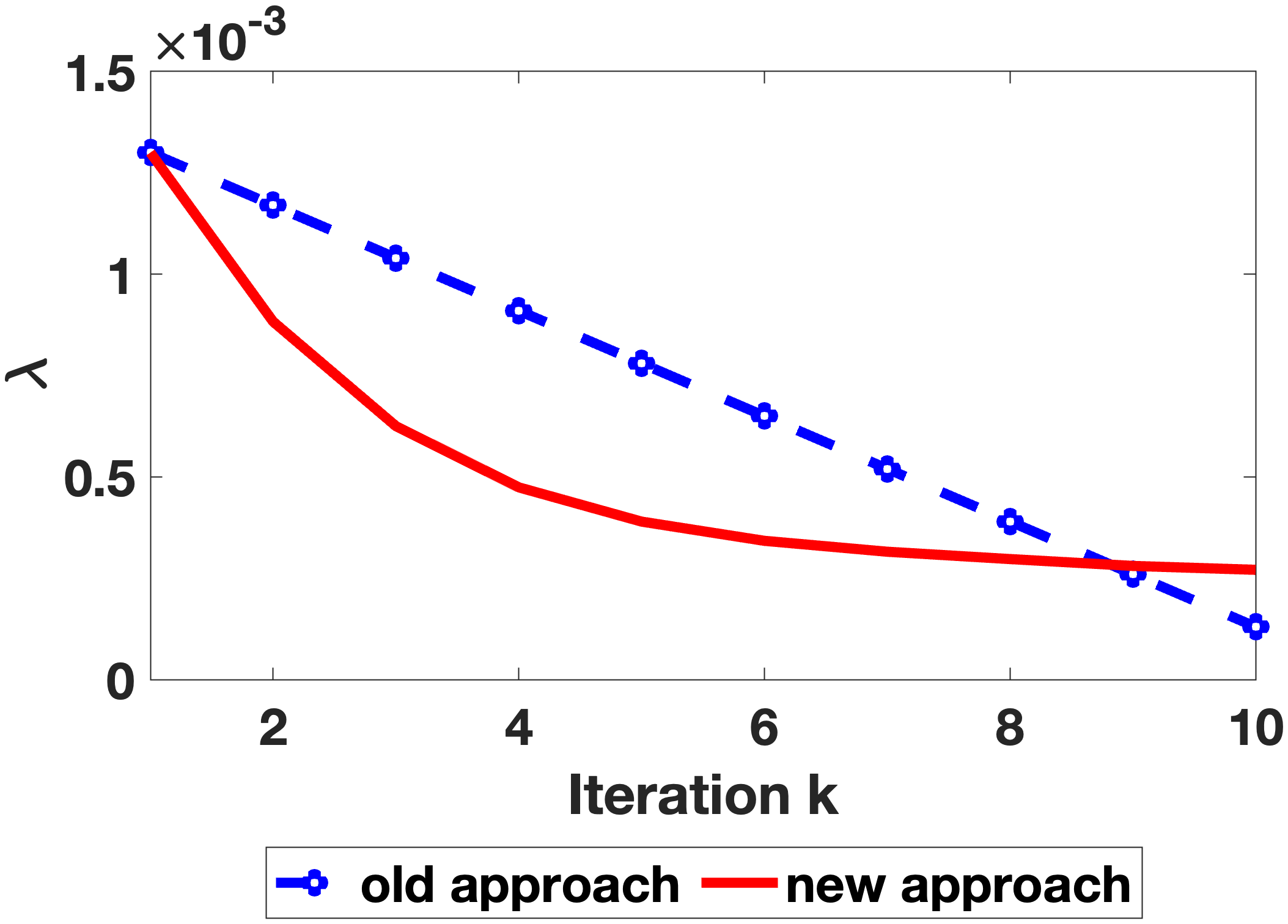} \\

        \includegraphics[width=0.48\linewidth]{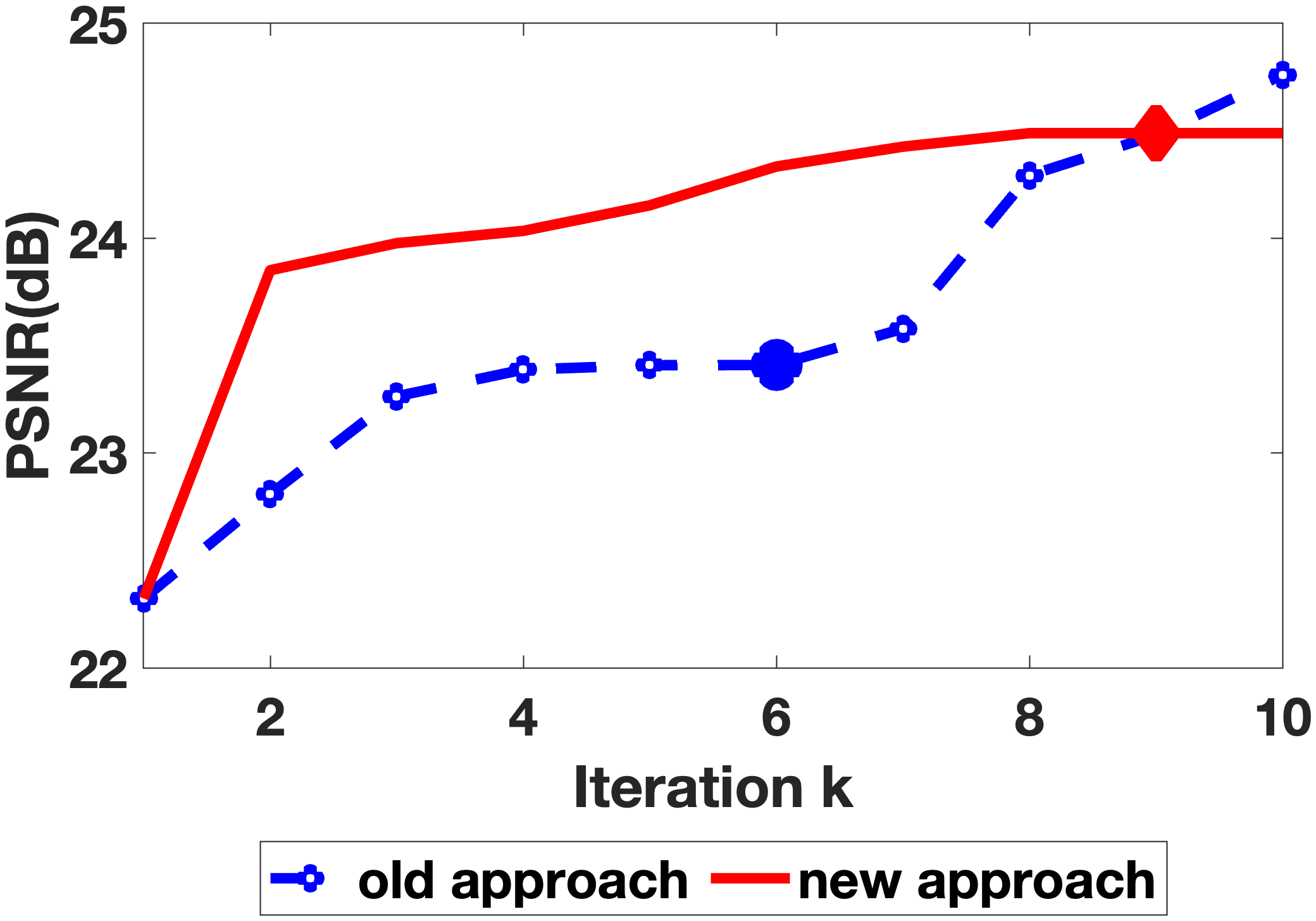} &
        \includegraphics[width=0.48\linewidth]{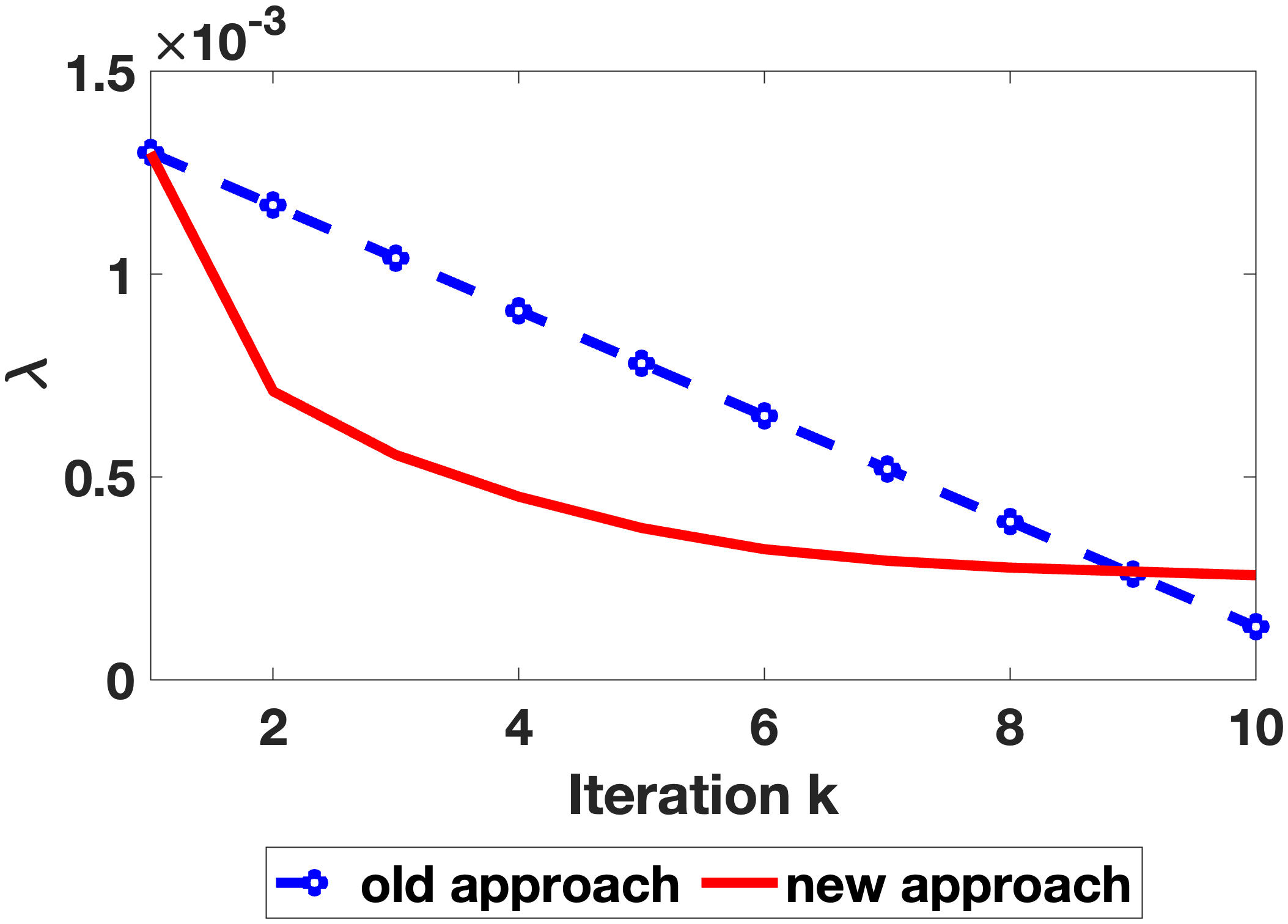} \\

        \includegraphics[width=0.48\linewidth]{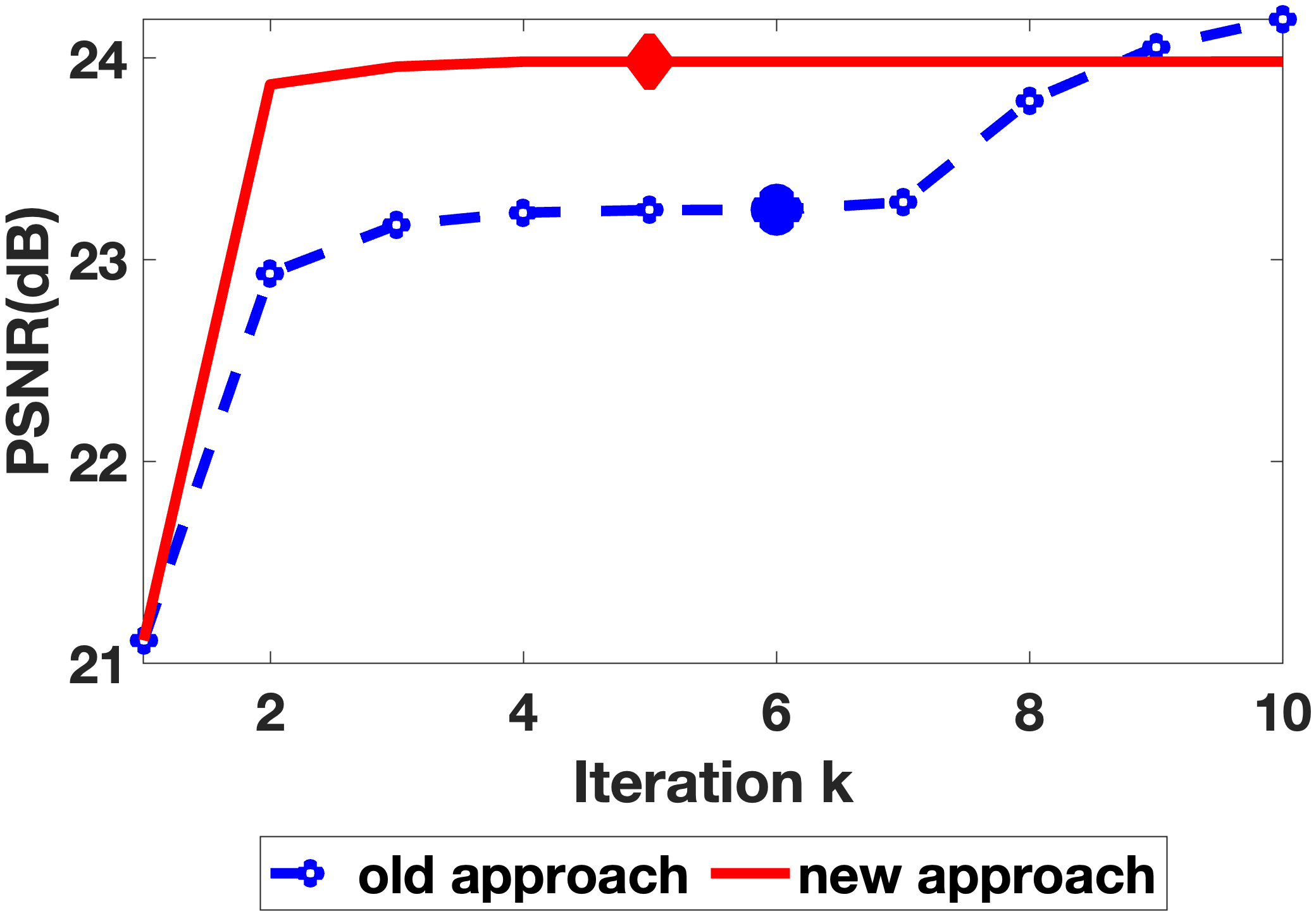} &
        \includegraphics[width=0.48\linewidth]{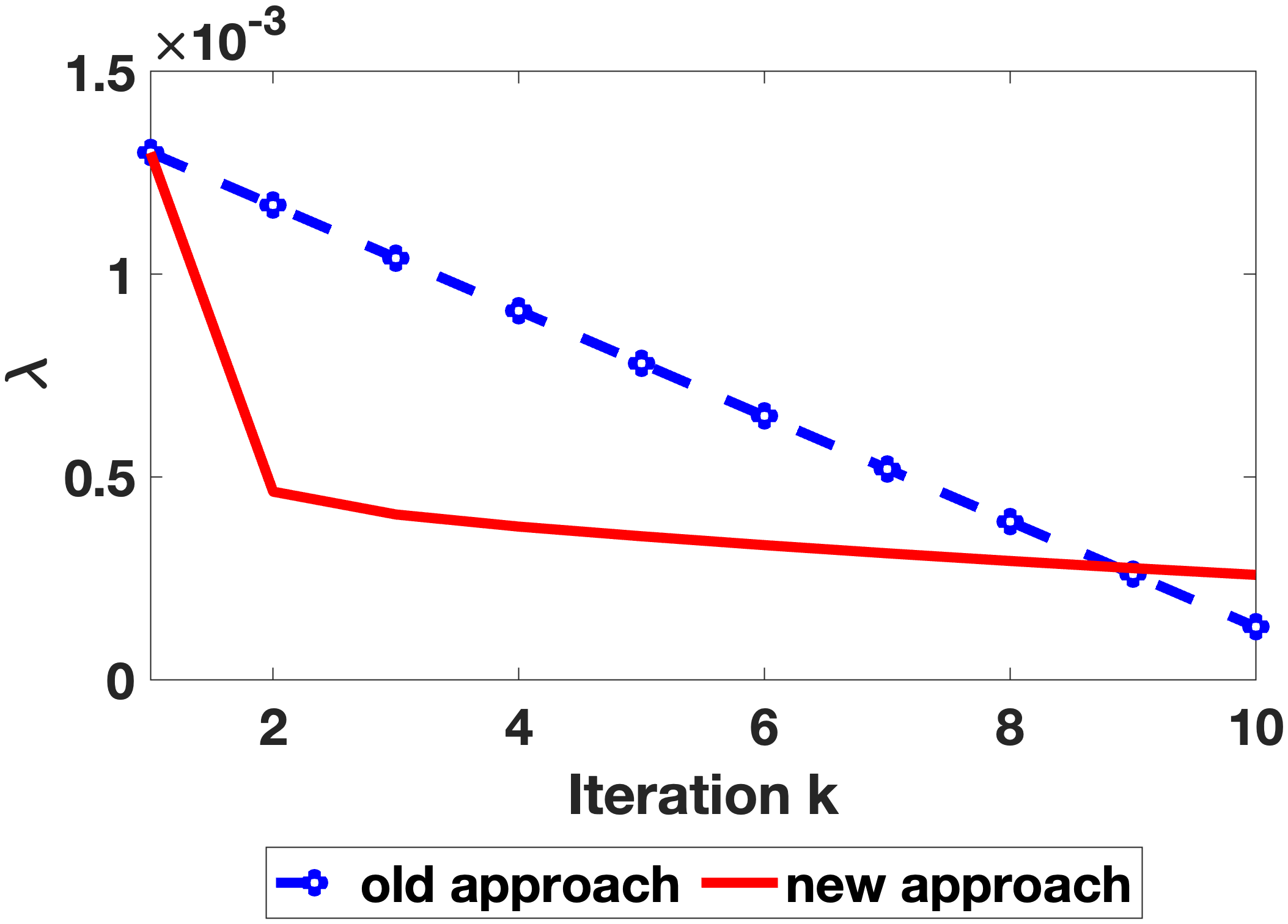} \\
    \end{tabular}
    \caption{256 by 256 brain. Top row: maxDrift = 1, noise = 0.02; middle row: maxDrift = 3, noise = 0.02; bottom row: maxDrift = 5, noise = 0.02. Left column: PSNR; right column: $\lambda_k$'s}
    \label{fig:braintype1conv1-002}
\end{figure}

We visualize the reconstructions of the phantom in Fig.~\ref{fig:phantype1recs1-002} and the brain image in Fig.~\ref{fig:braintype1recs1-001} for different drift levels. Drift calibration recovers sharper features than the baseline reconstruction, and the adaptive $\lambda_k$ update often improves contrast relative to the diminishing $\lambda_k$ update. We further investigate the iterative behavior of the two $\lambda_k$ update rules in Fig.~\ref{fig:phantype1conv1-002} and Fig.~\ref{fig:braintype1conv1-002}. The adaptive rule changes $\lambda_k$ until it stabilizes, while the empirical diminishing rule follows the prescribed schedule from \cite{huang2019calibrating}. As mentioned earlier, for a fair comparison, we allow both algorithms to continue until reaching the maximum number of iterations. We indicate the stopping point of the adaptive approach with a red diamond according to criterion~\eqref{eqn:stop}.

\subsection{Type II Drift Result}

We now evaluate the method on the more challenging Type II drift, where each individual measurement possesses a unique drift value, $\delta_{\theta,\tau}$. To simulate this general drift scenario, we use the pixel intensities of the ``cameraman'' image to create a 2D drift pattern, which is then applied to the sinogram of the test objects. This setup introduces a complex, non-uniform error field. The experimental configuration, including noise levels and maximum drift magnitudes, remains consistent with the Type I experiments. The quantitative results are summarized in Table~\ref{tab:PSNR_SSIM_Type2}. Both calibrated approaches improve over the uncalibrated baseline in many settings, and the adaptive $\lambda_k$ approach is strongest for several brain cases, but the Type II problem remains sensitive to the drift regularization weight.

\begin{figure}[!tb]
    \centering
    \setlength{\tabcolsep}{2pt}
    \begin{tabular}{ccc}
    \includegraphics[width=0.31\linewidth]{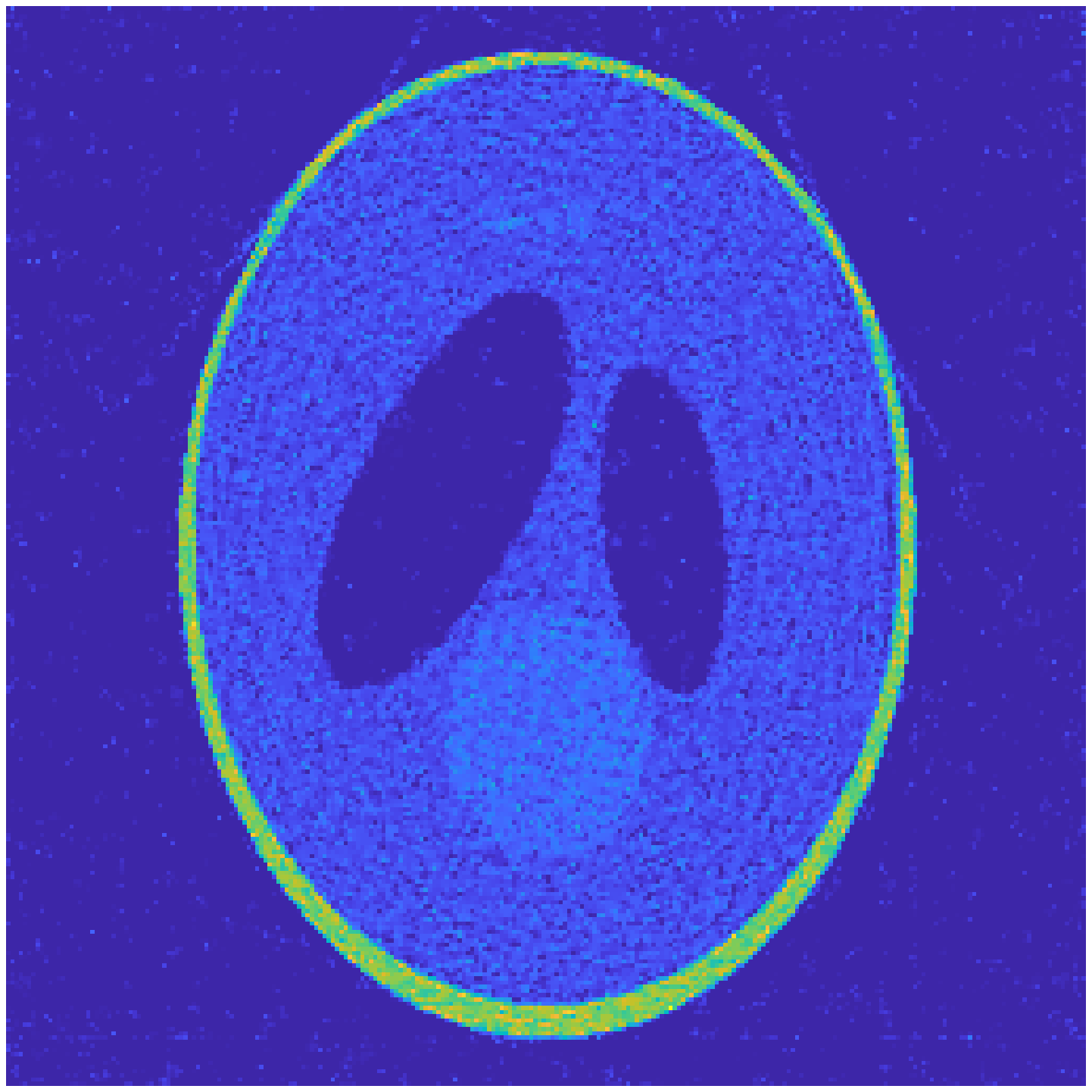} &
    \includegraphics[width=0.31\linewidth]{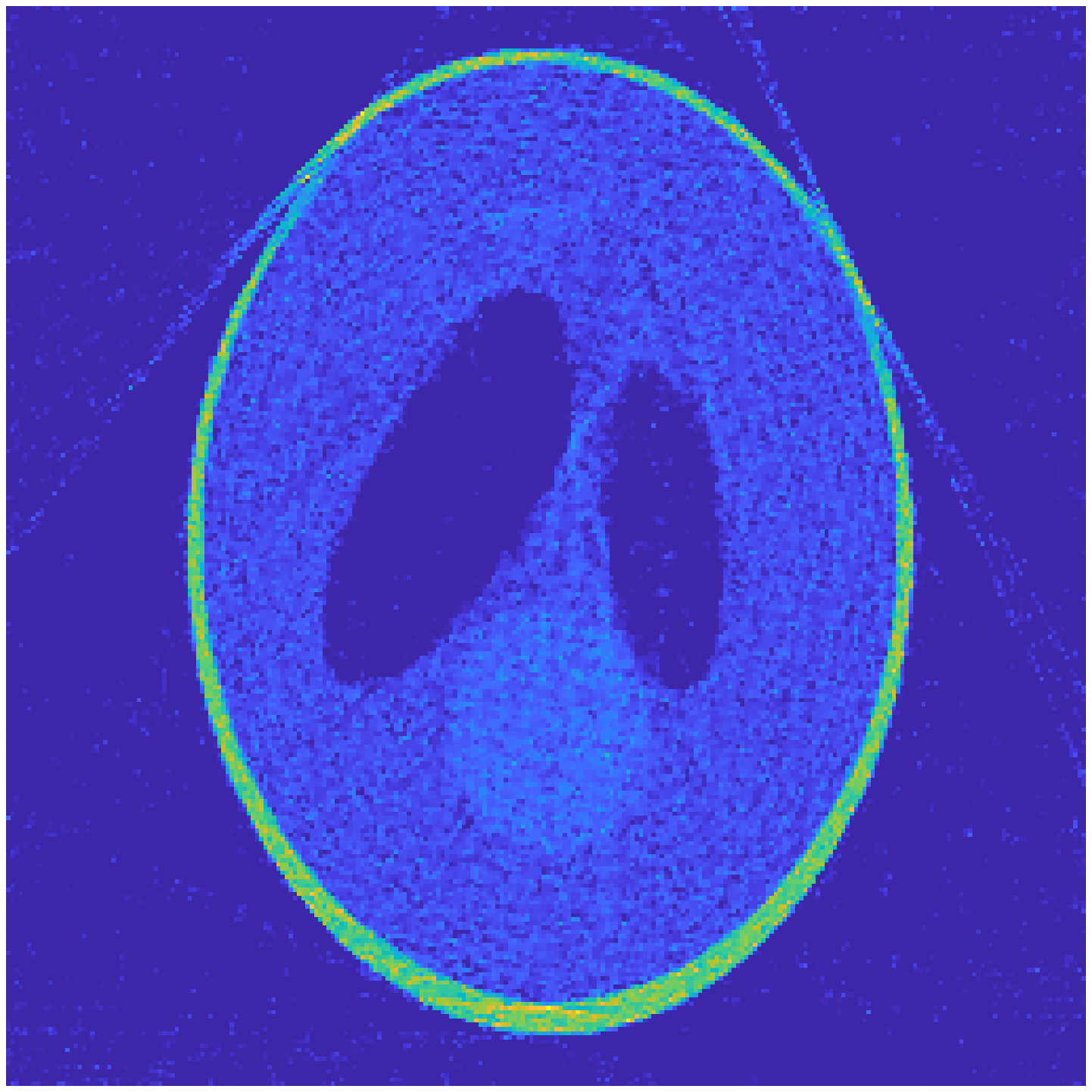} &
    \includegraphics[width=0.31\linewidth]{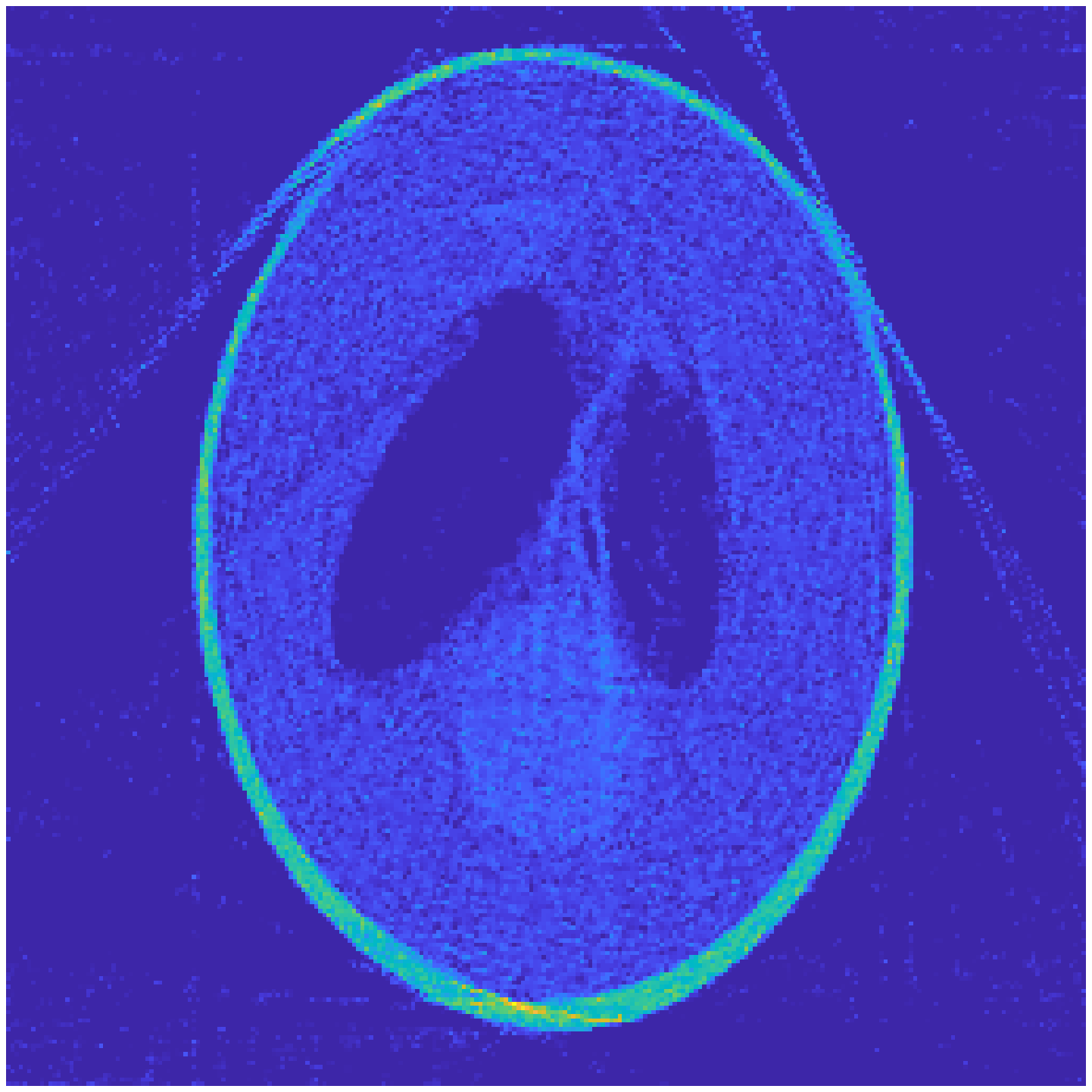} \\

    \includegraphics[width=0.31\linewidth]{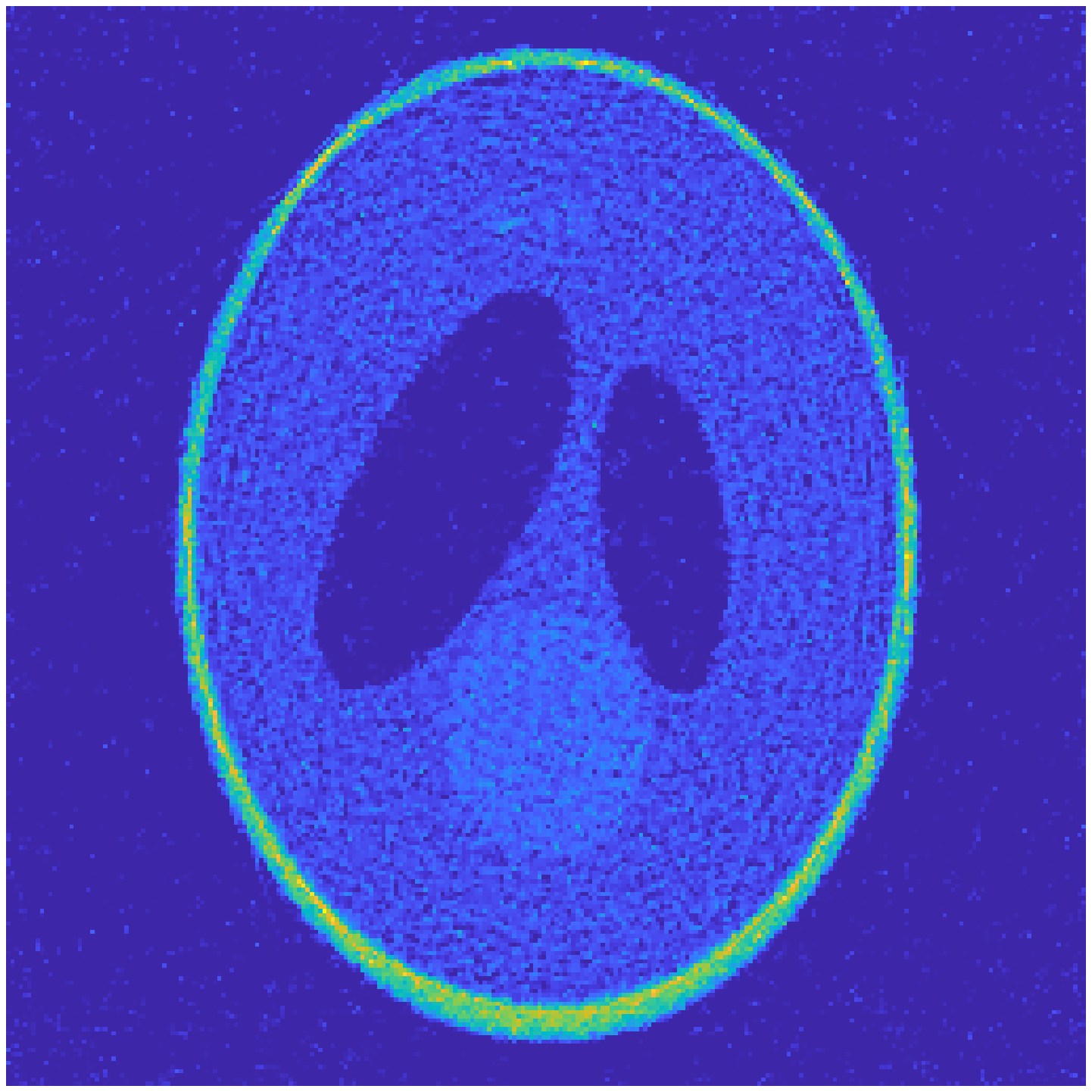} &
    \includegraphics[width=0.31\linewidth]{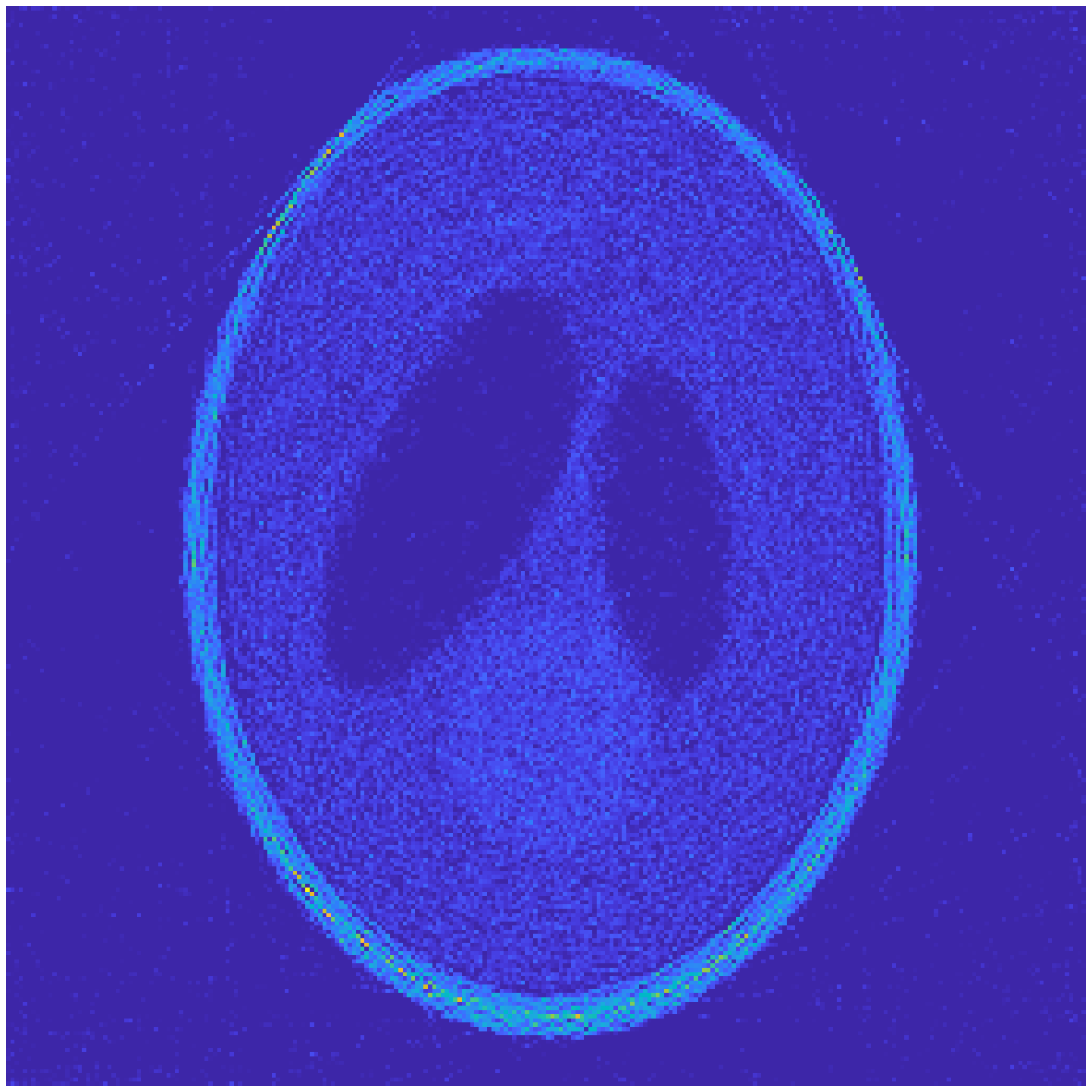} &
    \includegraphics[width=0.31\linewidth]{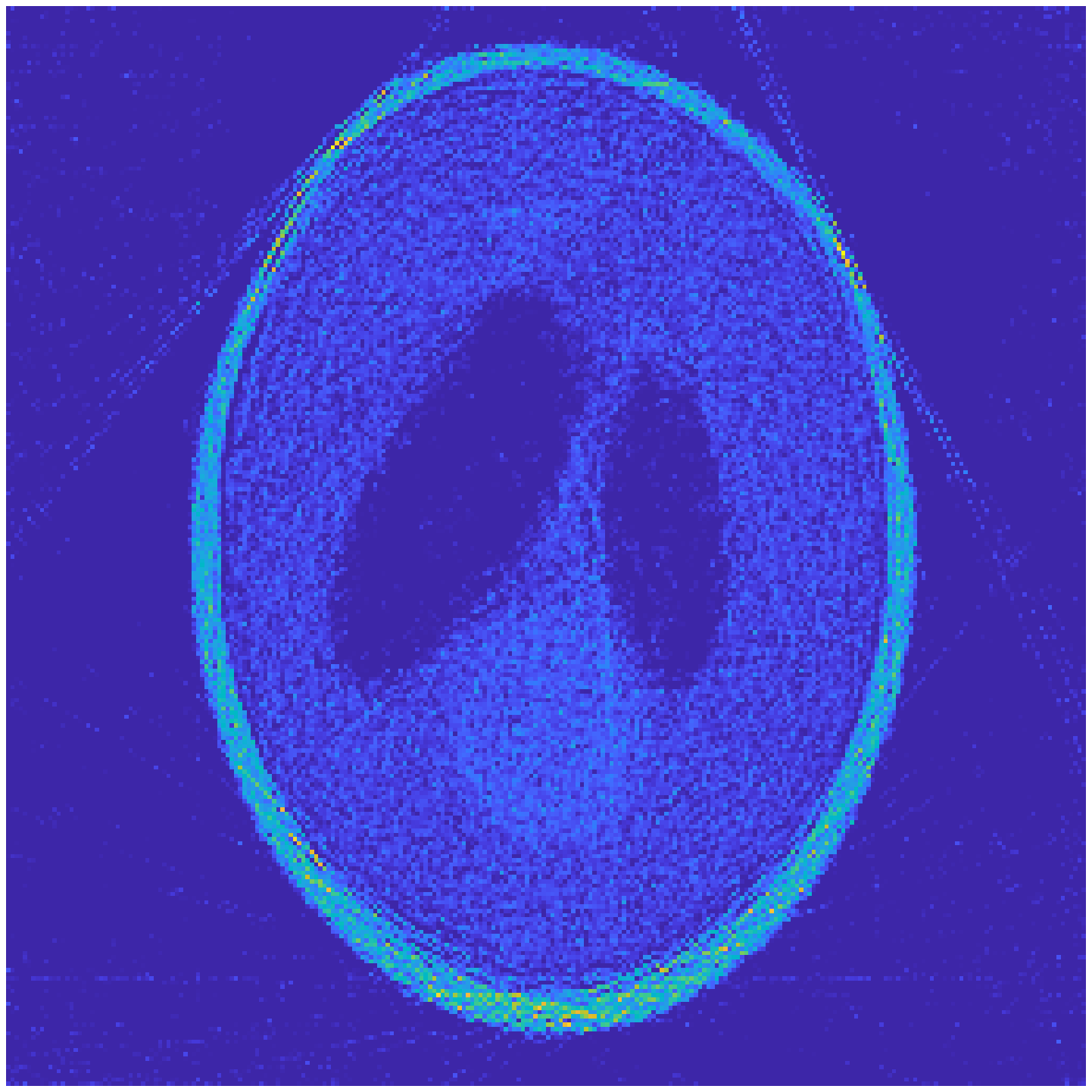} \\

    \includegraphics[width=0.31\linewidth]{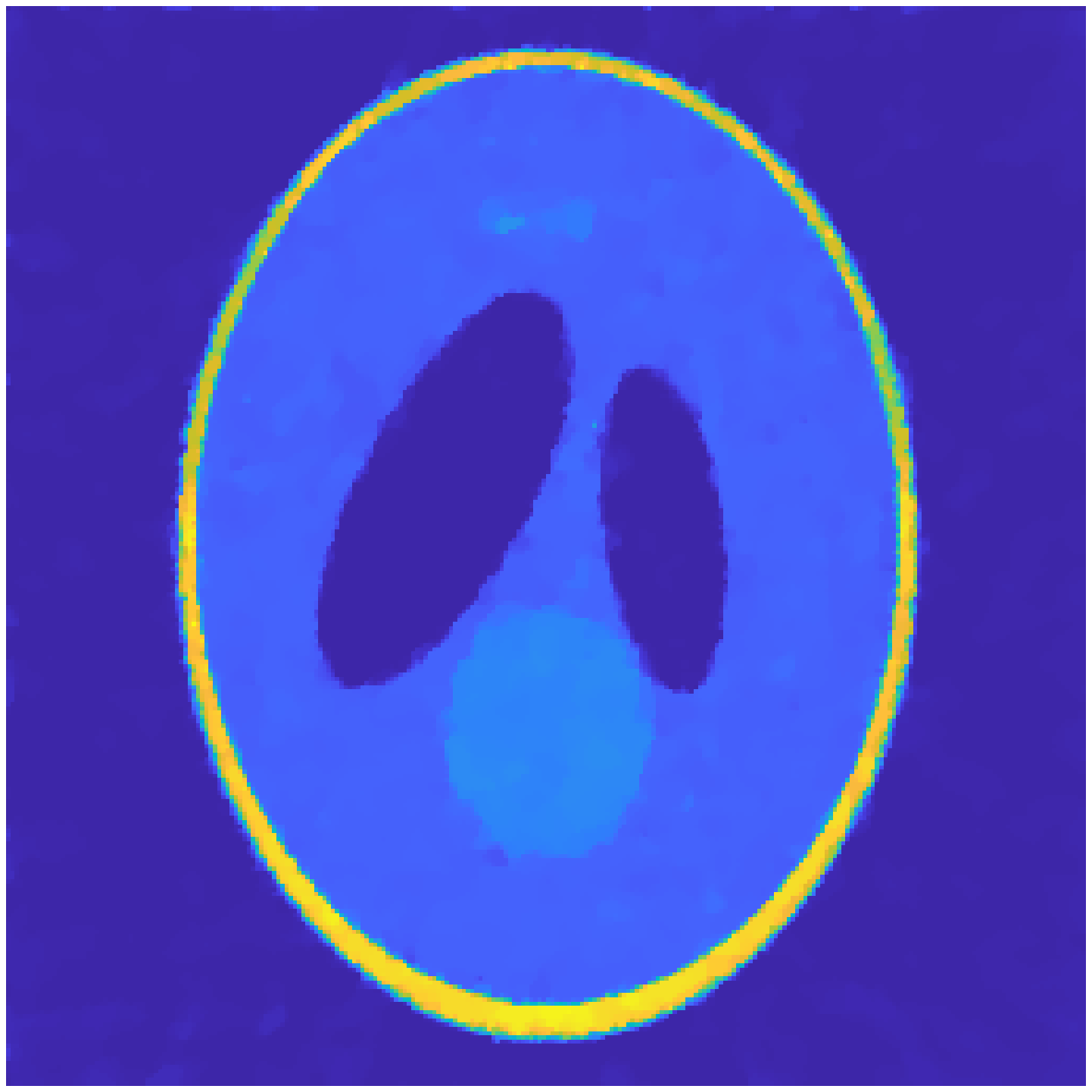} &
    \includegraphics[width=0.31\linewidth]{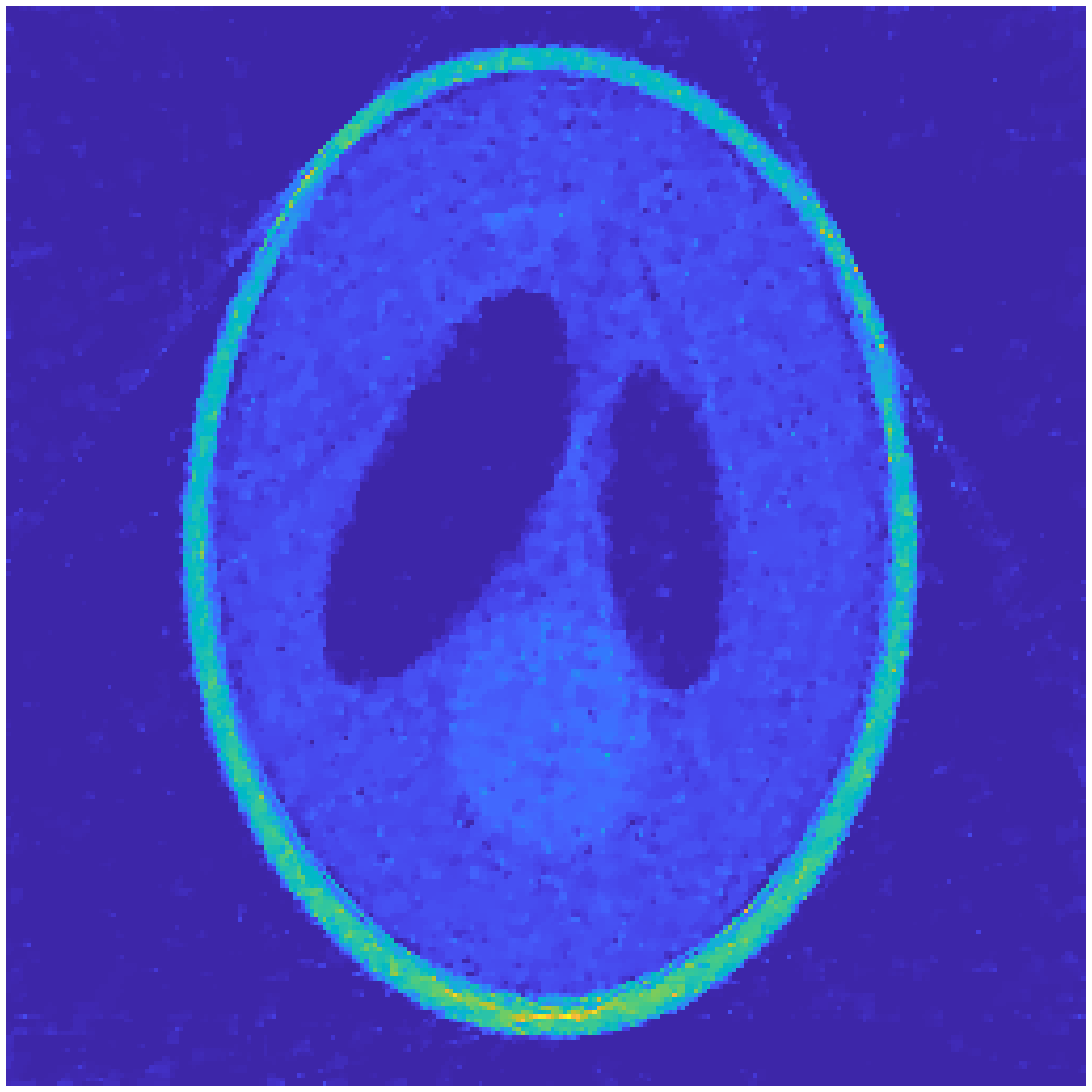} &
    \includegraphics[width=0.31\linewidth]{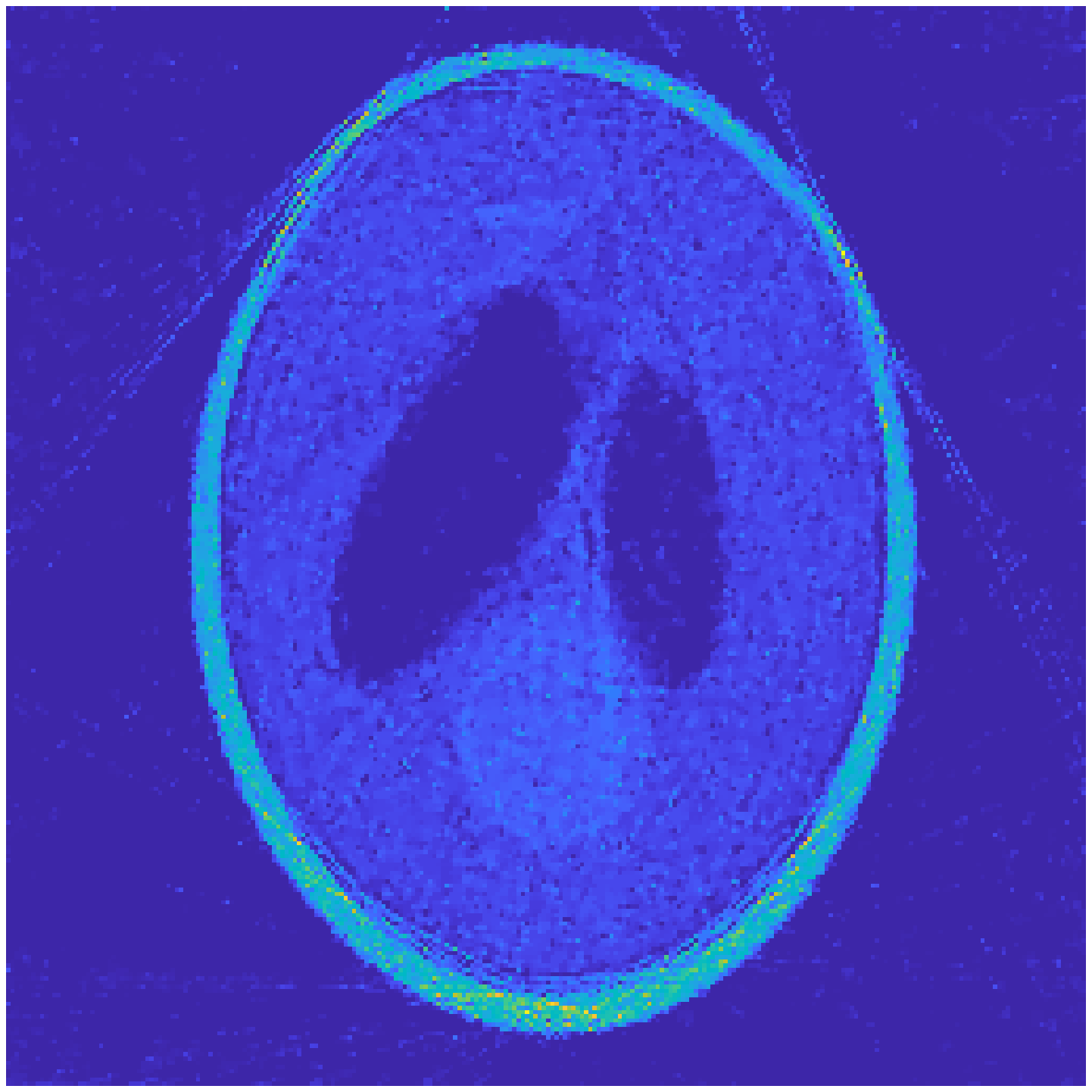} \\

    {\small $\Delta_{\max}=1,\, \sigma=0.02$} &
    {\small $\Delta_{\max}=3,\, \sigma=0.02$} &
    {\small $\Delta_{\max}=5,\, \sigma=0.02$} \\
    \end{tabular}

    \caption{256 by 256 phantom. Top row: baseline; middle row: diminishing $\lambda_k$ approach; bottom row: adaptive $\lambda_k$ approach}
    \label{fig:phantype2recs1-002}
\end{figure}

\begin{figure}[!tb]
    \centering
    \setlength{\tabcolsep}{2pt}
    \begin{tabular}{ccc}
    \includegraphics[width=0.31\linewidth]{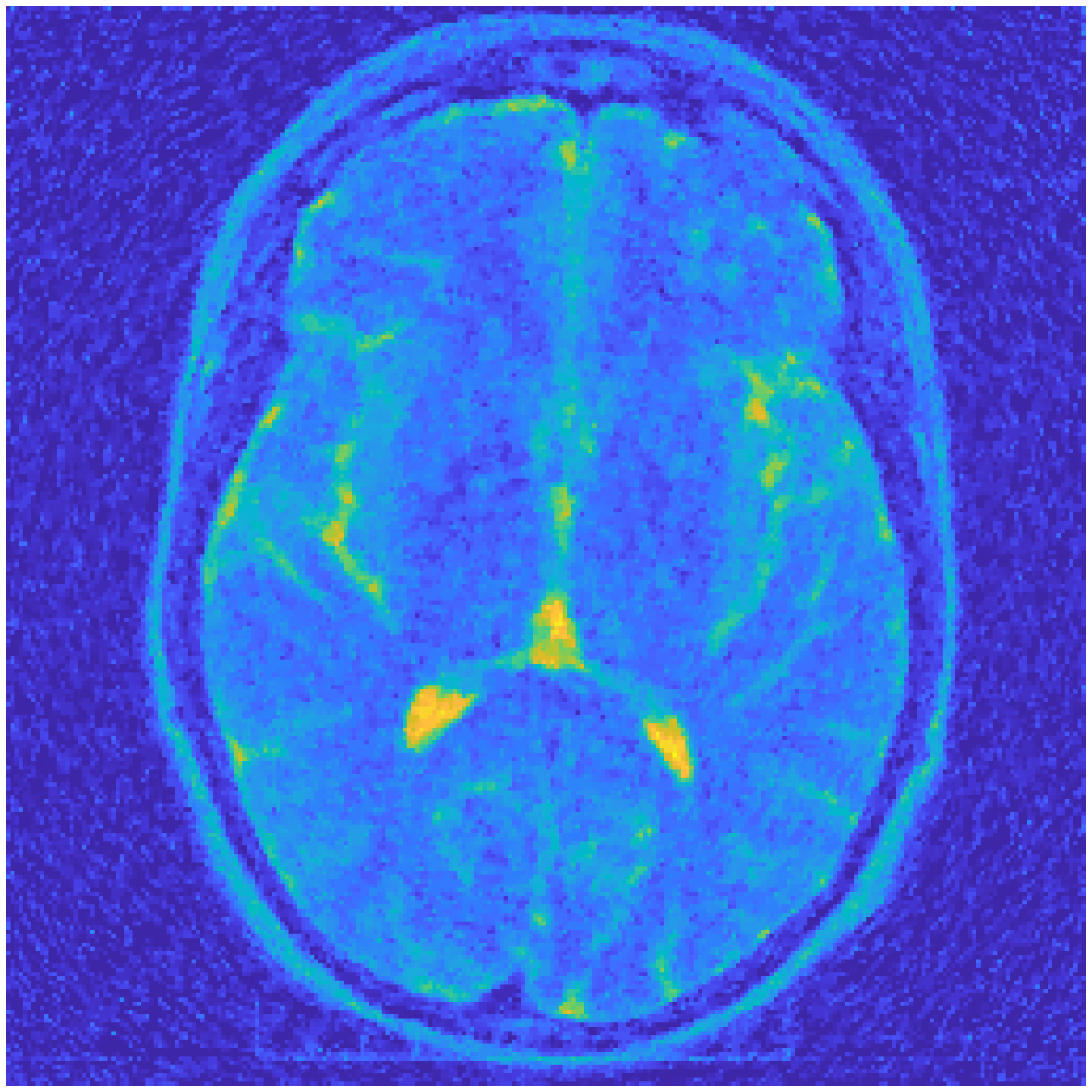} &
    \includegraphics[width=0.31\linewidth]{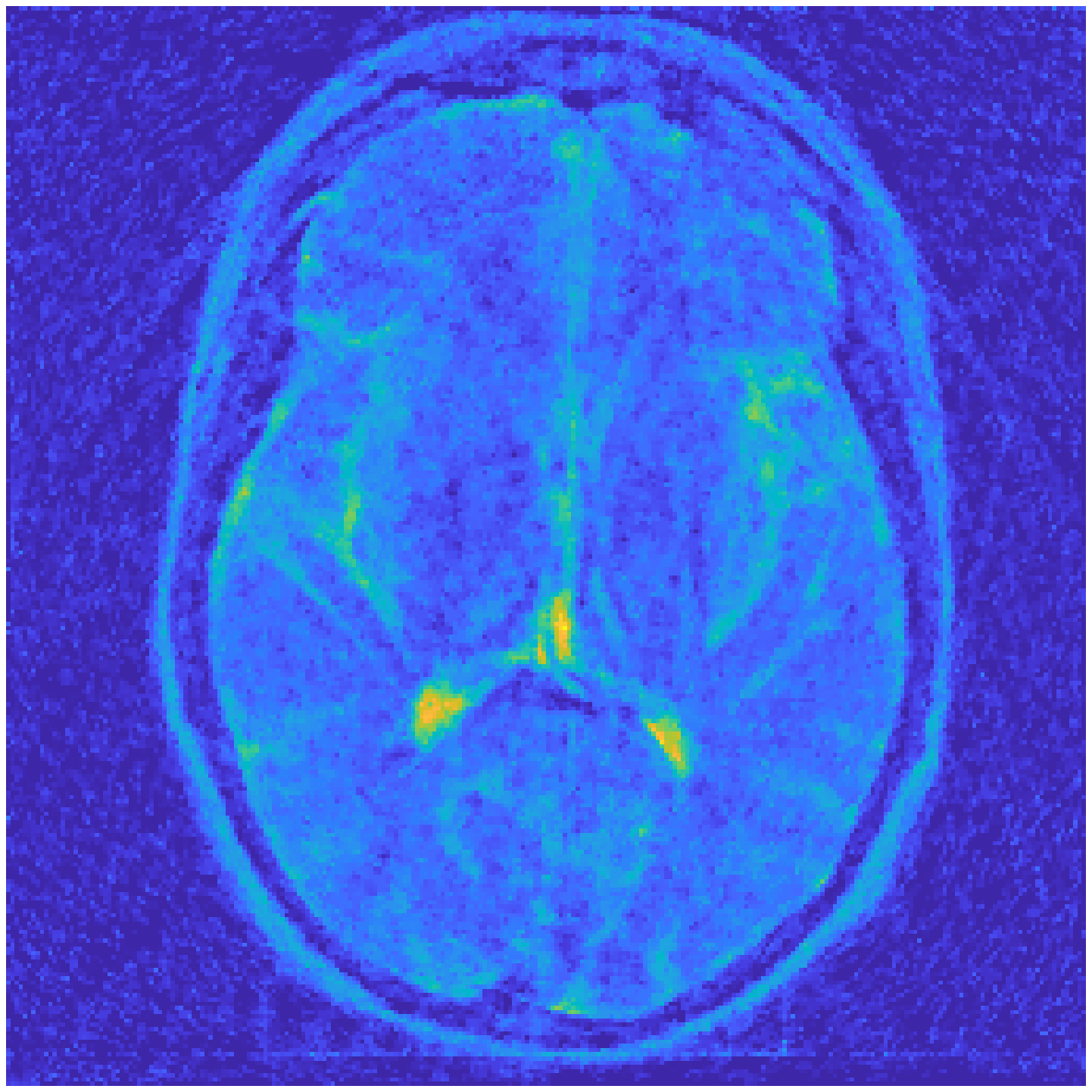} &
    \includegraphics[width=0.31\linewidth]{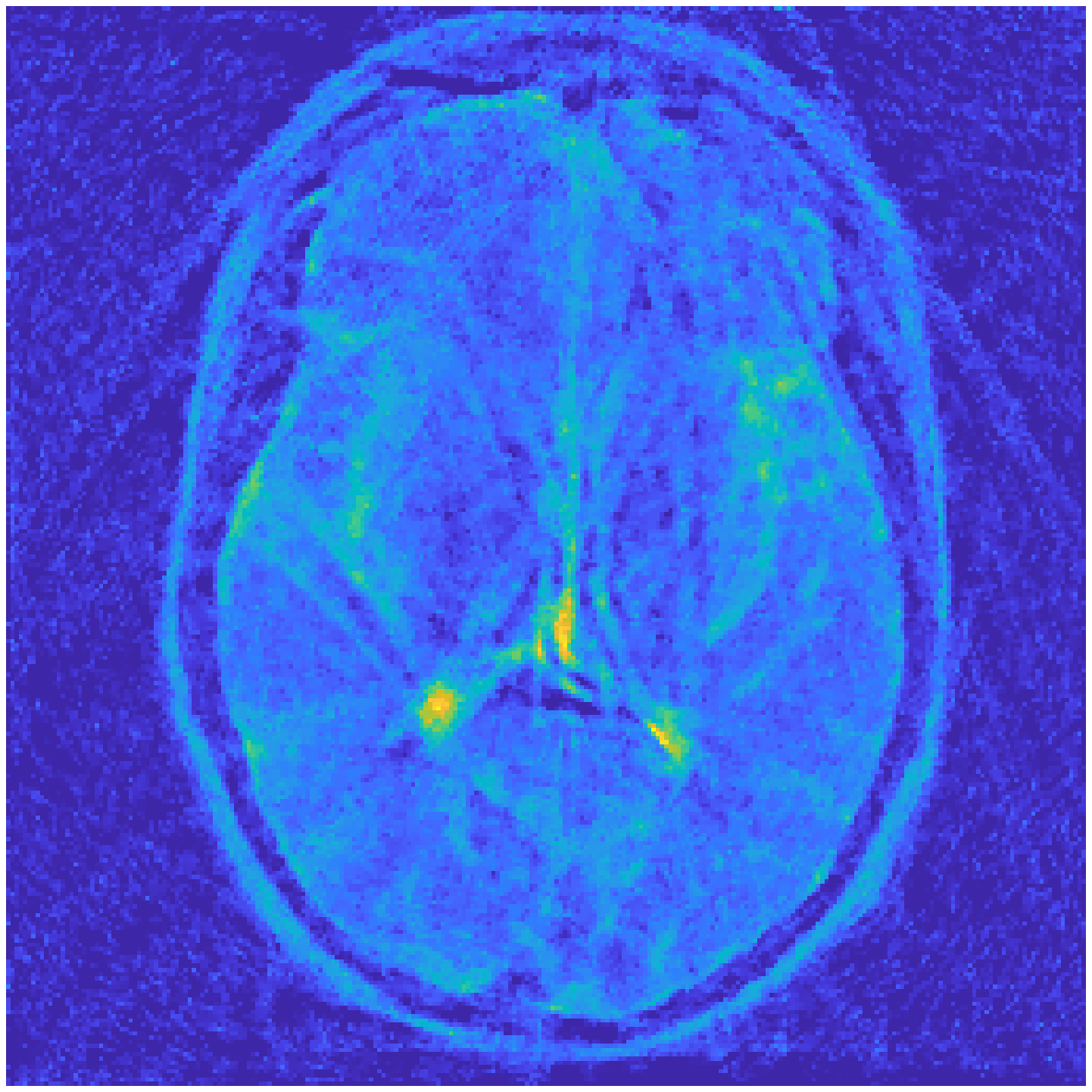} \\

    \includegraphics[width=0.31\linewidth]{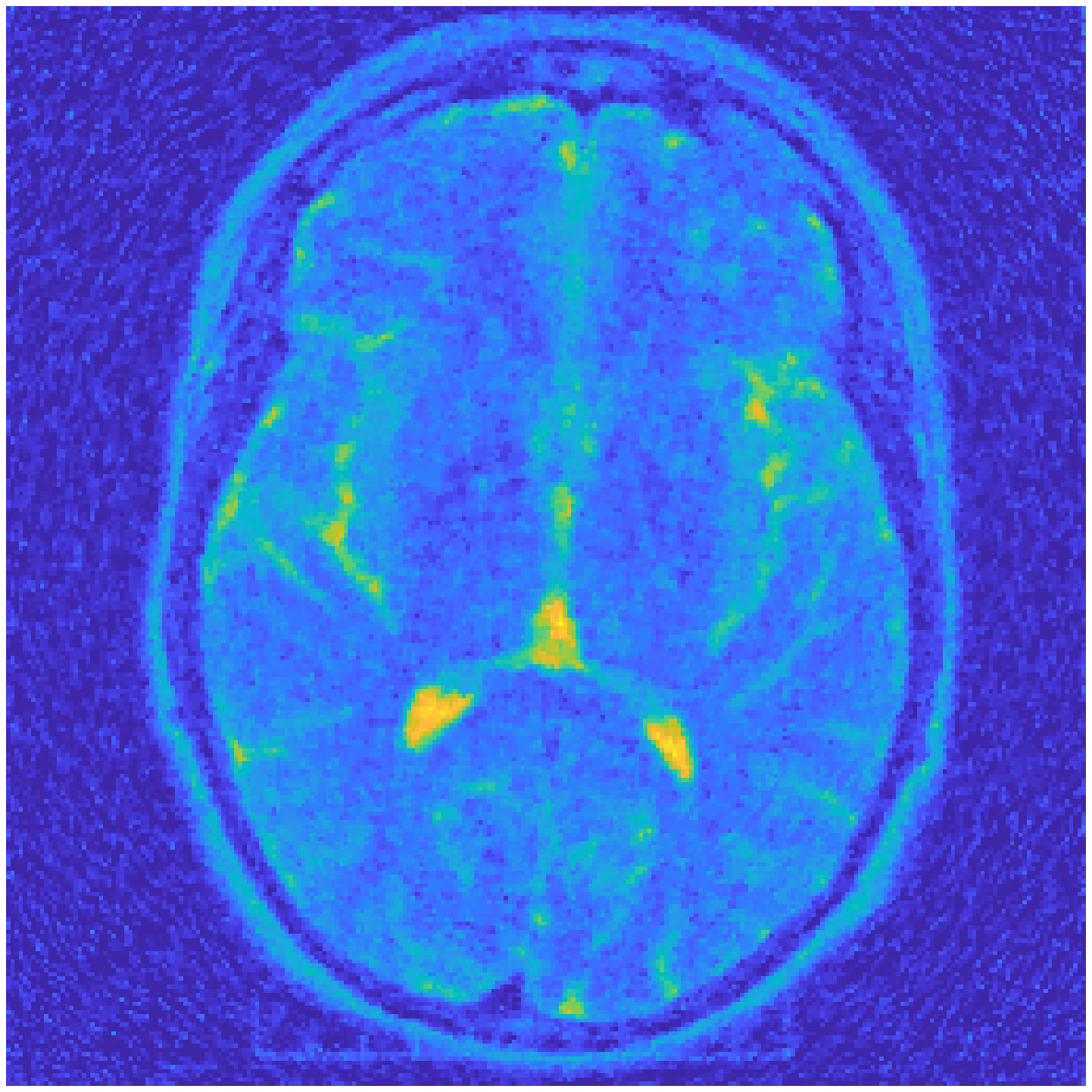} &
    \includegraphics[width=0.31\linewidth]{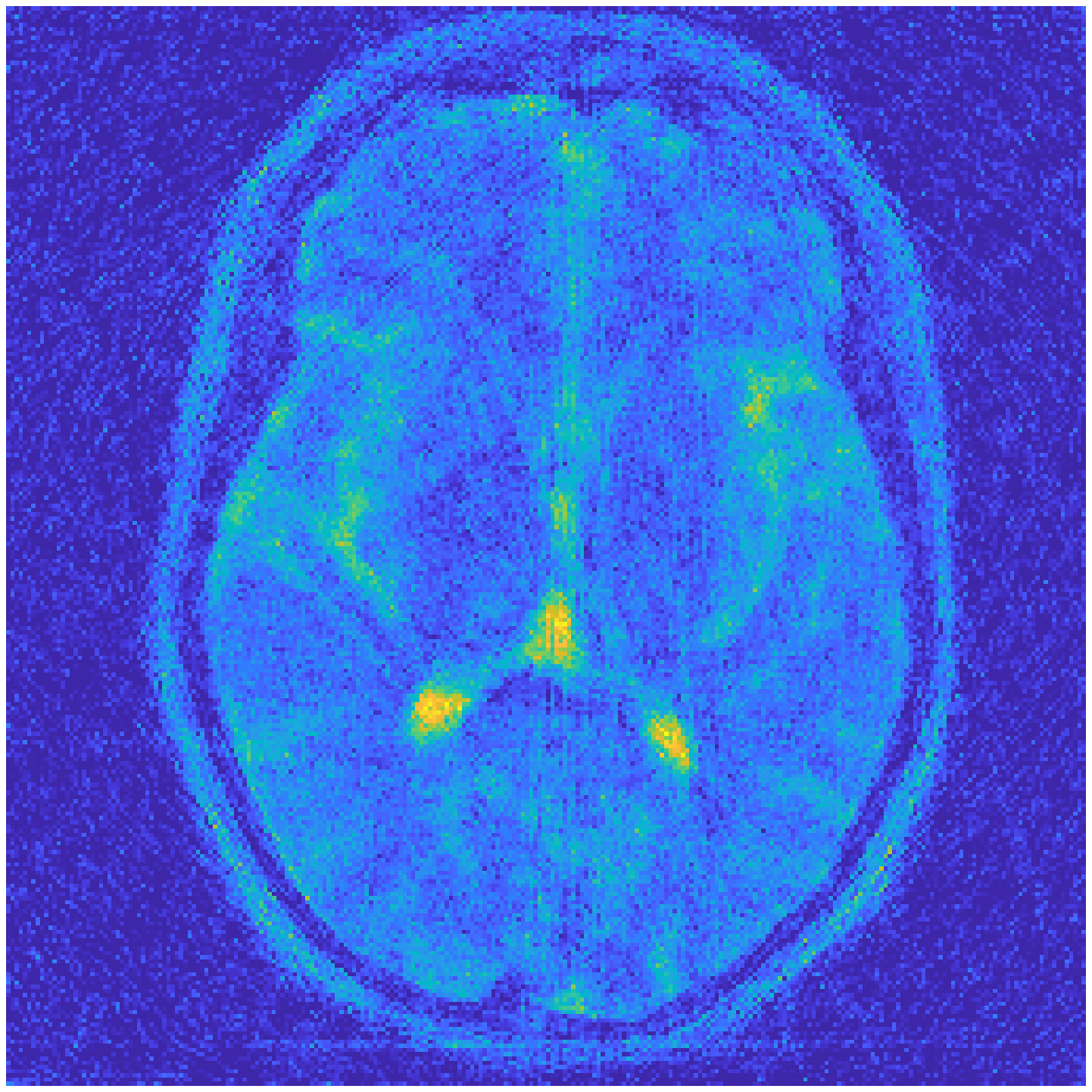} &
    \includegraphics[width=0.31\linewidth]{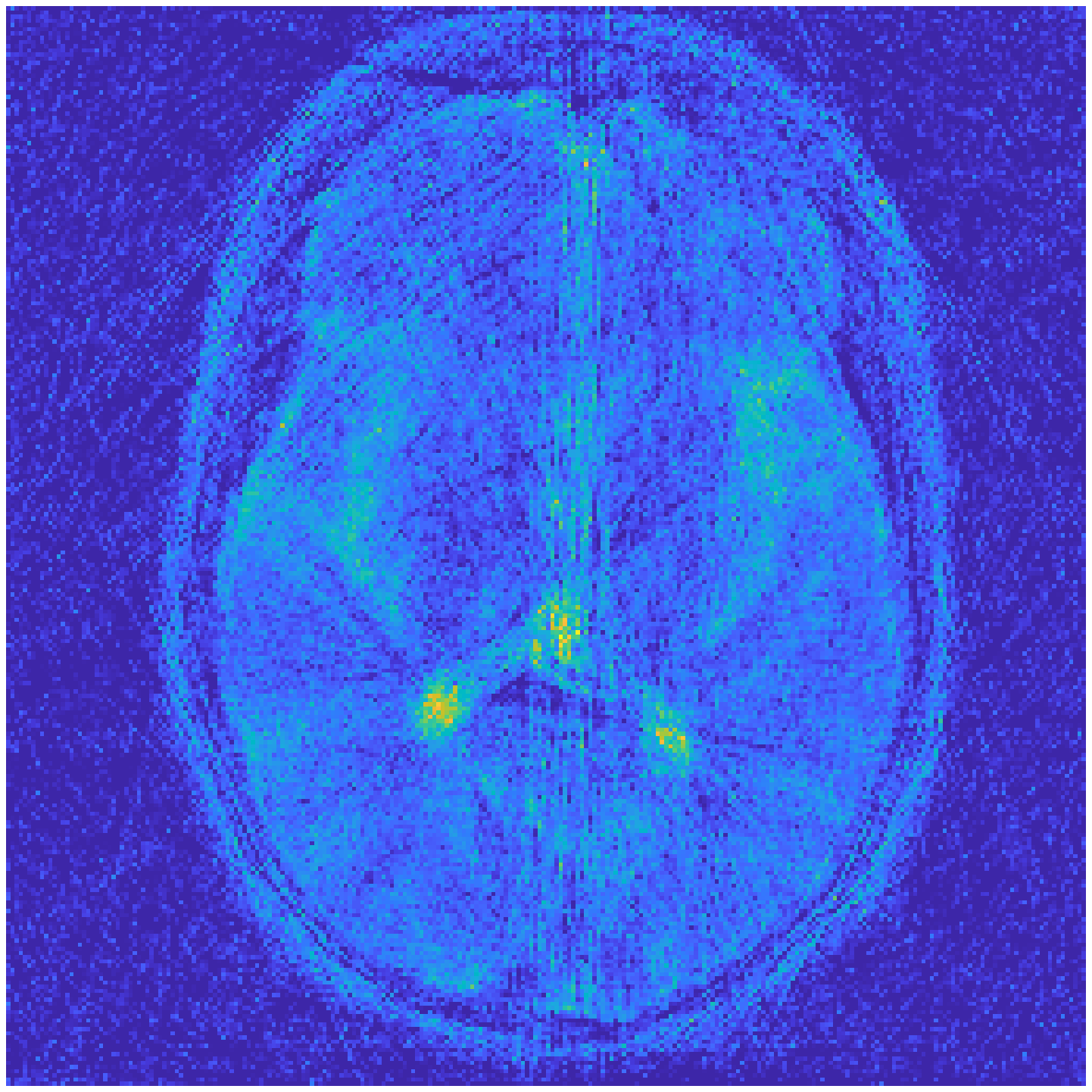} \\

    \includegraphics[width=0.31\linewidth]{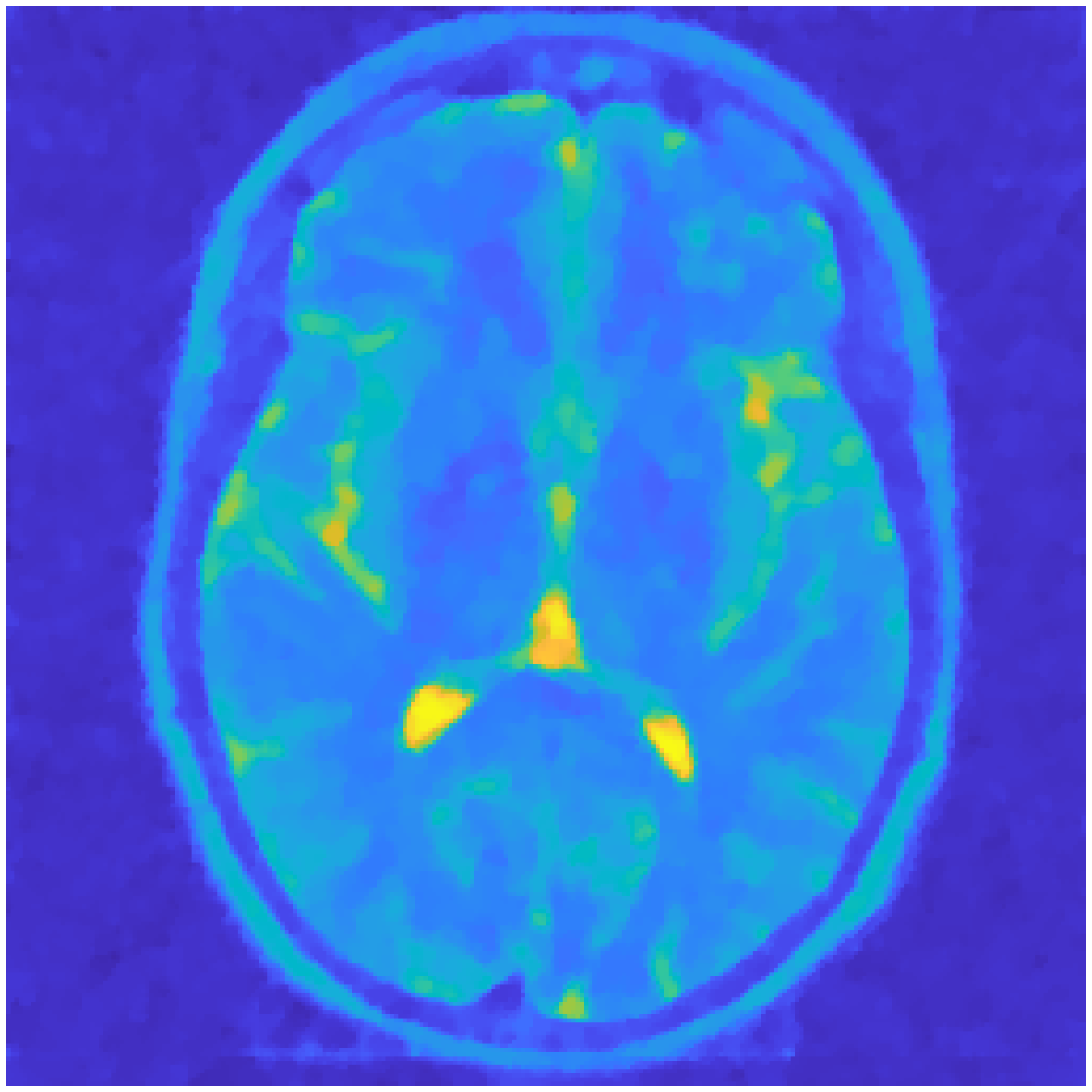} &
    \includegraphics[width=0.31\linewidth]{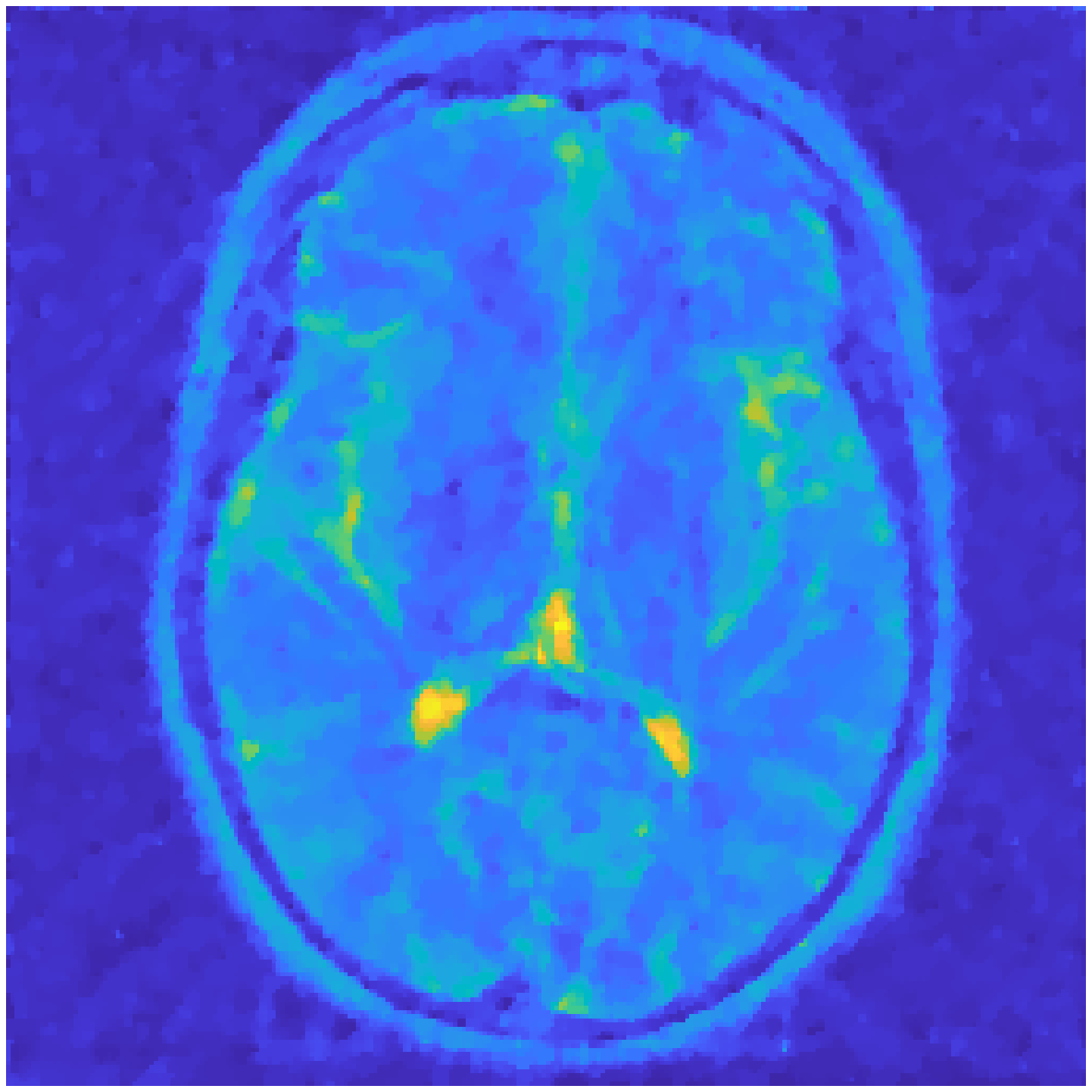} &
    \includegraphics[width=0.31\linewidth]{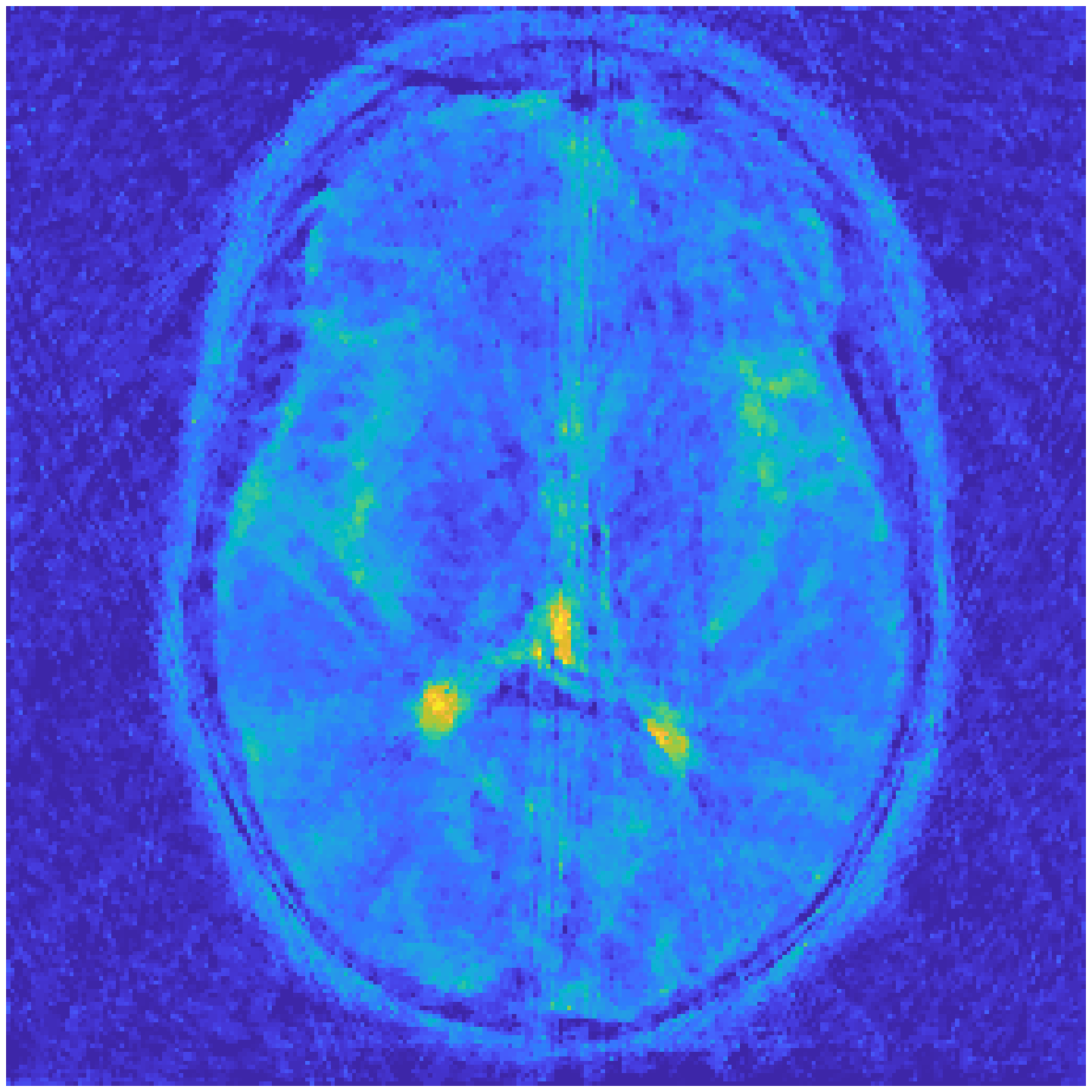} \\

    {\small $\Delta_{\max}=1,\, \sigma=0.01$} &
    {\small $\Delta_{\max}=3,\, \sigma=0.01$} &
    {\small $\Delta_{\max}=5,\, \sigma=0.01$} \\
    \end{tabular}

    \caption{256 by 256 brain. Top row: baseline; middle row: diminishing $\lambda_k$ approach; bottom row: adaptive $\lambda_k$ approach}
    \label{fig:braintype2recs1-001}
\end{figure}

\begin{figure}[!tb]
    \centering
    \setlength{\tabcolsep}{2pt}
    \begin{tabular}{cc}
        \includegraphics[width=0.48\linewidth]{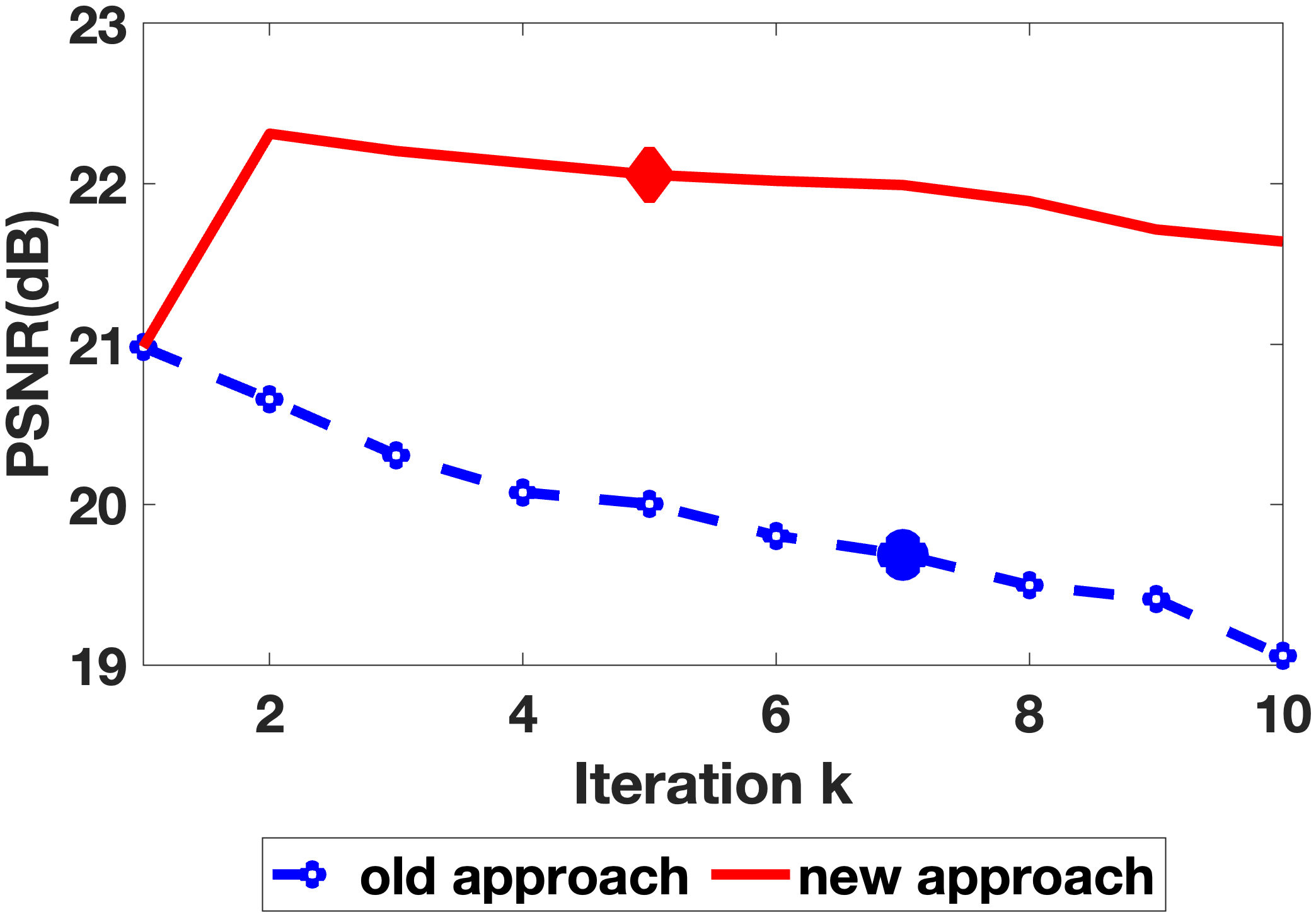} &
        \includegraphics[width=0.48\linewidth]{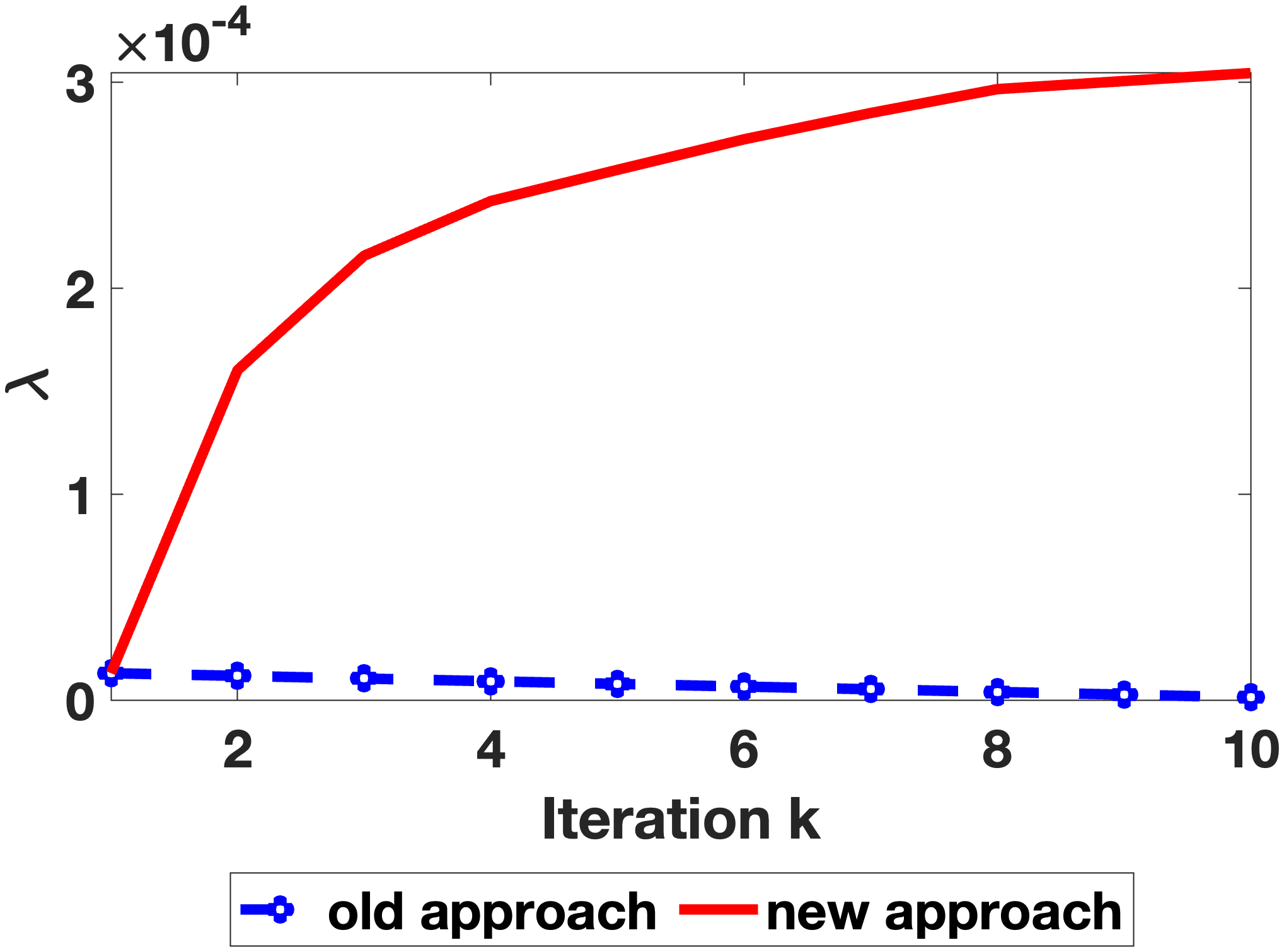} \\

        \includegraphics[width=0.48\linewidth]{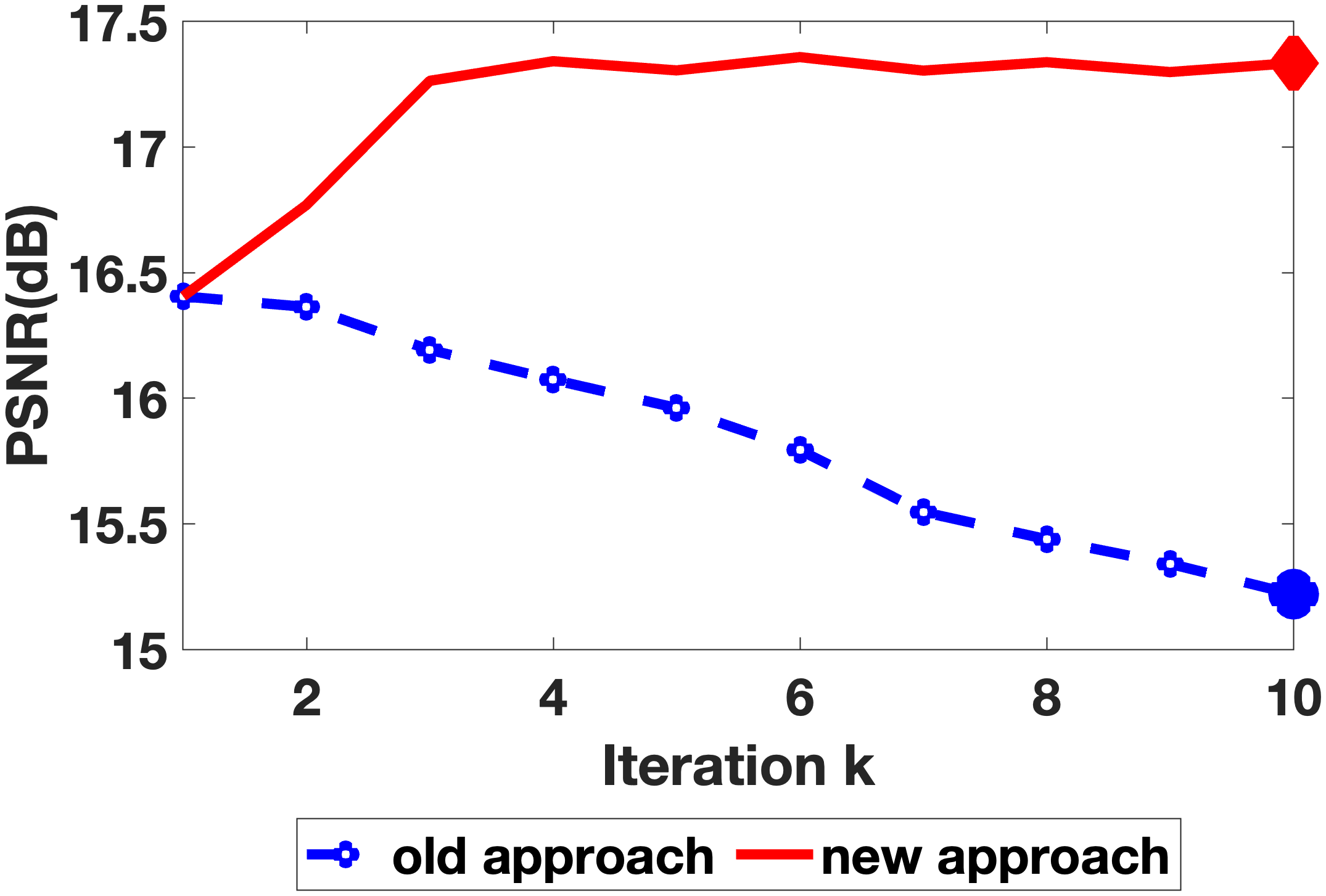} &
        \includegraphics[width=0.48\linewidth]{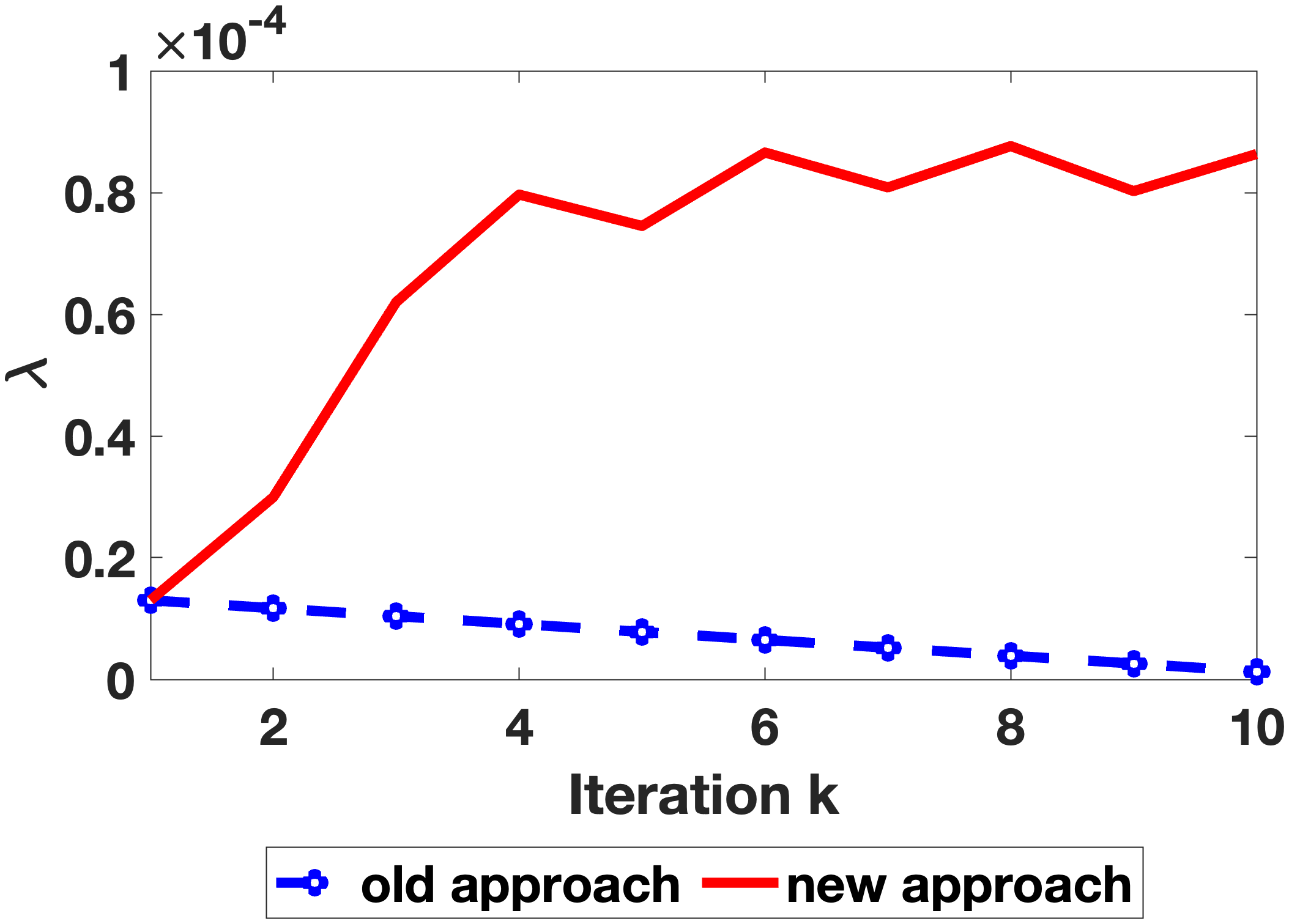} \\

        \includegraphics[width=0.48\linewidth]{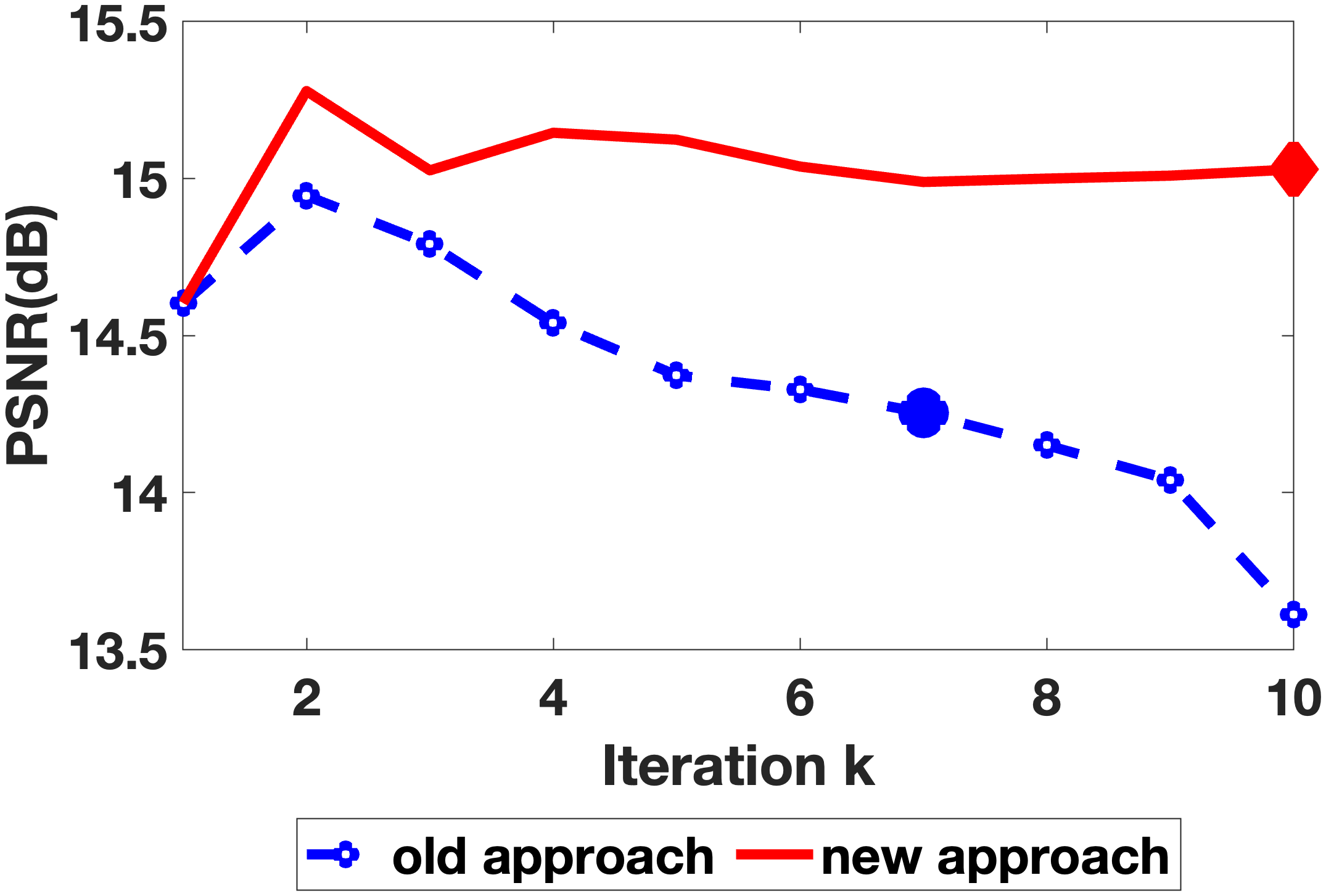} &
        \includegraphics[width=0.48\linewidth]{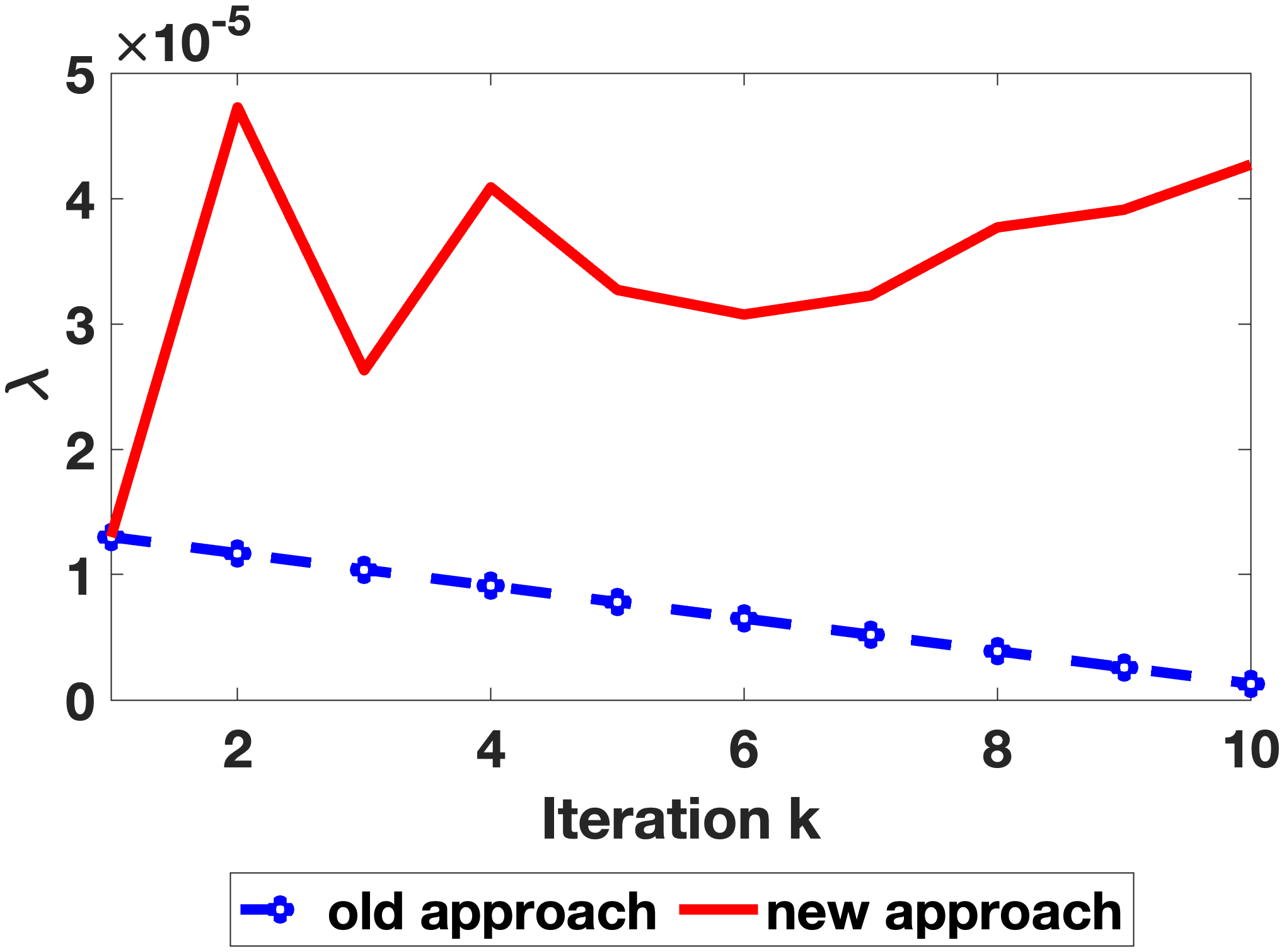} \\
    \end{tabular}
    \caption{256 by 256 phantom. Top row: maxDrift = 1, noise = 0.02; middle row: maxDrift = 3, noise = 0.02; bottom row: maxDrift = 5, noise = 0.02. Left column: PSNR; right column: $\lambda_k$'s}
    \label{fig:phantype2conv1-002}
\end{figure}

\begin{figure}[!tb]
    \centering
    \setlength{\tabcolsep}{2pt}
    \begin{tabular}{cc}
        \includegraphics[width=0.48\linewidth]{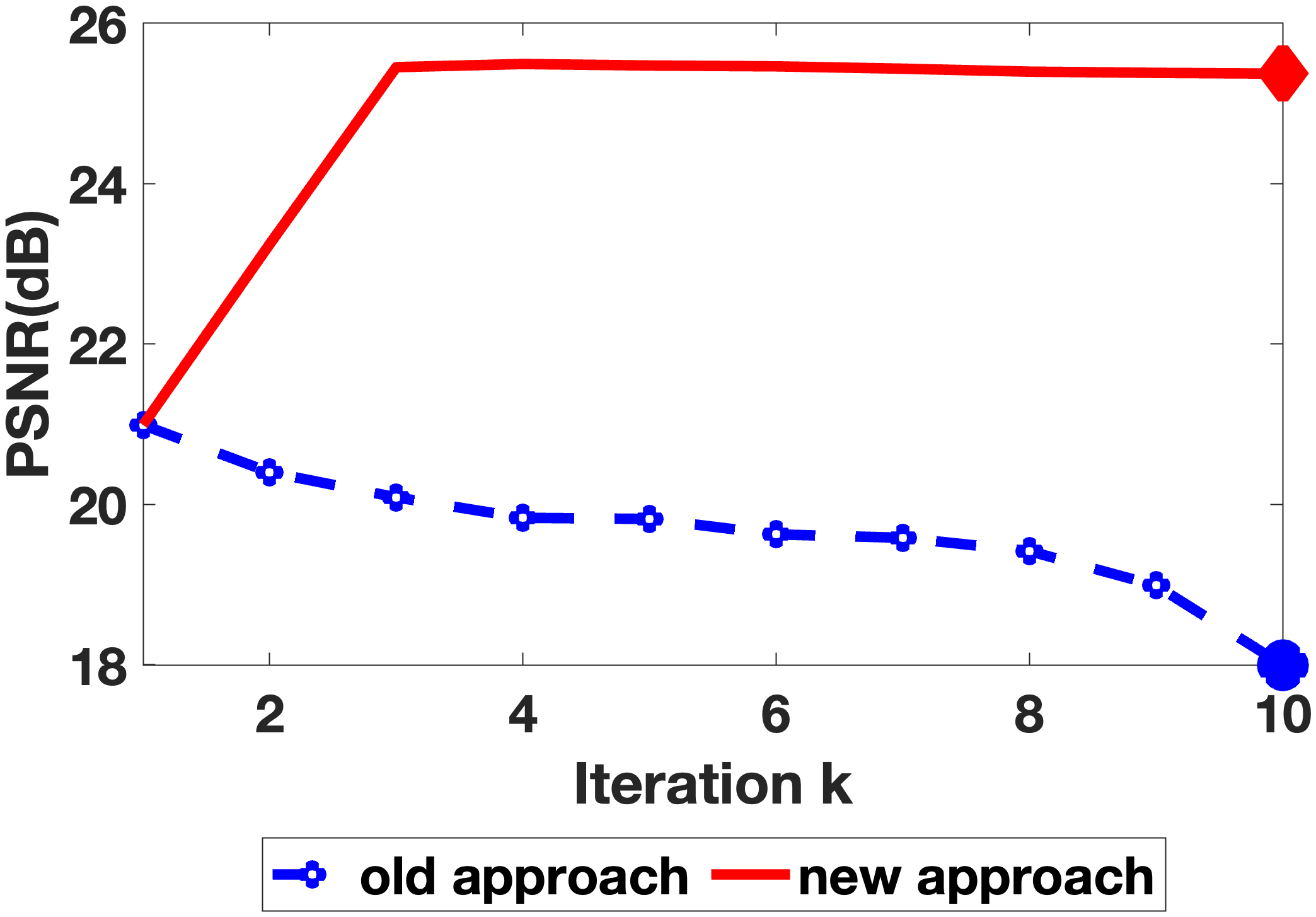} &
        \includegraphics[width=0.48\linewidth]{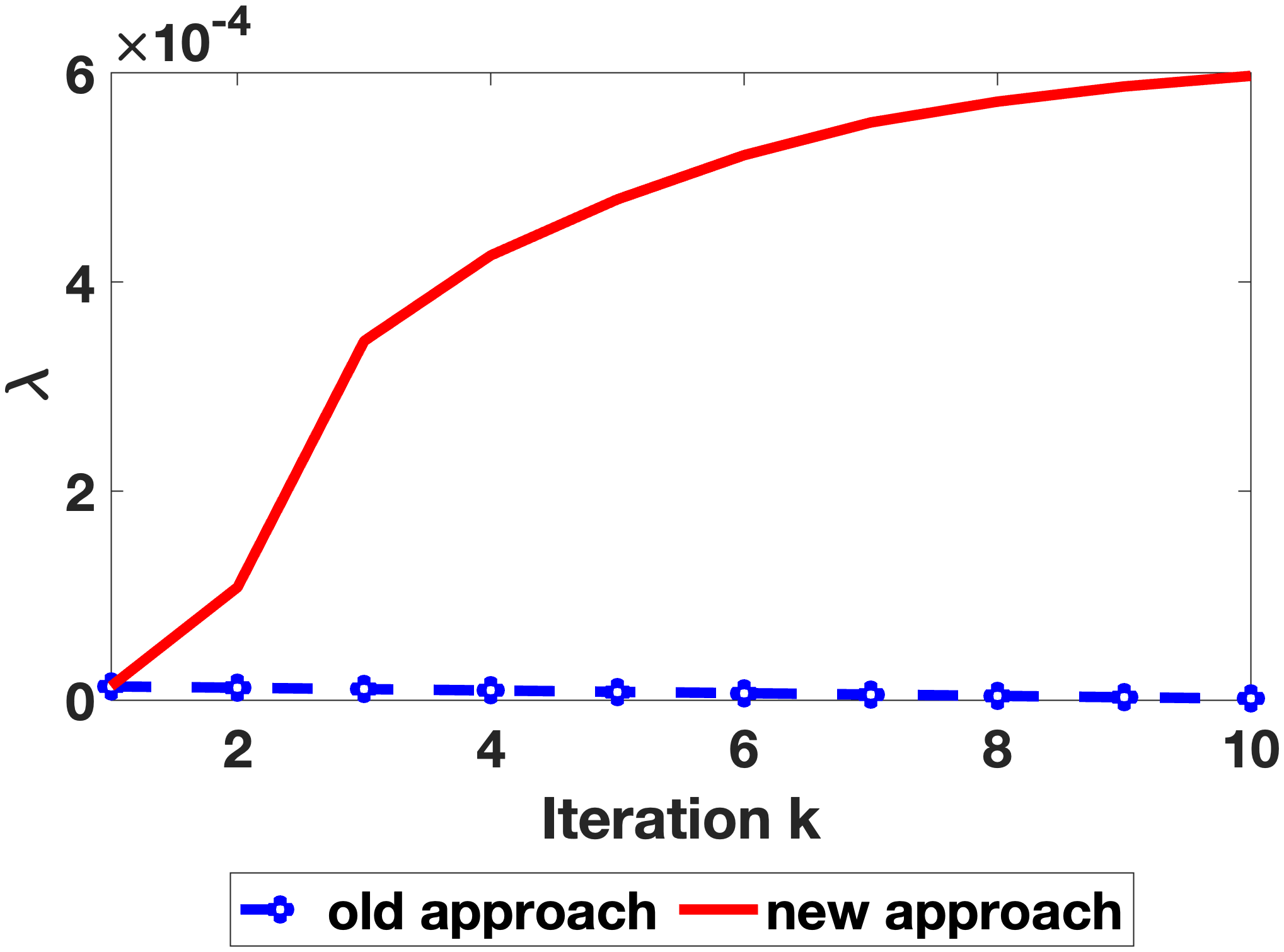} \\

        \includegraphics[width=0.48\linewidth]{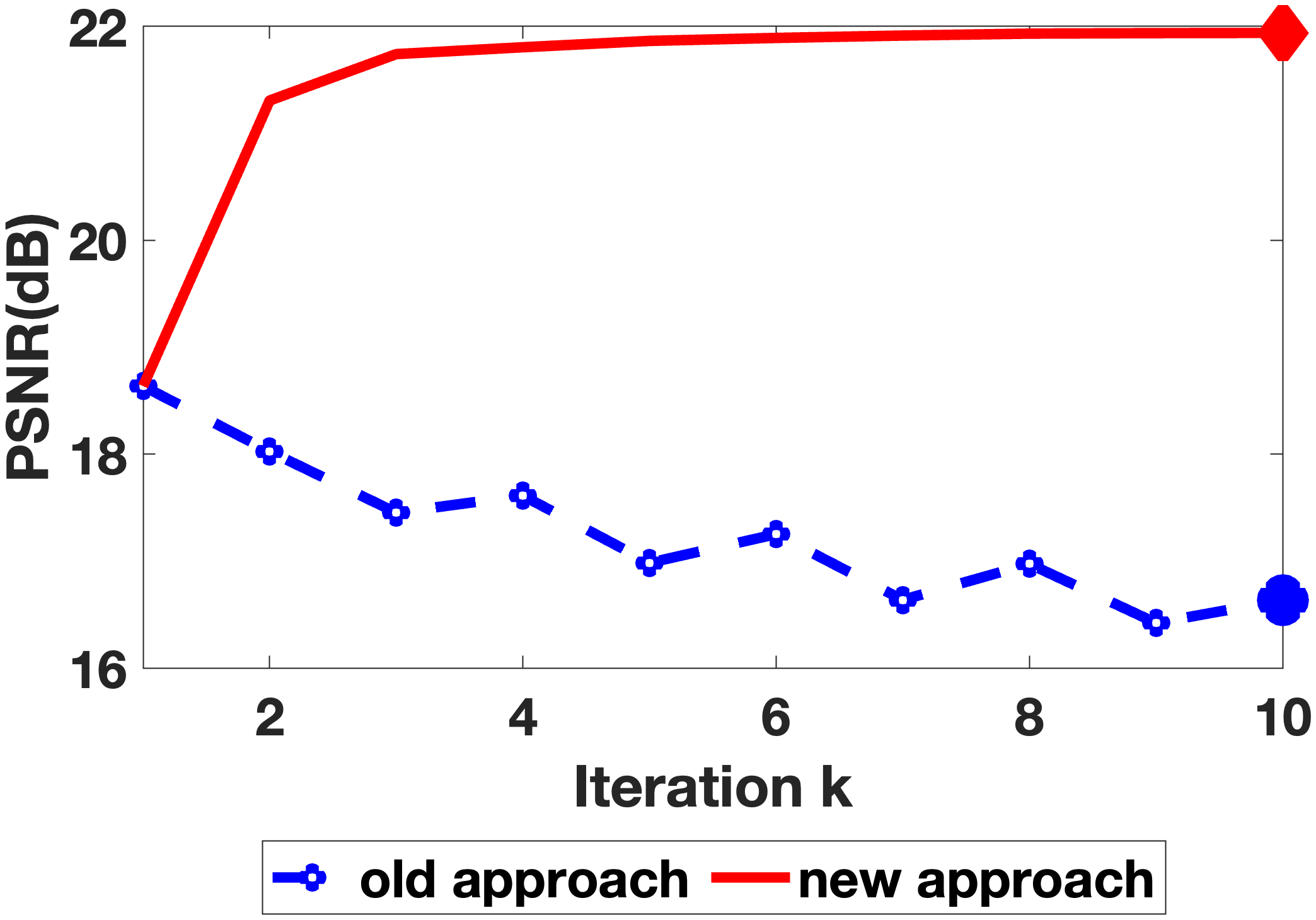} &
        \includegraphics[width=0.48\linewidth]{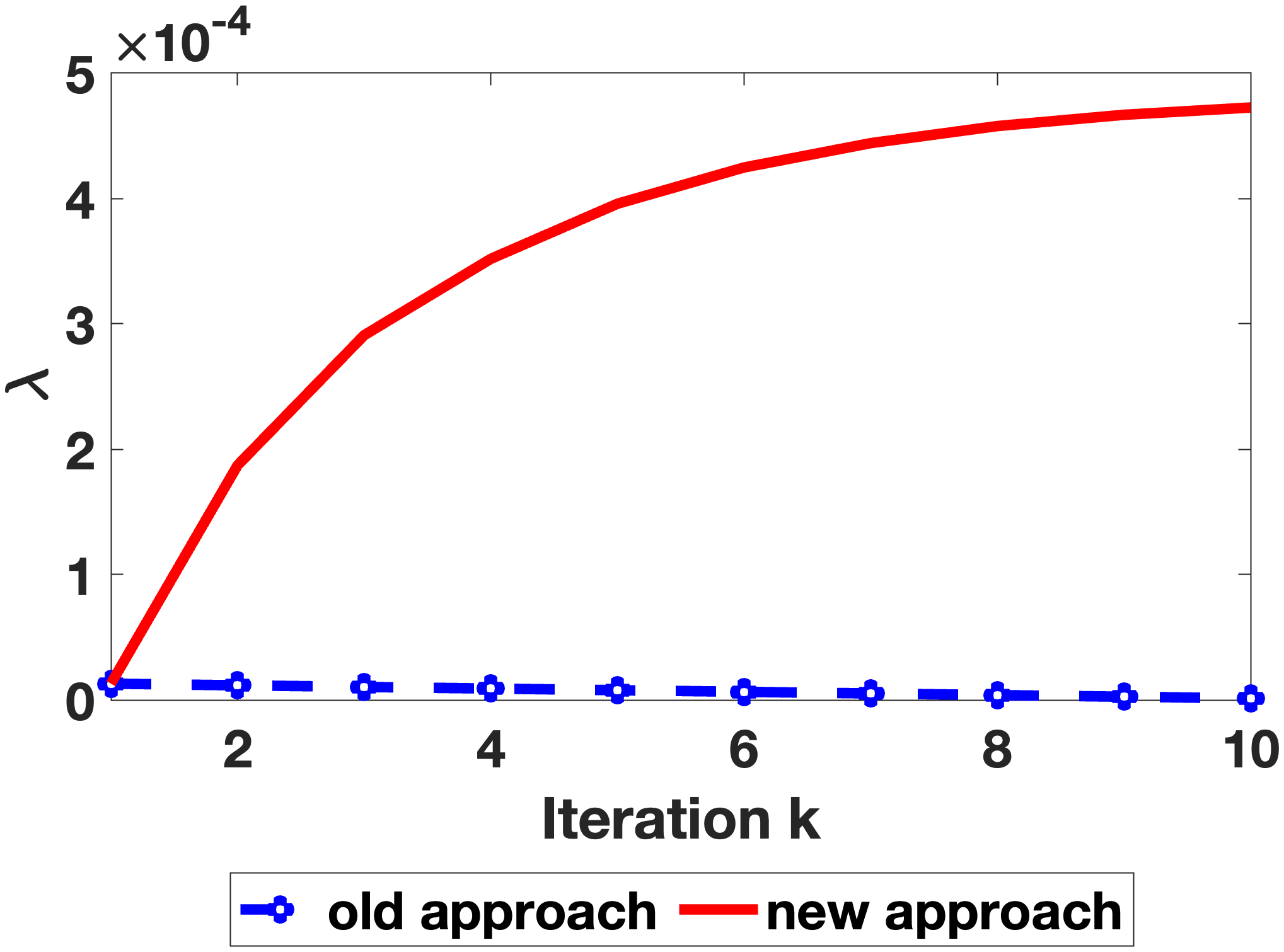} \\

        \includegraphics[width=0.48\linewidth]{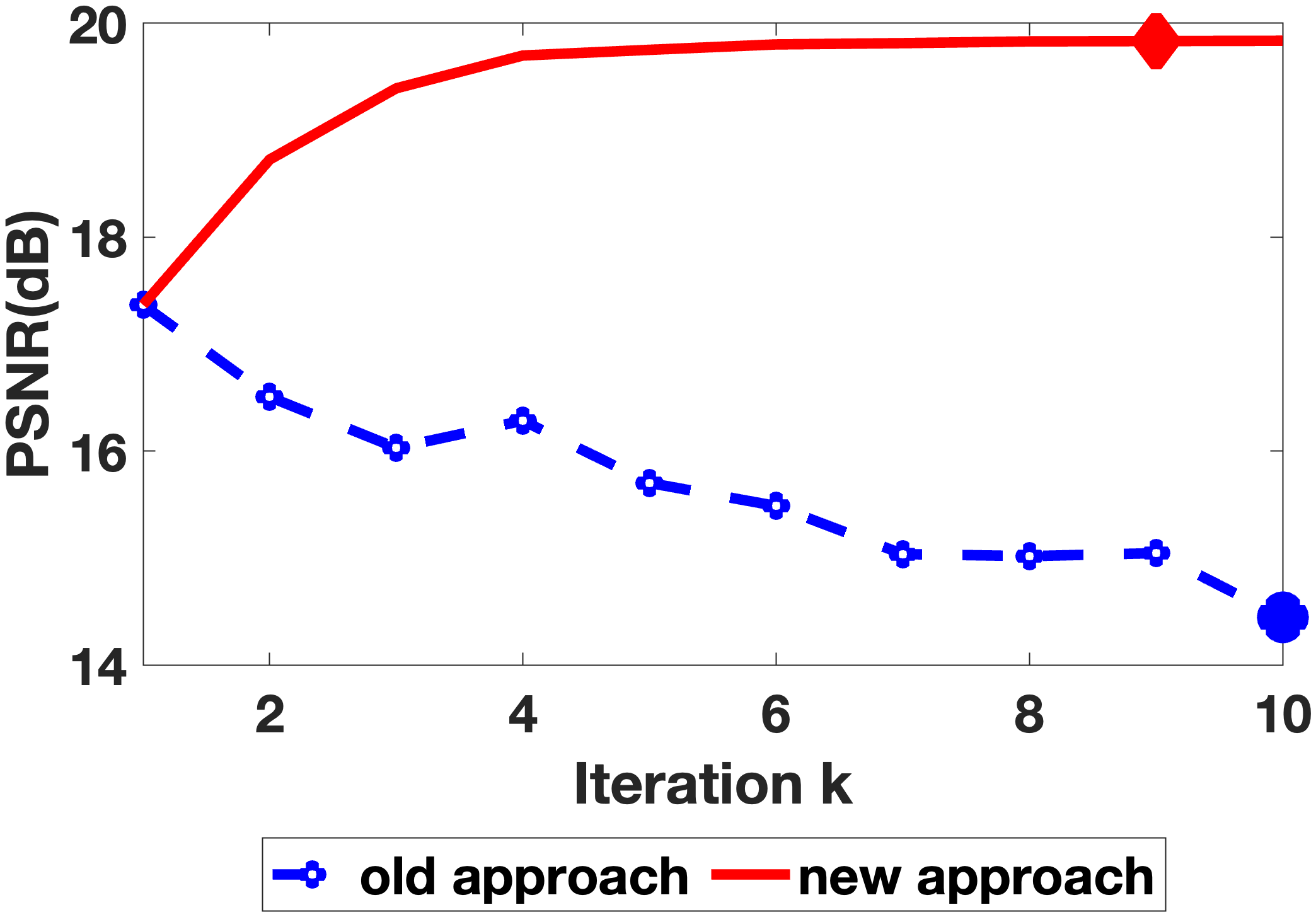} &
        \includegraphics[width=0.48\linewidth]{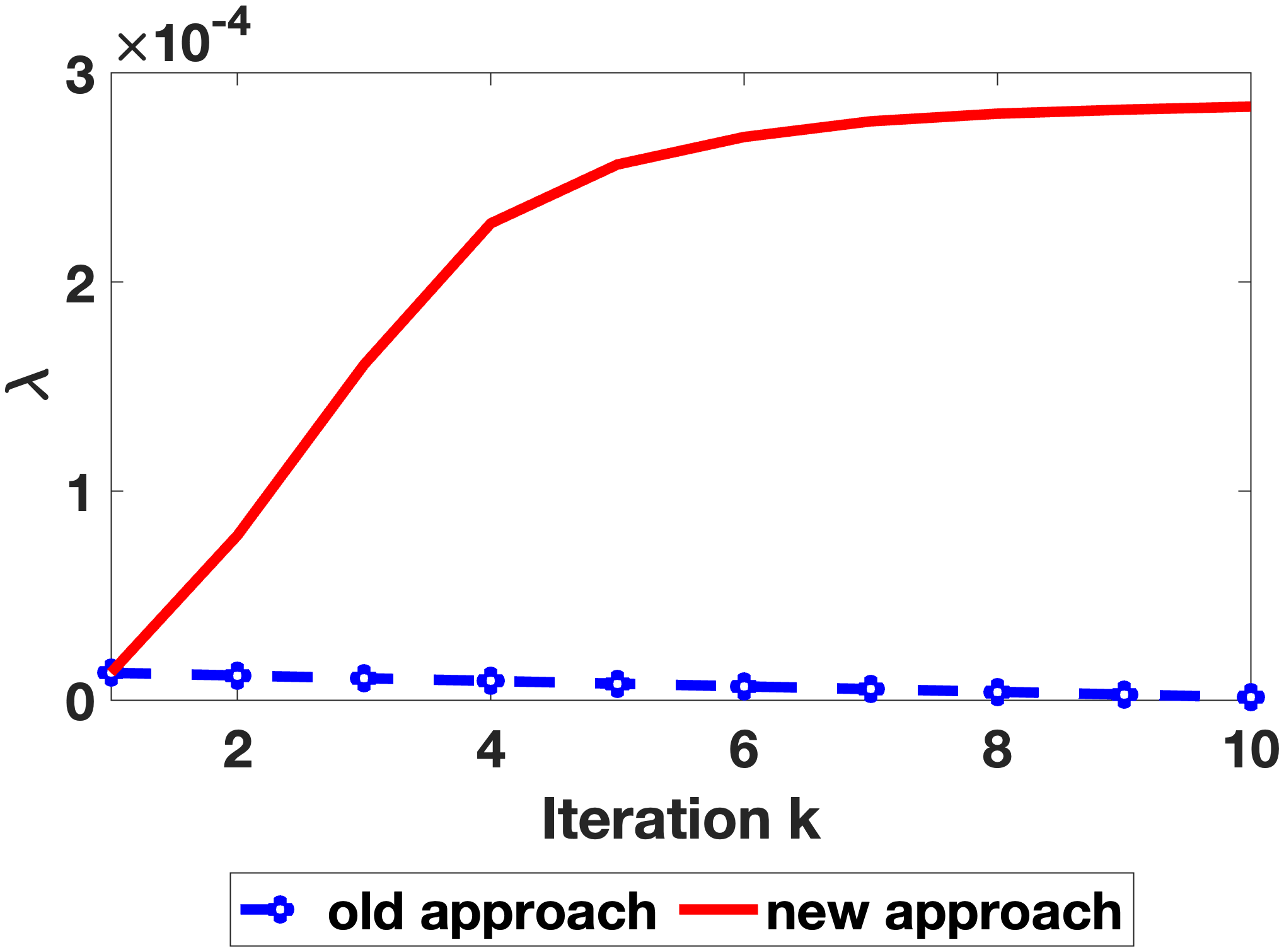} \\
    \end{tabular}
    \caption{256 by 256 brain. Top row: maxDrift = 1, noise = 0.02; middle row: maxDrift = 3, noise = 0.02; bottom row: maxDrift = 5, noise = 0.02. Left column: PSNR; right column: $\lambda_k$'s}
    \label{fig:braintype2conv1-002}
\end{figure}

Visual comparisons of the reconstructed images are provided for the phantom in Fig.~\ref{fig:phantype2recs1-002} and the brain image in Fig.~\ref{fig:braintype2recs1-001}. The baseline reconstructions are visibly blurred and contain artifacts that obscure fine details. Drift calibration mitigates these effects in many settings, although the visual and tabulated comparisons do not always rank the two $\lambda_k$ update rules identically. The iterative performance of the algorithm for the Type II drift case is illustrated in Fig.~\ref{fig:phantype2conv1-002} and Fig.~\ref{fig:braintype2conv1-002}.

For these Type II experiments, the drift regularization parameter $\mu$ is held fixed. Optimizing this parameter, perhaps with an adaptive scheme similar to the one proposed for $\lambda_k$, may further improve the Type II results and is a promising direction for future work.

\begin{table*}[!t]
\caption{Type II drift performance comparison. Blue row: PSNR (dB); gray row: SSIM. }

\label{tab:PSNR_SSIM_Type2}
  \begin{center}
    \begin{tabular}{|c|ccc|ccc|ccc|ccc|}
\hline
Max Drift &   & 1 &   &  & 2 &  &  & 3 &  &  & 5 &  \\
\hline
Noise & 0 & 0.01 & 0.02 & 0 & 0.01 & 0.02 & 0 & 0.01 & 0.02 & 0 & 0.01 & 0.02 \\
\hline \hline
\rowcolor{highlight}Brain (adaptive $\lambda_k$) & 27.50 & 26.37 & 25.37 & 24.12 & 23.77 & 23.52 & 22.32 & 21.89 & 21.94 & 20.68 & 19.13 & 19.83\\
\rowcolor{Gray} & 0.8636 & 0.8009 & 0.7608 & 0.7624 & 0.7307 & 0.7084 & 0.6855 & 0.6311 & 0.6432 & 0.5918 & 0.3167 & 0.4849\\
\hline

\rowcolor{highlight} Brain (diminishing $\lambda_k$) & 27.34 & 24.11 & 18.00 & 23.72 & 21.26 & 17.26 & 20.87 & 18.74 & 16.63 & 16.95 & 16.91 & 14.45\\
\rowcolor{Gray} & 0.8275 & 0.5090 & 0.2354 & 0.6970 & 0.3609 & 0.1888 & 0.4412 & 0.2266 & 0.1520 & 0.1748 & 0.1455 & 0.0806\\

\hline
\rowcolor{white} Brain (baseline) & 27.30 & 23.97 & 19.36 & 23.13 & 21.55 & 18.34 & 20.98 & 19.97 & 17.43 & 18.57 & 17.98 & 16.23\\
\rowcolor{Gray} & 0.8273 & 0.5051 & 0.2799 & 0.6766 & 0.4129 & 0.2249 & 0.5533 & 0.3400 & 0.1799 & 0.3790 & 0.2342 & 0.1197\\
\hline \hline

\rowcolor{highlight} Phantom (adaptive $\lambda_k$) & 22.23 & 22.08 & 22.06 & 19.24 & 18.40 & 19.07 & 17.53 & 16.81 & 17.34 & 16.22 & 15.50 & 15.03\\
\rowcolor{Gray} & 0.9318 & 0.7340 & 0.8557 & 0.8637 & 0.5700 & 0.7442 & 0.8155 & 0.5765 & 0.6185 & 0.7240 & 0.6330 & 0.4664\\
\hline
\rowcolor{highlight} Phantom (diminishing $\lambda_k$) & 22.02 & 21.20 & 19.69 & 17.45 & 17.05 & 17.32 & 16.51 & 15.48 & 15.22 & 14.36 & 13.69 & 14.25\\
\rowcolor{Gray} & 0.8269 & 0.6567 & 0.4884 & 0.6079 & 0.5100 & 0.4428 & 0.5976 & 0.4628 & 0.3857 & 0.4910 & 0.4255 & 0.3835\\
\hline
\rowcolor{white} Phantom (baseline) & 22.06 & 21.64 & 20.51 & 18.50 & 18.31 & 17.74 & 16.61 & 16.48 & 16.09 & 14.75 & 14.67 & 14.39\\
\rowcolor{Gray} & 0.8288 & 0.6841 & 0.5305 & 0.7385 & 0.6226 & 0.4899 & 0.6856 & 0.5813 & 0.4611 & 0.6052 & 0.5198 & 0.4219\\
\hline
\end{tabular}
\end{center}
\end{table*}


\section{Discussion and Conclusion}\label{sec:conclusion}
In this work, we address a pervasive source of experimental error, the scanning-position drift
in scanning-probe x-ray tomography, by means of a
numerical optimization approach. We propose a forward model that embeds the drift directly into the tomographic operator, and we develop an iterative
algorithm~(\cref{alg:driftedTomography}) that calibrates the drifts while reconstructing the object in an alternating fashion. To exploit the structure of
the drift subproblem, we distinguish two regimes: systematic drifts that repeat across projection angles (Type~I, \cref{subsec:SpecialDrift}), for which the
drifted operator reduces to an interpolated row selection of the nominal operator, and fully general per-measurement drifts (Type~II,
\cref{subsec:GeneralDrift}), for which the separability of the Radon transform allows the calibration to be recast as a bound-constrained, per-measurement
least-squares problem. In addition, we introduce an adaptive rule for the reconstruction regularization weight~$\lambda_k$~\eqref{eq:lambdaupdate}, derived
from the discrepancy principle, that automatically adjusts the balance between data fidelity and prior as the forward model is progressively calibrated.

Our numerical study supports three main observations. First, across both drift types and all tested drift magnitudes ($\Delta$ to $5\Delta$) and noise
levels ($\sigma\in\{0,0.01,0.02\}$), explicit drift calibration substantially improves reconstruction quality over the uncalibrated baseline in both PSNR
and SSIM, and the improvement grows as the drift magnitude increases. Second, the proposed adaptive~$\lambda_k$ update attains the best or near-best
accuracy in the majority of cases while removing the need for manual, trial-and-error tuning; the empirical diminishing schedule of~\cite{huang2019calibrating}
remains competitive mainly in some high-noise, high-drift Type~I settings. Third, as shown in the interpolation analysis of~\cref{sec:interp}, linear and
quadratic interpolation of the drifted operator yield nearly identical accuracy, so linear interpolation is preferred for its lower computational cost.

Several limitations point to natural directions for future work. In the Type~II experiments the drift regularization weight~$\mu$ is held fixed, and the
results remain sensitive to its choice; extending the adaptive strategy developed for~$\lambda_k$ to~$\mu$ is a promising avenue for improving the more
challenging general-drift case. The present formulation and experiments are confined to two-dimensional, parallel-beam geometry; the method extends naturally
to three-dimensional parallel-beam reconstruction, which we plan to validate on experimental datasets. Incorporating prior knowledge about the drift---such as
its temporal dynamics or spatial correlation across neighboring measurements---could further constrain the calibration and accelerate convergence. Finally,
the alternating scheme currently comes without theoretical guarantees; a rigorous convergence analysis of the coupled object--drift optimization would help
further stabilize and characterize the performance of the proposed approach.

\section*{Acknowledgments}
This material was based upon work supported by the U.S.\ Department of Energy, Office of Science,
Office of Advanced Scientific Computing Research's applied mathematics program under contract DE-AC02-06CH11357.

\clearpage

\newcommand{\BIBdecl}{\setlength{\itemsep}{0.25 em}}

\bibliographystyle{IEEEtran} 
\bibliography{references}


\end{document}